\documentclass[10pt,reqno]{amsart}
\usepackage{graphicx}
\usepackage{amsmath}
\usepackage{amssymb}
\usepackage{amsthm}

\renewcommand{\theequation}{\thesection.\arabic{equation}}
\newcommand{\sectie}[1]{\setcounter{equation}{0}\section{#1}}
\newcommand{\inlabel}[1]{\refstepcounter{equation} \label{#1}
\hspace{1.5ex} (\theequation)}

\newcommand{\Psis}[3]{\Psi\left(\begin{array}{c} #1
\\ #2 \end{array} ;  \,#3 \right)}
\newcommand{\vfisss}[4]{_2\vfi_1\left(\begin{array}{c} #1 \ \  #2
\\ #3   \end{array} ;  \,#4 \right)}

\newcommand{\R}{\mathbb R}
\newcommand{\C}{\mathbb C}
\newcommand{\Z}{\mathbb Z}
\newcommand{\T}{\mathbb T}
\newcommand{\N}{\mathbb N}
\newcommand{\NN}{\N_0}
\newcommand{\NE}{\N}
\newcommand{\cA}{{\mathcal A}}
\newcommand{\cB}{{\mathcal B}}
\newcommand{\cC}{{\mathcal C}}
\newcommand{\cD}{{\mathcal D}}
\newcommand{\cF}{{\mathcal F}}
\newcommand{\cH}{{\mathcal H}}
\newcommand{\cK}{{\mathcal K}}
\newcommand{\cL}{{\mathcal L}}
\newcommand{\cM}{{\mathcal M}}
\newcommand{\cN}{{\mathcal N}}
\newcommand{\cR}{{\mathcal R}}
\newcommand{\cU}{{\mathcal U}}
\newcommand{\cst}{\text{C}$\hspace{0.1mm}^*$}
\newcommand{\End}{\text{End}}

\newcommand{\al}{\alpha}
\newcommand{\be}{\beta}
\newcommand{\vfi}{\varphi}
\newcommand{\vep}{\varepsilon}
\newcommand{\sde}{\delta}
\newcommand{\de}{\Delta}
\newcommand{\io}{\iota}
\newcommand{\la}{\Lambda}
\newcommand{\si}{\sigma}
\newcommand{\Si}{\Sigma}
\newcommand{\vsi}{\varsigma}
\newcommand{\om}{\omega}
\newcommand{\Om}{\Omega}
\newcommand{\ze}{\zeta}
\newcommand{\flip}{ \chi }
\newcommand{\tde}{\tilde{\de}}

\newcommand{\Mfi}{{\mathcal M}_{\vfi}}
\newcommand{\Nfi}{{\mathcal N}_{\vfi}}
\newcommand{\Mpsi}{{\mathcal M}_{\psi}}

\newcommand{\ot}{\otimes}
\newcommand{\qz}{-q^\Z \cup q^\Z}
\newcommand{\sgn}{\text{\rm sgn}}
\newcommand{\Tr}{\text{\rm Tr}}
\newcommand{\pr}{\T \times I_q}
\newcommand{\spat}{\hspace{4ex}}
\newcommand{\nab}{\nabla}
\newcommand{\hoed}{\hat{\rule{0ex}{1.7ex}}\,}

\newtheorem{definition}{Definition}[section]
\newtheorem{proposition}[definition]{Proposition}
\newtheorem{lemma}[definition]{Lemma}
\newtheorem{corollary}[definition]{Corollary}

\newtheorem{theorem}[definition]{Theorem}

\newtheorem{notation}[definition]{Notation}
\newtheorem{result}[definition]{Result}

\begin{document}

\begin{center}
{\huge\bf A locally compact quantum group analogue of  }

\medskip

{\huge\bf the normalizer of $\mathbf{SU(1,1)}$ in $\mathbf{SL(2,\C)}$}

\bigskip\bigskip\bigskip

{\Large Erik Koelink and Johan Kustermans}

\bigskip\bigskip

\bigskip

Technische Universiteit Delft

Faculteit ITS

Afdeling Toegepaste Wiskundige Analyse

Mekelweg 4

2628CD Delft

The Netherlands

\medskip

\verb"h.t.koelink@its.tudelft.nl" \hspace{3ex} \& \hspace{3ex} \verb"j.kustermans@its.tudelft.nl"

\medskip

\end{center}

\bigskip\bigskip\medskip

\sc Abstract \rm

\medskip

S.L. Woronowicz proved in 1991 that quantum SU(1,1) does not exist as a locally compact quantum group. Results by L.I.
Korogodsky in 1994 and more recently by Woronowicz gave strong indications that the normalizer $\widetilde{SU}(1,1)$ of
$SU(1,1)$ in $SL(2,\C)$ is a much better quantization candidate than $SU(1,1)$ itself. In this paper we show that this
is indeed the case by constructing $\widetilde{SU}_q(1,1)$, a new example  of a unimodular locally compact quantum
group (depending on a parameter $0 < q < 1$) that is a deformation of $\widetilde{SU}(1,1)$. After defining the
underlying von Neumann algebra of $\widetilde{SU}_q(1,1)$ we use a certain class of $q$-hypergeometric functions and
their orthogonality relations to construct the comultiplication. The coassociativity of this comultiplication is the
hardest result to establish. We define the Haar weight and obtain simple formulas for the antipode and its polar
decomposition. As a final result we produce the underlying \cst-algebra of $\widetilde{SU}_q(1,1)$. The proofs of all
these results depend on various properties of $q$-hypergeometric $_1\vfi_1$ functions.

\bigskip\medskip\medskip

\renewcommand{\thefootnote}{}


\newpage

\section*{Introduction}

\bigskip

Arguably one of the most important and simplest non-compact Lie groups is the $SU(1,1)$ group, which is isomorphic to
$SL(2,\R)$. In 1990, one of the first attempts to construct a quantum version of $SU(1,1)$ was made in \cite{Jap1} and
\cite{Jap2} by T. Masuda, K. Mimachi, Y. Nakagami, M. Noumi, Y. Saburi and K. Ueno and independently, in \cite{KorVak}
by L.~Vaksman and L.~Korogodsky. We follow \cite{Jap1} and \cite{Jap2} because these expositions are more elaborate.
Their starting point is a real form $\cU_q(\mathfrak{su}(1,1))$ of the quantum universal enveloping algebra
$\cU_q(\mathfrak{sl}(2,\C))$ (defined in \cite[Eq.~(1.9)]{Jap1}). Intuitively, one should view upon
$\cU_q(\mathfrak{su}(1,1))$ as a quantum universal enveloping algebra of the \lq quantum Lie algebra\rq\ of the still
to be constructed locally compact quantum group $SU_q(1,1)$. The dual $\cA$ of $\cU_q(\mathfrak{su}(1,1))$ is turned
into some sort of topological Hopf$^*$-algebra (\cite[Sec.~2]{Jap1} and \cite[Eq.~(2.6)]{Jap2}) which is naively
referred to as a topological quantum group. In a next step the coordinate Hopf$^*$-algebra $A_q(SU(1,1))$ is given as a
$^*$-subalgebra of $\cA$ that also inherits the comultiplication, counit and antipode from $\cA$
(\cite[Eq.~(2.7)]{Jap1} and \cite[Eq.~(0.9)]{Jap2}).

In this philosophy, they first introduce infinite dimensional infinitesimal representations of quantum $SU_q(1,1)$ (see
\cite[Eq.~(1.1)]{Jap2}) which they then exponentiate to infinite dimensional unitary corepresentations of $\cA$ (see
\cite[Eq.~(1.2)]{Jap2}) providing hereby the quantum analogues of the discrete, continuous and complementary series of
$SU(1,1)$ but also a new strange series of corepresentations. In an attempt to get hold of the \cst-algebra $\cB$ that
should be viewed as the quantum analogue of the space $C_0(SU(1,1))$, they follow the  strategy of first transforming
the generators and relations of $A_q(SU(1,1))$ formally into other elements that produce a $^*$-algebra $B$ (not inside
$A_q(SU(1,1))$ !) definable by generators and relations that, unlike in the case of $A_q(SU(1,1))$, can be faithfully
represented by bounded operators. The \cst-algebra $\cB$ is then defined to be the universal enveloping \cst-algebra of
$B$.

\smallskip

Other important contributions (other than the ones that will be mentioned later on) to the study of different aspects
of quantum $SU(1,1)$ have been made by T.~Kakehi, T.~Masuda, K.~Ueno in \cite{Jap3}, by D. Shklyarov, S. Sinel'shchikov
\& L. Vaksman in \cite{Rus1}.

\smallskip

The authors of \cite{Jap2} seem to get a little bit overenthusiastic by claiming the existence of a comultiplication on
$\cB$ turning $\cB$ into a Hopf$^*$-algebra (see \cite[Prop.~7]{Jap2}). This is in fact in stark contrast with the
result, proven by  S.L. Woronowicz in 1991, that showed that quantum $SU(1,1)$ does not exist as a locally compact
quantum group (see \cite[Thm.~4.1 and Sec.~4.C]{Wor6}). Not surprisingly, this was considered to be quite a setback for
the theory of locally compact quantum groups in the operator algebra setting.

\smallskip

Recall that $SU(1,1)$ is the linear Lie group $\left\{\,X \in SL(2,\C) \mid X^* U X = U \, \right\}$, where $U = \left( \begin{array}{cc} 1 & 0 \\
0 & -1 \end{array}\right)$. In 1994, a breakthrough in this stalemate was forced by L.I. Korogodsky (see \cite{Kor})
who  studied deformations of the linear Lie group $\widetilde{SU}(1,1)$ instead of that of $SU(1,1)$. Here,
$\widetilde{SU}(1,1)$ denotes the  Lie group $\left\{\,X \in SL(2,\C) \mid X^* U X = U \text{ or } X^* U X = - U\,
\right\}$ which is in fact the normalizer of $SU(1,1)$ in $SL(2,\C)$. In this paper, we will follow up on Korogodsky's
and Woronowicz' ideas to introduce a locally compact quantum group $\widetilde{SU}_q(1,1)$, depending on a parameter $0
< q < 1$,  that is the quantum version of this group $\widetilde{SU}(1,1)$. In doing so, we believe we also make a
strong case for the use  of $q$-hypergeometric functions in the operator algebra approach to quantum groups.

\smallskip

But what is actually meant by a locally compact quantum group? This question has kept a lot of people busy for the last
20 years and the most satisfying answer has been given by S. Vaes and the second author in \cite{VaKust} (see the
introduction of this same paper for an extensive account of the history and the importance of different people in the
development of the theory of locally compact quantum groups). The paper \cite{VaKust} is written in the \cst-algebra
setting but it will be easier to use the von Neumann algebraic approach as introduced in \cite{VaKust2}. Both
approaches are completely equivalent and the last result of this paper produces the generic \cst-algebraic quantum
group out of the von Neumann algebraic one. In order to see what we are aiming for, we formulate the definition of a
von Neumann algebraic quantum group as defined in \cite{VaKust2}.

\smallskip\smallskip

\newpage

\bf Definition\! 1.\!\! \it Consider a von Neumann algebra $M$ together with a unital normal $^*$-homomorphism $\de : M
\rightarrow M \ot M$ such that $(\de \ot \io)\de = (\io \ot \de)\de$. Assume moreover the existence of
\begin{enumerate}
\item  a normal semifinite faithful weight $\vfi$ on $M$ that is left invariant:
$\vfi((\om \ot \io) \de(x)) = \vfi(x) \om(1)$ for all $\om \in M_*^+$ and $x \in \Mfi^+$.
\item  a normal semifinite faithful weight $\psi$ on $M$ that is right invariant:
$\psi((\io \ot \om) \de(x)) = \psi(x) \om(1)$ for all $\om \in M_*^+$ and $x \in \Mpsi^+$.
\end{enumerate}
Then we call the pair $(M,\de)$ a von Neumann algebraic quantum group. \rm

\smallskip\smallskip

For a discussion about the consequences of this definition and related notations we refer to \cite[Sec.~1]{VaKust2}.
For quite a while, the major drawback of the theory of non-compact quantum groups was the lack of a whole array of
different examples of non-compact locally compact quantum groups. So the importance of this paper lies mainly in the
fact that we add an example to the rather short list of existing atomic non-compact quantum groups like quantum E(2),
quantum $ax+b$ and quantum $az+b$. Also note that this new locally compact quantum group is the first analogue of a
non-compact semi-simple Lie group.  By atomic we mean that these examples (up till now) can not be constructed out of
simpler quantum groups through different existing theoretical construction procedures of which the most important one
is arguably the double crossed product construction. In turn, this new example $\widetilde{SU}_q(1,1)$ can now serve as
one of the ingredients in these construction procedures. It is moreover very conceivable  that the importance of
$SU(1,1)$ in the classical group theory  will be parallelled by that of $\widetilde{SU}_q(1,1)$ in the quantum group
setting.

\smallskip

In the immediate future we hope to calculate the dual of $\widetilde{SU}_q(1,1)$, unravel its corepresentation theory
and obtain a formula for the Plancherel measure. For one, we most definitely plan to fit  the corepresentations of
\cite{Jap2} rigorously in the framework of locally compact quantum groups. It is to be expected that new results for
$q$-hypergeometric functions can emerge this way.

\smallskip

Since the ideas of Korogodsky are a chief motivation for this paper, we will briefly discuss the most important results
of \cite{Kor}. Let $(\cA_q,\de)$ be the Hopf$^*$-algebra  associated to $\widetilde{SU}_q(1,1)$ (in \cite{Kor}, $\cA_q$
is denoted by ${\mathfrak S}$). We give a  precise definition of $(\cA_q,\de)$ in Eqs.~(\ref{eqset1}) and
(\ref{hopf.comult1}). The $^*$-algebra $\cA_q$ is generated by elements $\al$, $\gamma$ and a central self-adjoint
involution $e$. Define two central orthogonal projections $p_{\pm}$ in $\cA_q$ as $p_{\pm} = \frac{1}{2}(1 \pm e)$ and
define the $^*$-subalgebras $\cA_q^{\pm}$ of $\cA_q$ as $\cA_q^{\pm} = p_{\pm} \cA_q$. Thus, $\cA_q  = \cA^+_q \oplus
\cA^-_q$ (in \cite{Kor}, $\cA_q^\pm$ are denoted by ${\mathfrak R}^\pm$). Let us also set $\gamma_{\pm} = p_{\pm}
\gamma$.

Define a $^*$-homomorphism $\de_+ : \cA_q^+ \rightarrow \cA_q^+ \ot \cA_q^+$  such that $\de_+(a) = (p_+ \ot
p_+)\de(a)$ for all $a \in \cA_q^+$. Then $(\cA_q^+,\de_+)$ is a Hopf$^*$-algebra which should be thought of as the
$q$-deformation of the coordinate Hopf$^*$-algebra associated to the group $SU(1,1)$.

Most of the relevant representations of $\cA_q$ are infinite dimensional so that a lot of care has to be taken to make
the notion of representations of $\cA_q$ more precise (see \cite{Schm}). Consider a Hilbert space $H$, a dense subspace
$D$ of $H$ and a unital algebra representation $\pi$ of $\cA_q$ on $D$ such that $\langle \pi(a) v , w \rangle =
\langle v , \pi(a^*) w \rangle$ for all $a \in \cA_q$ and $v,w \in D$. Then we call $\pi$ a $^*$-representation of
$\cA_q$ {\it in} $H$. We set $\cD(\pi) = D$.  We call $\pi$ well-behaved if  $\pi(\gamma^* \gamma)$ is an essentially
self-adjoint operator in $H$ and if the spectrum $\si\bigl(\,\overline{\pi(\gamma^* \gamma)}\,\bigr) \subseteq
q^{2\Z}$, where $\overline{\pi(\gamma^* \gamma)}$ denotes the closure of $\pi(\gamma^* \gamma)$ as an operator in $H$.
Similar definitions are used for $^*$-representations of $\cA_q^\pm$ in $H$ (where $\gamma$ is replaced by
$\gamma_{\pm}$). The spectral condition is not present in \cite{Kor}, a weaker spectral condition was introduced in
\cite{Wor7}.

It is clear that we can identify $^*$-representations of $\cA^\pm_q$  with $^*$-representations $\pi$ of $\cA_q$ such
that $\pi(p_{\mp}) = 0$, or equivalently $\pi(e) = \pm 1$.

Since $e$ belongs to the center of $\cA_q$, each irreducible $^*$-representation of $\cA_q^{\pm}$ is necessarily an
irreducible $^*$-representation of either $\cA_q^+$ or $\cA_q^-$. Korogodsky classified all well-behaved
$^*$-representations of $\cA_q$ in \cite[Prop.~2.4]{Kor}, most of them are infinite dimensional.  K. Schm\"{u}dgen
pointed us to the necessity of imposing some sort of spectral condition for this classification result to be true (and
which is completed in \cite{Wor7}).

Given two well-behaved infinite dimensional irreducible $^*$-representations $\pi_1$, $\pi_2$ of $\cA_q$ in Hilbert
spaces $H_1$, $H_2$ respectively. Korogodsky showed in \cite[Thm.~6.1]{Kor} that there does not exist a well-behaved
$^*$-representation $\pi$ of $\cA_q$ in $H_1 \ot H_2$ such that $\cD(\pi_1) \odot \cD(\pi_2) \subseteq \cD(\pi)$ and
$\pi(a)\, v = (\pi_1 \odot \pi_2)(\de(a))\,v$ for all $a \in \cA_q$ and $v \in \cD(\pi_1) \odot \cD(\pi_2)$, where
$\odot$ denotes the algebraic tensor product. Thus, loosely speaking, the tensor product of two well-behaved infinite
dimensional irreducible $^*$-representations $\pi_1$ and $\pi_2$ of $\cA_q$ does not exist. Applied to two well-behaved
infinite dimensional $^*$-representations of $\cA_q^+$, this result corresponds to the non-existence of quantum
$SU(1,1)$ as a locally compact quantum group as proven in \cite{Wor6}.

However, Korogodsky also showed that the situation is not as bad as the above result on first sight suggests. Therefore
consider two well-behaved $^*$-representations $\pi_1$, $\pi_2$ of $\cA_q$ in Hilbert spaces $H_1$, $H_2$ such that for
$i=1,2$, $\pi_i$ is the finite direct sum of well-behaved infinite dimensional irreducible $^*$-representations of
$\cA_q$ for which the number of direct summands of $\pi_i$ that are $^*$-representations of $\cA_q^+$ equals the number
of direct summands of $\pi_i$ that are $^*$-representations of $\cA_q^-$. Then \cite[Thm.~6.2 and 6.3]{Kor} imply the
existence of a well-behaved $^*$-representation $\pi$ of $\cA_q$ in $H_1 \ot H_2$ such that $\cD(\pi_1) \odot
\cD(\pi_2) \subseteq \cD(\pi)$ and $\pi(a)\, v = (\pi_1 \odot \pi_2)(\de(a))\,v$ for all $a \in \cA_q$ and $v \in
\cD(\pi_1) \odot \cD(\pi_2)$. It should however be pointed out that $\pi$ is not unique!

\smallskip

Around 1996, Woronowicz picked up on this observation to start his study of quantum $\widetilde{SU}_q(1,1)$ on the
Hilbert space level (see \cite{Wor7}). Instead of focusing on all well-behaved  $^*$-representation of $\cA_q$, he
redirected his attention to $^*$-representations that are direct sums (or more precisely, direct integrals) of the sort
described above. For this purpose, Woronowicz introduces in   \cite[Def.~1.1]{Wor7} the notion of a
$\widetilde{SU}_q(1,1)$-quadruples on a Hilbert space. Let us use some slightly different terminology in order to stay
closer to Korogodsky's viewpoint.

Let $\pi$ be a well-behaved $^*$-representation of $\cA_q$ in some Hilbert space $H$ and $Y$ a closed linear operator
in $H$ such that $\text{Phase}\ \overline{\pi(\gamma)}$ commutes with $Y$, the domain of $Y$ and $Y^*$ agree and
$\pi\bigl(\sqrt{q}\,e\,\gamma^* - \al\bigr) \subseteq Y$. Then we call $(\pi,Y)$ a
$\widetilde{SU}_q(1,1)$-representation in $H$. It follows that for such a $\widetilde{SU}_q(1,1)$-representation, the
quadruple $(\overline{\pi(\al)},\overline{\pi(\gamma)},\overline{\pi(e)},Y)$ forms a $\widetilde{SU}_q(1,1)$ of
unbounded type in the sense of \cite{Wor7}.

If $(\tilde{\al},\tilde{\gamma},\tilde{e},\tilde{Y})$ is a $\widetilde{SU}_q(1,1)$-quadruple of unbounded type, it
follows from \cite[Thm.~3.7]{Wor7} that it arises from such a $\widetilde{SU}_q(1,1)$-representation. Moreover,
\cite[Thm.~1.3]{Wor7} gives a precise meaning to the statement that $\pi$ contains as much irreducible well-behaved
$^*$-representations of $\cA_q^+$ as of $\cA_q^-$. In  \cite[Thm.~1.4]{Wor7} it is then shown that, in accordance with
the above formal discussion, the tensor product of two $\widetilde{SU}_q(1,1)$-quadruples can be defined as a new
$\widetilde{SU}_q(1,1)$-quadruple. Recasted in the terminology of this introduction this means that given two
$\widetilde{SU}_q(1,1)$-representations $(H_1,\pi_1,Y_1)$, $(H_2,\pi_2,Y_2)$ it is possible to construct a new
$\widetilde{SU}_q(1,1)$-representation $(H_1 \ot H_2,\pi,Y)$ such that $\cD(\pi_1) \odot \cD(\pi_2) \subseteq \cD(\pi)$
and $\pi(a)\, v = (\pi_1 \odot \pi_2)(\de(a))\,v$ for all $a \in \cA_q$ and $v \in \cD(\pi_1) \odot \cD(\pi_2)$. The
extra information contained in $Y$ moreover allows to get rid of the uniqueness problems mentioned before.

There is however one drawback to \cite{Wor7}. At the moment of writing this paper, there was still a gap in the proof
of the associativity of the tensor product construction $\widetilde{SU}_q(1,1)$-quadruples. Although it is difficult to
estimate the seriousness of this gap, it corresponds to the most difficult part of the proof of the coassociativity of
the  comultiplication in this paper.

\smallskip

In 1999, J. Stokman together with the first author studied the properties of a left invariant weight $h$ on the
\cst-algebra $C$ of quantum $SU(1,1)$ (see \cite{ErikJasper2}). However, remember that it is impossible to define a
comultiplication on $C$ turning it into a locally compact quantum group. Instead, one looks at a suitable dense
subspace $B$ of $C^*$ so that for $\om_1$, $\om_2$ it is possible to define  the product $\om_1 \star \om_2$ so that on
a formal level, $\om_1 \star \om_2 = (\om_1 \ot \om_2)\de$. It turns out that the invariant weight $h$ can be written
as a sum of positive functionals in $B$ and as a consequence $(\om \star h)(x)$ can be defined as a pointwise
convergent sum for suitable $\om \in B$ and suitable $x \in C$. The left invariance statement then becomes $(\om \star
h)(x) = \om(1)\,h(x)$. The benefit of \cite{ErikJasper2} to this paper lies in the fact that the definition of $\om_1
\star \om_2$ depends on a class of special functions  that arise formally as Clebsch-Gordan coefficients. These special
functions can be extended in such a way that they serve as the principal set of data to construct the locally compact
quantum group $\widetilde{SU}_q(1,1)$. Moreover, the formula for the left Haar weight of $\widetilde{SU}_q(1,1)$ is a
straightforward generalization of the formula defining $h$.

\smallskip

A naive response to this all would be to presume that given all these existing results it should be not so difficult to
construct $\widetilde{SU}_q(1,1)$ as a full blown locally compact quantum group. But this seems to be far from the
truth. The family of special functions produced in \cite{ErikJasper2} can be easily extended from the context of
$SU(1,1)$ to the context of $\widetilde{SU}(1,1)$, but the orthogonality and completeness of the relevant family of
functions needs some non-trivial special function theory that can be found in \cite{Salam} and \cite{Cic} (see
Proposition \ref{comult.ortho1}). This is in fact the only place where we heavily rely on the theory of
$q$-hypergeometric not covered in the appendix of this paper. Given these special functions it is not hard to define
the comultiplication of $\widetilde{SU}_q(1,1)$, the hard part of the construction of $\widetilde{SU}_q(1,1)$ lies in
the proof of the coassociativity of the comultiplication. Similarly, the left Haar weight is pretty easy to define, but
checking that it is left invariant in the sense of \cite{VaKust2} requires some non-trivial quantum group techniques.

\smallskip

Whereas \cite{Kor} is mainly a motivational paper for this one, \cite{Wor7} is not only motivational (especially the
use of the reflection operator $u$ introduced in notation \ref{neumann.not1}) but also gives an alternative approach to
$\widetilde{SU}_q(1,1)$. The advantage of the approach of \cite{Wor7} lies mainly in the fact that the role of the
generators and relations is more explicit and constructional. In order to do so, Woronowicz also presents a beautiful
operator theoretic theory of balanced extensions, similar to the theory of self-adjoint extension of symmetric
operators. As a consequence of this all, it should not be so difficult to prove that the \cst-algebra introduced in
proposition \ref{mult.prop5} is generated by the generators and relations, introduced in \cite[Def.~1.1]{Wor7}, and in
the sense of \cite{Wor4}. We on the other hand take a more pragmatic stance and approach the subject by using  special
functions. As a disadvantage, it follows that the role of generators and relations is less clear (but still lurking in
the background), but as a clear advantage we can prove the coassociativity of the comultiplication, prove the left
invariance of the left Haar weight and also give an explicit formula for the multiplicative unitary. It also allows us
to circumvent the theory of balanced extensions and stick to more standard operator theoretic techniques. As a
consequence, the proofs of this paper are different from the ones in \cite{Wor7}. A notable exception is the use of an
operator, introduced in \cite{Wor7}, in the proof of the last (albeit simplest) step in the proof of Theorem
\ref{coass.thm1}.

\smallskip\smallskip

The paper is organized as follows. In the first section we introduce the Hopf$^*$-algebra associated to
$\widetilde{SU}_q(1,1)$ and realize it as a $^*$-algebra of unbounded operators on some Hilbert space. In the
subsequent section the underlying von Neumann algebra is introduced. We define the comultiplication in section 3 and
claim its coassociativity at the end of the section but the proof of this fact is given in  section 5 in order to
enhance the readability of the paper. In section 4 we construct the Haar weights and prove their invariance, thereby
proving that $\widetilde{SU}_q(1,1)$ is indeed a locally compact quantum group. We end section 4 by calculating the
underlying \cst-algebraic quantum group. The appendix contains most of the basic results that we will need from the
theory of special functions.

\bigskip

\bf Acknowledgement. \rm  We would like to thank S.~L.~Woronowicz for insightful discussions on quantum
$\widetilde{SU}(1,1)$, providing the preliminary version of \cite{Wor7} and giving a series of lectures on \cite{Wor7}
in Trondheim (April 2000), which the second author attended.

\bigskip\medskip

\section*{Notations and conventions.}

\bigskip

The set of all natural numbers, not including $0$, is denoted by $\NE$. Also, $\NN = \NE \cup \{0\}$.

Fix a number $q >0$.  Let $a \in \C$. If $n \in \NN$, the $q$-shifted factorial $(a;q)_n \in \C$ is defined as $(a;q)_n
= \prod_{i=0}^{n-1} (1 - q^i a)$ \ \ (so $(a;q)_0 = 1$). From now on, we assume that $q <1$. Then $(a;q)_\infty \in \C$
is defined as $(a;q)_\infty = \lim_{n\rightarrow\infty}\,(a;q)_n$. Using some basic infinite product theory, one checks
that this limit exists and that the function $a \in \C \mapsto (a;q)_\infty$ is analytic. Also, $(a;q)_\infty = 0$ if
and only if $a \in q^{-\NN}$. We also use the notation $(a_1,\ldots,a_m;q)_k = (a_1;q)_k\,\ldots\,(a_m;q)_k$ if
$a_1,\ldots,a_m \in \C$ and $k \in \NN \cup \{\infty\}$.

\smallskip

If $a,b,z \in \C$, we define
$$
\hspace{18ex}\Psis{a}{b}{q,z} = \Psi(a;b;q,z) := \sum_{n=0}^\infty
\,\frac{(a;q)_n\,(b\,q^n;q)_\infty}{(q\,;q)_n}\,(-1)^n\,q^{\frac{1}{2}n(n-1)}\, z^n \ . \hspace{12ex} (1)
$$
We collected some further basic information about these $\Psi$ function in the appendix. See \cite{Gas} for an
extensive treatment on $q$-hypergeometric functions.

\medskip

Let $B_1,\ldots,B_n$ be sets such that $B_i = \T$ or $B_i = \qz$. Set  $I = \{\,i \in \{1,\ldots,n\} \mid B_i =
\qz\,\}$. Consider a set $K \subseteq (\qz)^I$ and a function $f : \{\,(x_1,\ldots,x_n) \in B_1 \times \ldots \times
B_n \mid (x_i)_{i \in I} \in K \,\} \rightarrow \C $. Then we set $f(x) := 0$ for $(x_1,\ldots,x_n) \in B_1 \times
\ldots \times B_n$ such that $(x_i)_{i \in I} \not\in K$. \ \ (2)

\smallskip

If $f$ is a function, the domain of $f$ will be denoted by $D(f)$. If $X$ is a set, the identity mapping on $X$ will be
denoted by $\iota_X$ and most of the time even by $\io$. If $H$ is a Hilbert space, $1_H$ denotes $\io_H$.  The set of
all complex valued functions on $X$ is denoted by $\cF(X)$, the set of all elements in $\cF(X)$ having finite support,
is denoted by $\cK(X)$. Let $V$ be a vector space and $S$ a subset of $V$. Then $\langle S \rangle$ denotes the linear
span of $S$ in $V$.

Consider a locally compact space $\Om$ with a regular Borel measure $\mu$ on it. If $g \in \cL^2(\Om,\mu)$, we denote
the class of $g$ in $L^2(\Om,\mu)$ by $[g]$. If $f$ is a measurable function on $\Om$, we define the linear operator
$M_f$ in $L^2(\Om,\mu)$ such that $D(M_f) = \{\,[g] \mid g \in \cL^2(\Om,\mu) \text{ such that } f g \in \cL^2(\Om,\mu)
\,\}$ and $M_f([g]) = [f g]$ for all such classes $[g] \in D(M_f)$.

Let $S$, $T$ be two linear operators acting in a Hilbert space $H$. We say that $S \subseteq T$ if $D(S) \subseteq
D(T)$ and $S(v) = T(v)$ for all $v \in D(S)$. Following \cite{Wor7} we call $S$  balanced if $D(S) = D(S^*)$.

The symbol $\odot$ will be used to denote the algebraic tensor product of vector spaces and linear mappings. The symbol
$\ot$ on the other hand will denote the tensor product of Hilbert spaces,  von Neumann algebras and sufficiently
continuous linear mappings.

If $n \in \NE$ and  $X$ is locally compact space, $X^n$ will denote the $n$-fold product set $X^n = X \times \ldots
\times X$.

\smallskip

Consider von Neumann algebras $M$,$N$ acting on Hilbert spaces $H$,$K$ respectively and $\pi : M \rightarrow N$ a
normal $^*$-homomorphism. Let $T$ be a densely defined, closed linear operator in $H$ affiliated with $M$ (in the von
Neumann algebraic sense) and $T = U\,|T|$ the polar decomposition of $T$. Then there exists a unique positive operator
$P$ in $K$ such that $f(P) = \pi(f(|T|))$ for all $f \in \cL^\infty(\R^+)$. Now we set $\pi(T) = \pi(U)\,P$.

\bigskip\medskip

\sectie{The Hopf$^*$-algebra underlying quantum $\widetilde{SU}(1,1)$.} \label{section1}

\bigskip

In order to resolve the problems surrounding quantum $SU(1,1)$, Korogodsky proposed in \cite{Kor} to construct the
quantum version of $\widetilde{SU}(1,1)$ instead of constructing the quantum version of $SU(1,1)$ itself. He suggested
that the Hopf$^*$-algebra of this quantum group should be the one that we describe now.

\smallskip

Throughout this paper, we fix a number $0 < q < 1$. Define $\cA_q$ to be the unital $^*$-algebra generated by elements
$\al_0$, $\gamma_0$ and $e_0$ and relations
\begin{equation} \label{eqset1}
\begin{array}{lcl}
\al_0^\dag \al_0 - \gamma_0^\dag \gamma_0 = e_0 & & \al_0 \al_0^\dag - q^2 \, \gamma_0^\dag \gamma_0 = e_0
\\ \gamma_0^\dag \gamma_0 = \gamma_0\, \gamma_0^\dag & &
\\ \al_0 \,\gamma_0 = q \, \gamma_0\, \al_0 & \hspace{6ex} & e_0^\dag = e_0
\\ \al_0 \,\gamma_0^\dag = q \, \gamma_0^\dag \al_0  & \hspace{6ex} & e_0^2 = 1
\\  \al_0 \, e_0 = e_0 \,\al_0
\\ \gamma_0 \, e_0 = e_0 \,\gamma_0
\end{array}
\end{equation}
where $^\dag$ denotes the $^*$-operation on $\cA_q$ (in order to distinguish this kind of adjoint with the adjoints of
possibly unbounded operators in Hilbert spaces).

There exists a unique unital $^*$-homomorphism $\de_0 : \cA_q \rightarrow \cA_q \odot \cA_q$ such that
\begin{eqnarray}
\de_0(\al_0) & = & \al_0 \ot \al_0 + q \, (e_0 \, \gamma_0^\dag) \ot \gamma_0 \nonumber
\\ \de_0(\gamma_0) & = & \gamma_0 \ot \al_0 + (e_0 \al_0^\dag) \ot \gamma_0 \label{hopf.comult1}
\\ \de_0(e_0) & = & e_0 \ot e_0 \ . \nonumber
\end{eqnarray}
The pair $(\cA_q,\de_0)$ turns out to be a Hopf$^*$-algebra with counit $\vep_0$ and antipode $S_0$ determined by
$$\begin{array}{lcl}
S_0(\al_0) = e_0\,\al_0^\dag & \hspace{10ex} & \vep_0(\al_0) = 1
\\ S_0(\al_0^\dag) = e_0\, \al_0 & & \vep_0(\gamma_0) = 0
\\ S_0(\gamma_0) = - q \, \gamma_0 & & \vep_0(e_0) = 1
\\ S_0(\gamma_0^\dag) = - \frac{1}{q}\,\gamma_0^\dag
\\ S_0(e_0) = e_0 \ .
\end{array} $$

\bigskip

As always we want to represent this Hopf $^*$algebra $\cA_q$ by possibly unbounded operators in some Hilbert space in
order to produce  a locally compact quantum group in the sense of Definition 1 of the Introduction. In Proposition 2.4
of \cite{Kor}, Korogodsky produces a family of irreducible $^*$-representations of the $^*$-algebra $\cA_q$. We glue
part of this family of irreducible $^*$-representations together to a $^*$-representation of $\cA_q$ on one Hilbert
space that will be the Hilbert space that our locally compact quantum group will act upon.

\smallskip

For this purpose we define
$$I_q = \{\,-q^k \mid k \in \NE\,\} \cup \{\,q^k \mid k \in \Z\,\} \ .$$ Set $I_q^- = \{\,x \in I_q \mid x < 0\,\}$
and $I^+_q = \{\,x \in I_q \mid x > 0\,\}$. We will use the discrete topology on $I_q$.

Let $\T$ denote the group of complex numbers of modulus 1. We will consider the uniform measure on $I_q$ and the
normalized Haar measure on $\T$. Our $^*$-representation of $\cA_q$ will act in the Hilbert space $H$ defined by $H =
L^2(\T) \ot L^2(I_q)$.

\smallskip

We let $\xi$ denote the identity function on $I_q$, whereas $\ze$ denotes the identity function on $\T$.

If $p \in \qz$, we define $\sde_p \in \cF(I_q)$ such that $\sde_p(x) = \sde_{x,p}$ for all $x \in I_q$ (note that
$\sde_p = 0$ if $p \not\in I_q$). The family $(\,\sde_p \mid p \in I_q\,)$ is the natural orthonormal basis of
$L^2(I_q)$. Also recall the natural orthonormal basis $(\,\ze^m \mid m \in \Z\,)$ for $L^2(\T)$.

\medskip\smallskip

Instead of looking at the algebra $\cA_q$ as the abstract algebra generated by generators and relations we will use an
explicit realization of this algebra as linear operators on the dense subspace $E$ of $H$ defined by $E = \langle\,
\ze^m \ot \sde_x \mid m \in \Z, x \in I_q, \rangle \subseteq H$. Of course, $E$ inherits the inner product from $H$.
Let $\cL^+(E)$ denote the $^*$-algebra of adjointable operators on $E$ (see \cite[Prop.~2.1.8]{Schm}), i.e.
$$\cL^+(E) = \{\,T \in \End(E) \mid \exists S \in
\End(E),\forall v,w \in E : \langle T v , w \rangle = \langle v , S w \rangle\,\}\ .$$ The $^*$-operation in $\cL^+(E)$
(and $\cL^+(E \odot E)$ for that matter) will be denoted by $^\dag$. So if $T \in \cL^+(E)$, the operator $T^\dag \in
\cL^+(E)$ is defined to be the operator $S$ in the above definition.  It follows that $T^\dag \subseteq T^*$ where
$T^*$ is the usual adjoint of $T$ as an operator in the Hilbert space $H$. It also follows that $T$ is a closable
operator in $H$.

\smallskip

Define linear operators $\al_0$, $\gamma_0$, $e_0$ on $E$ such that
\begin{eqnarray}
\al_0(\ze^m \ot \sde_p ) & = & \sqrt{\sgn(p)+p^{-2}}\,\,\ze^m \ot \sde_{qp} \nonumber
\\ \gamma_0(\ze^m \ot \sde_p) & = &  p^{-1}\,\,\ze^{m+1} \ot  \sde_p \label{hopf.eq1}
\\ e_0(\ze^m \ot \sde_p) & = & \sgn(p) \,\, \ze^m \ot \sde_p \nonumber
\end{eqnarray}
for all $p \in I_q$, $m \in \Z$.

Then $\al_0$, $\gamma_0$ and $e_0$ belong to $\cL^+(E)$, $e_0^\dag = e_0$ and
\begin{equation} \label{hopf.eq2}
\begin{array}{clc}
\al_0^\dag\,\bigl(\ze^m \ot \sde_p \bigr) & = &  \sqrt{\sgn(p)+q^2\,p^{-2}}\,\,\ze^m \ot \sde_{q^{-1}p} \vspace{1.5ex}
\\ \gamma_0^\dag\,\bigl(\ze^m \ot \sde_p\bigr) & = &  p^{-1}\,\,\ze^{m-1} \ot  \sde_p
\end{array}
\end{equation}
for all $p \in I_q$, $m \in \Z$. Note that $\al_0^\dag\,\bigl(\ze^m \ot \sde_{-q} \bigr) = 0$.

\medskip

Then $\cA_q$ is the $^*$-subalgebra of $\cL^+(E)$ generated by $\al_0$, $\gamma_0$ and $e_0$. Since $\cL^+(E) \odot
\cL^+(E)$ is canonically embedded in $\cL^+(E \odot E)$, we obtain $\cA_q \odot \cA_q$ as a $^*$-subalgebra of $\cL^+(E
\odot E)$. As such, the comultiplication $\de$ is, according to Eqs.~(\ref{hopf.comult1}), given by
\begin{eqnarray} \label{hopf.comult2}
\de_0(\al_0) & = & \al_0 \odot \al_0 + q \, (e_0 \, \gamma_0^\dag) \odot \gamma_0 \nonumber
\\ \de_0(\gamma_0) & = & \gamma_0 \odot \al_0 + (e_0 \al_0^\dag) \odot \gamma_0
\\ \de_0(e_0) & = & e_0 \odot e_0 \ . \nonumber \,
\end{eqnarray}
where $\odot$ denotes the algebraic tensor product of linear mappings.

\bigskip\medskip

\sectie{The von Neumann algebra underlying quantum $\widetilde{SU}(1,1)$.}

\bigskip

In this section we introduce the von Neumann algebra acting on $H$ that underlies the von Neumann algebraic version of
the quantum group $\widetilde{SU}_q(1,1)$. The reason for looking first at the von Neumann algebraic picture is
twofold. The most important one is the lack of density conditions in the definition of von Neumann algebraic quantum
groups (Definition 1 of the introduction), these are automatically satisfied. Moreover, in this case the von Neumann
algebra turns out to be very simple. Towards the end of this paper we will use \cite[Prop.~1.6]{VaKust2} to produce the
\cst-algebraic version of $\widetilde{SU}_q(1,1)$.

\medskip

In order to get into the framework of operator algebras, we need to introduce the topological versions of the algebraic
objects $\al_0$, $\gamma_0$ and $e_0$ as possibly unbounded operators in the Hilbert space $H$.  So let $\al$ denote
the closure of $\al_0$, $\gamma$ the closure of $\gamma_0$ and $e$ the closure of $e_0$, all as linear operators in
$H$. So $e$ is a bounded linear operator on $H$, whereas $\al$ and $\gamma$ are unbounded, closed, densely defined
linear operators in $H$. Note also that $\al^*$ is the closure of $\al_0^\dag$ and that $\gamma^*$ is the closure of
$\gamma_0^\dag$.

\smallskip

\begin{lemma} \label{neumann.lem1}
Let $\cH$ be a Hilbert space. Then the  domains of $1_\cH \ot \al$, $1_\cH \ot \al^*$, $1_\cH \ot \gamma$ and $1_\cH
\ot \gamma^*$~coincide.
\end{lemma}
\begin{proof}
Since $1_\cH \ot \gamma$ is normal, $D(1_\cH \ot \gamma^*) = D(1_\cH \ot \gamma)$. We know that $\al_0^\dag \al_0 =
\gamma_0^\dag \gamma_0 + e_0$, thus $\|(1_\cH \odot \al_0) v\|^2 = \|(1_\cH \odot \gamma_0) v\|^2 + \|(1_\cH \ot
e)v\|^2$ for all $v \in \cH \odot E$. Since $1_\cH \ot e$ is bounded, $1_\cH \ot \al$ is the closure of $1_\cH \odot
\al_0$ and $1_\cH \ot \gamma$ is the closure of $1_\cH \odot \gamma_0$, it follows that $D(1_\cH \ot \al) = D(1_\cH \ot
\gamma)$. Using the equality $\al_0 \al_0^\dag = q^2\,\gamma_0^\dag \gamma_0 + e_0$, one proves in a similar way that
$D(1_\cH \ot \al^*) = D(1_\cH \ot \gamma)$.
\end{proof}

\medskip

All these operators can easily be realized as a combination of shift and multiplication operators.

\smallskip

Consider $p \in \qz$. We define a translation operator $T_p$ on $\cF(\T \times I_q)$ such that for $f \in \cF(\T \times
I_q)$, $\lambda \in \T$ and $x \in I_q$, we have that $(T_p f)(\lambda,x) =  f(\lambda,px)$. By discussion (2) in
 Notations and conventions, we get that $(T_p f)(\lambda,x) = 0$ if $p x \not\in I_q$. If $p,t \in I_q$ and $g \in
\cF(\T)$, then $T_p(g \ot \sde_t) = g \ot \sde_{p^{-1} t}$, thus,  $T_p(g \ot \sde_t) = 0$ if $p^{-1} t \not\in I_q$.

\smallskip

For instance, note that if $f \in \cF(\T \times I_q)$, $g \in \cF(\T)$ and $\lambda \in \T$,
\begin{trivlist}
\item[\ \,$\bullet$]\ $(T_{q^{-1}} f)(\lambda,-q) = 0$  and $T_q(g \ot \sde_{-q}) = 0$,
\vspace{0.5mm}
\item[\ \,$\bullet$]\ $(T_{-1} f)(\lambda,x) = 0$ for all  $x \in I_q$ such that $x \geq 1$ ;
$T_{-1}(g \ot \sde_p) = 0$ for all $p \in I_q$ such that $p \geq 1$.
\end{trivlist}

\smallskip

\begin{notation} \label{neumann.not1}
Define the mapping $\rho : \qz \rightarrow B(H)$ such that $\rho_p$ equals the partial isometry on $H$ induced by $T_p$
for all $p \in \qz$. Let us also single out the following special cases:
\begin{trivlist}
\item[\ \,1)] Define $w = \rho_{q^{-1}}$, which is an isometry on $B(H)$.
\item[\ \,2)] Define $u = \rho_{-1}$, which is a self-adjoint partial isometry on $B(H)$.
\end{trivlist}
\end{notation}

\smallskip

In terms of multiplication and shift operators, our closed linear operators in $H$ are easily recognized as
\begin{equation} \label{neumann.eq1}
\al = w\,(1 \ot M_{\sqrt{\sgn(\xi) + \xi^{-2}}}) \hspace{8ex}  \gamma = M_{\ze} \ot M_{\xi^{-1}} \hspace{8ex} e = 1 \ot
M_{\sgn{\xi}}\ .
\end{equation}
These tensor products are obtained by closing the algebraic tensor product mappings with respect to the norm topology
on $H$.

\medskip

Let us recall the following natural terminology. If $T_1, \ldots, T_n$ are closed, densely defined linear operators in
$H$, the von Neumann algebra $N$ on $H$ generated by $T_1, \ldots ,T_n$ is the one such that
$$N' = \{\,x \in B(H) \mid x T_i \subseteq T_i x \text{ and } x T_i^* \subseteq T_i^* x \text{ for } i =
1,\ldots,n\,\}\ .$$ Almost by definition, $N$ is the smallest von Neumann algebra acting on $H$ so that
$T_1,\ldots,T_n$ are affiliated with $M$ in the von Neumann algebraic sense. If $w_1,\ldots,w_n$ are the partial
isometries obtained from the polar decompositions of $T_1,\ldots,T_n$ respectively, then $N$ is also the von Neumann
algebra on $H$ generated by
\begin{equation} \label{neumann.eq2}
\bigcup_{i=1}^n\,\, \{w_i\}\, \cup\, \{\,f(T_i^* T_i) \mid f \in \cL^\infty(\si(T_i^* T_i))\,\} \ .
\end{equation}

It is now very tempting to define the von Neumann algebra underlying quantum $\widetilde{SU}_q(1,1)$ as the von Neumann
algebra generated by $\al$, $\gamma$ and $e$. However, for reasons that will become clear later (see the discussion in
the beginning of the next section and the remark after Proposition \ref{comult.prop7}), the underlying von Neumann
algebra will be the one generated by $\al$, $\gamma$, $e$ and $u$.

\smallskip

\begin{definition} \label{neumann.def1}
We define $M_q$ to be the von Neumann algebra on $H$ generated by $\al$, $\gamma$, $e$ and $u$.
\end{definition}

\smallskip

So $M_q$ will be the von Neumann algebra underlying $\widetilde{SU}_q(1,1)$. For convenience, we will also introduce
$M_q^\diamond$ as the von Neumann algebra on $H$ generated by $\al$, $\gamma$ and $e$.

\smallskip

\begin{lemma} \label{neumann.lem2}
\begin{trivlist}
\item[\ \,1)] $M_q^\diamond$ is generated by $\{w\} \cup \{\,M_f \mid f \in \cL^\infty(\T \times I_q)\,\}$,
\item[\ \,2)] $M_q$ is generated by $\{w,u\} \cup \{\,M_f \mid f \in \cL^\infty(\T \times I_q)\,\}$,
\item[\ \,3)] $M_q = L^\infty(\T) \ot B(L^2(I_q))$.
\end{trivlist}
\end{lemma}
\begin{proof}
By Eq.~(\ref{neumann.eq2}) we know that $M_q^\diamond$ is generated by the set $\{w,M_{\ze \ot 1},M_{1 \ot
\sgn(\xi)}\}\,\cup $\newline$\{\,f(\al^* \al) | f \in \cL^\infty(\si(\al^* \al))\,\}\,\cup\, \{\,f(\gamma^* \gamma) | f
\in \cL^\infty(\si(\gamma^* \gamma))\,\}$. Since $\gamma^* \gamma = 1 \ot M_{\xi^{-2}}$, it follows that the von
Neumann algebra generated by $\{M_{1 \ot \sgn(\xi)}\}\,\cup\,  \{\,f(\gamma^* \gamma) | f \in \cL^\infty(\si(\gamma^*
\gamma))\,\}$ equals the von Neumann algebra generated by $\{M_{1 \ot \sgn(\xi)}\}\,\cup\,  \{\,M_{1 \ot f(|\xi|)} | f
\in \cL^\infty(I_q^+)\,\}$. Thus, the von Neumann algebra generated by $\{M_{1 \ot \sgn(\xi)}\}\,\cup\,  \{\,f(\gamma^*
\gamma) | f \in \cL^\infty(\si(\gamma^* \gamma))\,\}$ equals $\{\,M_{1 \ot f} \mid f \in \cL^\infty(I_q)\,\}$. Because
$\al^* \al = 1 \ot M_{\sgn(\xi) +  \xi^{-2}}$, it now follows that $M_q^\diamond$ is generated by $\{w,M_{\ze \ot 1}\}
\, \cup \, \{\,M_{1 \ot f} \mid f \in \cL^\infty(I_q)\,\}$. Therefore the first statement holds, and as a consequence
also the second one holds.

\smallskip

Define the partial isometries $w_0,u_0$ on $L^2(I_q)$ such that $(w_0 f)(x) = f(q^{-1} x)$ and $(u_0 f)(x) = f(-x)$ for
all $f \in L^2(I_q)$ and $x \in I_q$. Letting $M_0$ denote the von Neumann algebra generated by
$\{w_0,u_0\}\,\cup\,L^\infty(I_q)$, the second statement implies that $M_q = L^\infty(\T) \ot M_0$. It is not so
difficult to see that every rank one projector of the form $L^2(I_q) \rightarrow L^2(I_q) : v \mapsto \langle v ,
\sde_p \rangle\,\sde_t$, where $p,t \in I_q$, belongs to $M_0$, thus $M_0 = B(L^2(I_q))$.
\end{proof}

\medskip

So we get in particular that $\rho_p \in M_q$ for all $p \in \qz$. Let us collect some further elementary results about
the map $\rho$. If $p \in \qz$, we define $i_p = \rho_p^* \rho_p$, the initial projection of $\rho_p$. Consider $f \in
\cL^2(\pr)$ and define $g \in  \cL^2(\pr)$ such that for $\lambda \in \T$, $x \in I_q$, we have that $g(\lambda,x) =
f(\lambda,x)$  if $p^{-1} x \in I_q$ and $g(\lambda,x) = 0$ if $p^{-1} x \not\in I_q$.  Then $i_p([f]) = [g]$.

\smallskip

\begin{lemma} \label{neumann.lem3}
Consider $p,t \in \qz$ and $f \in \cL^\infty(\T \times I_q)$, then
$$\rho_p\,\rho_t = \rho_{p\,t}\,i_t \hspace{4ex} \text{and} \hspace{4ex} \rho_p^* = \rho_{p^{-1}}\ $$
and
$$\rho_p\,\,M_f = M_{T_p f}\,\,\rho_p \hspace{4ex} \text{and} \hspace{4ex} \rho_p\,\,M_{T_{p^{-1}}f} =
M_f\,\,\rho_p\ .$$
\end{lemma}
\begin{proof}
We only prove the first equality, the other ones are even more straightforward to check. Take $f \in \cL^2(\pr)$ such
that $f(\lambda,x) = 0$ for all $\lambda \in \T$ and $x \in I_q$ such that $t^{-1}\,x \not\in I_q$. Now $(\rho_p\,
\rho_t)([f]) = [T_p(T_t(f))]$, whereas $\rho_{p\,t}([f]) = [T_{p\,t}(f)]$. Let $y \in I_q$ and $\nu \in \T$. We
consider 4 different cases:

\vspace{-1ex}

\begin{trivlist}
\item[\ \,(1)] If $p y \in I_q$ and $p\,t y \in I_q$, then $(T_p(T_t(f))(\nu,y) =
T_t(f)(\nu,py) = f(\nu,p\,ty)$. Thus $T_{p\,t}(f)(\nu,y) = f(\nu,p\,t y) = (T_p(T_t(f))(\nu,y)$.
\item[\ \,(2)] Now suppose that $p y \not\in I_q$ and $p\,t y \in I_q$. Then $(T_p(T_t(f)))(\nu,y)
= 0$. On the other hand, $T_{p\,t}(f)(\nu,y) = f(\nu,p\,t y)$. Since $t^{-1}(p\,t y) = p y \not\in I_q$, we have that
$f(\nu,p\,t y) = 0$. Thus $T_{p\,t}(f)(\nu,y) = 0 = (T_p(T_t(f))(\nu,y)$.
\item[\ \,(3)] Suppose that $p y \in I_q$ and $p\,t y \not\in I_q$  . Then $\bigl(T_p(T_t(f))\bigr)(\nu,y)
= T_t(f)(\nu,p y) = 0$ since $t(py) \not\in I_q$. Thus $T_{p\,t}(f)(\nu,y) = 0 = \bigl(T_p(T_t(f))\bigr)(\nu,y)$.
\item[\ \,(4)] If $p y \not\in I_q$ and $p\,t y \not\in I_q$, then clearly
$T_{p\,t}(f)(\nu,y) = 0 = \bigl(T_p(T_t(f))\bigr)(\nu,y)$.
\end{trivlist}
So we see that $T_p(T_t(f)) = T_{p\,t}(f)$ for this kind of function $f$. Thus $\rho_p\,\rho_t\,i_t =
\rho_{p\,t}\,i_t$.
\end{proof}

\smallskip

So we see that $\rho$ is almost a group representation, but not quite. It is possible to use the results in the
previous lemma to perform  a partial cross product construction. Since we do not need this approach in the rest of the
paper, we will not pursue this matter any further. Instead we focus on a more useful and simpler picture of $M_q$:

\smallskip

For every $p,t \in I_q$ and $m \in \Z$ we define $\Phi(m,p,t) = \rho_{p^{-1} t}\,M_{\zeta^m \ot \sde_t} \in M_q$. So if
$x \in I_q$ and $r \in \Z$, then
\begin{equation}
\Phi(m,p,t)\,(\ze^r \ot \sde_x) = \sde_{x,t}\,\, \ze^{m+r} \ot \sde_p\ .\label{neumann.eq3}
\end{equation}

Define $M_q^\circ = \langle \,\Phi(m,p,t) \mid m \in \Z, p,t \in I_q\,\rangle$. Using Eq.~(\ref{neumann.eq3}), it is
easy to see that the family $(\,\Phi(m,p,t) \mid m \in \Z, p,t \in I_q\,)$ is a linear basis of $M_q^\circ$.

\medskip

The multiplication and $^*$-operation are easily expressed in terms of these basis elements:
\begin{eqnarray}
& & \Phi(m_1,p_1,t_1)\, \Phi(m_2,p_2,t_2) = \sde_{p_2,t_1}\,\,\Phi(m_1+m_2,p_1,t_2) \nonumber
\\ & & \Phi(m,p,t)^* = \Phi(-m,t,p) \label{neumann.eq4}
\end{eqnarray}
for all $m,m_1,m_2 \in \Z$, $p,p_1,p_2,t,t_1,t_2 \in I_q$. So we see that $M_q^\circ$ is a $\si$-weakly dense
sub$^*$-algebra of $M_q$.

\bigskip\medskip

\sectie{The comultiplication of quantum $\widetilde{SU}(1,1)$.}

\bigskip

In this section we introduce the comultiplication of $\widetilde{SU}_q(1,1)$. In the first part we start with a
motivation for the formulas appearing in Definition \ref{comult.def1}. Although the discussion is not really needed in
the build up of $\widetilde{SU}_q(1,1)$, it is important and clarifying to know how we arrived at the formulas in
Definition \ref{comult.def1}.

\medskip

But first we introduce two auxiliary functions
\begin{trivlist}
\item[\ \,\,(1)] $\chi : \qz \rightarrow \Z$ such that $\chi(x) = \log_q(|x|)$ for all $x \in \qz$,
\smallskip
\smallskip
\item[\ \,\,(2)] $\kappa : \R \rightarrow \R$ such that $\kappa(x) = \sgn(x)\,x^2$ for all $x \in \R$.
\end{trivlist}

\medskip\smallskip

Our purpose is to define a comultiplication $\de : M_q \rightarrow M_q \ot M_q$. Assume for the moment that this has
already been done. It is natural to require $\de$ to be closely related to the comultiplication $\de_0$ on $\cA_q$ as
defined in Eqs.~(\ref{hopf.comult2}). The least that we expect is $\de_0(T_0) \subseteq \de(T)$ and $\de_0(T_0^\dag)
\subseteq \de(T)^*$ for $T = \al,\gamma$. In the rest of this discussion we will focus on the inclusion
$\de_0(\gamma_0^\dag \gamma_0) \subseteq \de(\gamma^* \gamma)$, where $\de_0(\gamma_0^\dag \gamma_0) \in \cL^+(E \odot
E)$.

\smallskip

Because $\gamma^* \gamma$ is self-adjoint, the element $\de(\gamma^* \gamma)$ would also be self-adjoint. So the hunt
is on for self-adjoint extensions of the explicit operator $\de_0(\gamma_0^\dag \gamma_0)$. Unlike in the case of
quantum $E(2)$ (see \cite{Wor6}), the operator $\de_0(\gamma_0^\dag \gamma_0)$ is not essentially self-adjoint. But it
was already known in  \cite{Kor}  that $\de_0(\gamma_0^\dag \gamma_0)$ has self-adjoint extensions (this follows easily
because the operator in (\ref{comult.eq1}) commutes with complex conjugation, implying that the deficiency spaces are
isomorphic).

Let us make a small detour to quantum $SU(1,1)$.  In this case, the operators $\al_0$ and $\gamma_0$ are replaced by
their restrictions to $L^2(\T \ot I_q^+)$. Then $\de_0(\gamma_0^\dag \gamma_0)$ still has self-adjoint extensions for
each of which the domain is obtained by imposing a boundary condition on functions in its domain (this boundary
condition is a simple relation between the function and its Jackson derivative in the limit towards 0, see \cite[Eq.~
(6.2)]{Kor}). However, in this case the closure of the operator $\de_0(\al_0)$ does not leave the domain of such a
self-adjoint extension invariant because it distorts the boundary condition. From a practical point of view, it should
also be said that an explicit manageable spectral decomposition in terms of special functions for these self-adjoint
extensions is missing, cf.~\cite[Rem.~2.7]{ErikJasper1}.

\smallskip

Now we return to the case of quantum $\widetilde{SU}(1,1)$. Although $\de_0(\gamma_0^\dag \gamma_0)$ has a self-adjoint
extension, it is not unique. We have to make a choice for this self-adjoint extension, but we cannot extract the
information necessary to make this choice from $\al$ and $\gamma$ alone. This is why we do not work with $M_q^\diamond$
but with $M_q$ which has the above extra extension information contained in the element $u$. These kind of
considerations were already present in \cite{Wor8} and were also introduced in \cite{Wor7} for quantum
$\widetilde{SU}_q(1,1)$. In this paper, this principle is only lurking in the background but it is treated in a
fundamental and rigorous way in \cite{Wor7}.

\medskip

Now we get into slightly more detail in our discussion about the extension of $\de_0(\gamma_0^\dag \gamma_0)$. Define a
linear map $L : \cF(\T \times I_q \times \T \times I_q) \rightarrow \cF(\T \times I_q \times \T \times I_q)$ such that
\begin{eqnarray*}
(L f)(\lambda,x,\mu,y) & = & [\,x^{-2}(\sgn(y) + y^{-2}) + (\sgn(x) + q^2\,x^{-2})\,y^{-2}\,]\, f(\lambda,x,\mu,y)
\\ &  & +\, \sgn(x)\,q^{-1}\,\bar{\lambda} \mu\, x^{-1} y^{-1} \,\sqrt{(\sgn(x)+x^{-2})(\sgn(y)+y^{-2})} \,
f(\lambda,qx,\mu,qy)
\\ &  & +\, \sgn(x)\,q\,\lambda \bar{\mu}\, x^{-1} y^{-1} \,\sqrt{(\sgn(x)+q^2 x^{-2})(\sgn(y)+q^2 y^{-2})} \,
f(\lambda,q^{-1} x,\mu,q^{-1} y)
\end{eqnarray*}
for all $\lambda,\mu \in \T$ and $x,y \in I_q$. A straightforward calculation reveals that if $f \in E \odot E$, then
$\de_0(\gamma_0^\dag \gamma_0)\,[f] = [L(f)]$. From this, it is a standard exercise to check that $[f] \in
D(\de_0(\gamma_0^\dag \gamma_0)^*)$ and $\de_0(\gamma_0^\dag \gamma_0)^* [f] = [L(f)]$ if $f \in \cL^2(\T \times I_q
\times \T \times I_q)$ and $L(f) \in \cL^2(\T \times I_q \times \T \times I_q)$ (without any difficulty, one can even
show that $D(\de_0(\gamma_0^\dag \gamma_0)^*)$ consists precisely of such elements $[f]$).

\medskip

If $\theta \in \qz$, we define $\ell_\theta' = \{\,(\lambda,x,\mu,y) \in \T \times I_q \times \T \times I_q \mid y =
\theta x \,\}$. We consider $L^2(\ell_\theta')$ naturally embedded in $L^2(\T \times I_q \times \T \times I_q)$. It
follows easily from the above discussion that $\de_0(\gamma_0^\dag \gamma_0)^*$ leaves $L^2(\ell_\theta')$ invariant.
Thus, if $T$ is a self-adjoint extension of $\de_0(\gamma_0^\dag \gamma_0)$, the obvious inclusion $T \subseteq
\de_0(\gamma_0^\dag \gamma_0)^*$ implies that $T$ also leaves $L^2(\ell_\theta')$ invariant.

Therefore every self-adjoint extension $T$ of $\de_0(\gamma_0^\dag \gamma_0)$ is obtained by choosing a self-adjoint
extension $T_\theta$ of the restriction of $\de_0(\gamma_0^\dag \gamma_0)$ to $L^2(\ell_\theta')$ for every $\theta \in
\qz$ and setting $T = \oplus_{\theta \in \qz} T_\theta$.

\medskip

Therefore fix $\theta \in \qz$. Define $J_\theta = \{\,z \in I_{q^2} \mid \kappa(\theta)\,z \in I_{q^2}\,\}$ which is a
$q^2$-interval around $0$. On $J_\theta$ we define a measure $\nu_\theta$ such that $\nu_\theta(\{x\}) = |x|$ for all
$x \in J_\theta$.

Now define the unitary transformation $U_\theta : L^2(\T \times \T \times J_\theta) \rightarrow L^2(\ell_\theta')$ such
that $U_\theta([f]) = [g]$ where $f \in \cL^2(\T \times \T \times J_\theta)$ and $g \in \cL^2(\ell_\theta')$ are such
that
\begin{equation} \label{comult.eq3}
g(\lambda,z,\mu,\theta z) = (\lambda \bar{\mu})^{\chi(z)} \,  (-\sgn(\theta z))^{\chi(z)}\, |z| \,
f(\lambda,\mu,\kappa(z))
\end{equation}
for all $\lambda,\mu \in \T$ and $z \in I_q$ such that $\theta z \in I_q$.

Define the linear operator $L_\theta : \cF(J_\theta) \rightarrow \cF(J_\theta)$ such that
\begin{eqnarray} \label{comult.eq1}
(L_\theta f)(x) & = & \frac{1}{\theta^2 x^2}\,\bigl(\, - \sqrt{(1+x)(1+\kappa(\theta)\, x)} \, f(q^2 x) - q^2\,
\sqrt{(1+q^{-2} x)(1+q^{-2} \kappa(\theta)\, x)}\, f(q^{-2} x) \nonumber
\\ & & \hspace{7ex} [(1+ \kappa(\theta)\,x) + q^2(1+ q^{-2} x)]\, f(x)\, \bigr)
\end{eqnarray}
for all $f \in \cF(J_\theta)$ and $x \in J_\theta$.

Then an easy calculation shows that $U_\theta^*\,\de_0(\gamma_0^\dag \gamma_0)\!\!\restriction_{E \odot E}\, U_\theta =
1 \odot (L_\theta \!\!\restriction_{\cK(J_\theta)})$. So our problem is reduced to finding self-adjoint extensions of
$L_\theta \!\!\restriction_{\cK(I_q)}$. This operator $L_\theta \!\!\restriction_{\cK(I_q)}$ is a second order
$q$-difference operator for which eigenfunctions in terms of $q$-hypergeometric functions are known.

\smallskip

We can use a reasoning similar to the one in  \cite[Sec.~2]{ErikJasper1} to get hold of the self-adjoint extensions of
$L_\theta \!\!\restriction_{\cK(I_q)}$:

Let $\be \in \T$. Then we define a linear operator $L^\be_\theta : D(L_\theta^\be)\subseteq L^2(J_\theta,\nu_\theta)
\rightarrow L^2(J_\theta,\nu_\theta)$ such that
$$D(L^\be_\theta) = \{\,f \in L^2(J_\theta,\nu_\theta) \mid L_\theta(f) \in L^2(J_\theta,\nu_\theta), f(0+) = \be\,f(0-) \text{ and }
(D_q f)(0+) = \be\,(D_q f)(0-)\,\}$$ and $L_\theta^\be$ is the restriction of $L_\theta$ to $D(L_\theta^\be)$. Here,
$D_q$ denotes the Jackson derivative, that is, $(D_q f)(x) = (f(qx) - f(x))/(q-1)x$ for $x \in J_\theta$. Also, $f(0+)
= \be\,f(0-)$ is an abbreviated form of saying that the limits $\lim_{x\uparrow0} f(x)$ and $\lim_{x\downarrow0} f(x)$
exist and $\lim_{x\downarrow0} f(x) = \be\,\lim_{x\uparrow0} f(x)$.

Then $L_\theta^\be$ is a self-adjoint extension of $L_\theta \!\!\restriction_{\cK(I_q)}$. If $\be,\be' \in \T$ and
$\be \not= \be'$, then $L_\theta^\be \not= L_\theta^{\be'}$.

\medskip

It is tempting to use the extension $L^1_\theta$ to construct our final self-adjoint extension for $\de_0(\gamma_0^\dag
\gamma_0)$ (although there is no apparent reason for this choice). However, in order to obtain a coassociative
comultiplication, it turns out that we have to use the extension $L^{\sgn(\theta)}_\theta$ to construct our final
self-adjoint extension. This is reflected in the fact that the expression $s(x,y)$ appears in the formula for $a_p$ in
Definition \ref{comult.def1}.

\medskip

This all would  be only a minor achievement if we could not go any further. But the results and techniques used in the
theory of $q$-hypergeometric functions will even allow us to find an explicit orthonormal basis consisting of
eigenvectors of $L^{\sgn(\theta)}_\theta$.  These eigenvectors are, up to a  unitary transformation, obtained by
restricting the functions $a_p$ in Definition \ref{comult.def1} to $\ell_\theta$, which is introduced after this
definition. The special case $\theta = 1$ was already known to Korogodsky (see \cite[Prop.~A.1]{Kor}), although the
proof in \cite{Kor} seems to contain quite a few gaps. For instance, the presentation in \cite{Cic} shows that a lot of
more care has to be taken to solve these kind of problems.

\bigskip

In order to compress the formulas even further, we introduce two other auxiliary functions
\begin{trivlist}
\item[\ \,\,(1)] $\nu : \qz  \rightarrow \R^+$ such that $\nu(t) = q^{\frac{1}{2}(\chi(t)-1)(\chi(t) -
2)}$ for all $t \in \qz$.
\smallskip
\item[\ \,\,(2)] Another auxiliary function $s : \R_0 \times \R_0 \rightarrow \{-1,1\}$ is defined such that
$$s(x,y) = \begin{cases} -1 & \text{if } x > 0 \text{ and } y < 0
\\ 1 & \text{ if } x < 0 \text{ or } y > 0  \end{cases}$$
for all $x,y \in \R_0$.
\end{trivlist}

\smallskip

Let us also collect some basic manipulation rules. Therefore take $x,y,z \in \R_0$.
\begin{trivlist}
\item[\ \,(1)] If $x > 0$ then $s(x,y) = \sgn(yz)\,s(x,z)$ and $s(y,x) = s(z,x)$. If $x < 0$, then $s(x,y) = s(x,z)$
and $s(y,x) = \sgn(yz)\,s(z,x)$. So if $\sgn(xyz) = 1$, then $s(x,y) = s(x,z)$.
\smallskip
\item[\ \,(2)] It is clear that $s(x,y) = -\sgn(x)\,s(x,-y)$ and $s(x,y) = \sgn(y)\,s(-x,y)$.
\smallskip
\item [\ \,(3)] One easily checks that $s(x,y) = \sgn(xy) \, s(y,x)$. As a consequence we get also that $s(-x,y) =
s(-y,x)$ and $s(x,-y) = s(y,-x)$.
\end{trivlist}

\bigskip

If $a,b,z \in \C$, we define $\Psis{a}{b}{z}$ and $\Psi(a;b;z)$ to be equal to $\Psis{a}{b}{q^2;z}$.

\medskip\smallskip

We will also use the normalization constant $c_q = (\sqrt{2}\,q\,(q^2,-q^2;q^2)_\infty)^{-1}$.

\smallskip

\begin{definition} \label{comult.def1}
If $p \in I_q$, we define a function $a_p : I_q \times I_q \rightarrow \R$ such that $a_p$ is supported on the set
$\{\,(x,y) \in I_q \times I_q \mid \sgn(xy) = \sgn(p)\,\}$ and $a_p(x,y)$ is given by
$$c_q\,s(x,y)\,(-1)^{\chi(p)}\,(-\sgn(y))^{\chi(x)}\,|y|\,\,\nu(py/x)\,\,
\,\, \sqrt{\frac{(-\kappa(p),-\kappa(y);q^2)_\infty}{(-\kappa(x);q^2)_\infty}}\,\, \Psis{-q^2/\kappa(y)}{q^2
\kappa(x/y)}{q^2 \kappa(x/p)}
$$
for all $(x,y) \in I_q \times I_q$ satisfying $\sgn(xy) = \sgn(p)$.
\end{definition}

\bigskip

The presence of the expression $(-\sgn(y))^{\chi(x)}$ in this same formula can be traced back to the defining formula
(\ref{comult.eq3}) for $U_\theta$.

\medskip

The extra vital information that we need is contained in the following proposition. For $\theta \in \qz$ we define
$\ell_\theta = \{\,(x,y) \in I_q \times I_q \mid y = \theta x \,\}$.

\smallskip

\begin{proposition} \label{comult.ortho1}
Consider $\theta \in \qz$. Then the family  $(\,a_p\!\!\restriction_{\ell_\theta}\, \mid p \in I_q \text{ such that }
\sgn(p) = \sgn(\theta)\,)$ is an orthonormal basis for $l^2(\ell_\theta)$.
\end{proposition}

\smallskip

This result implies also a dual result, stemming from the following simple duality principle (in special function
theory, these are referred to as dual orthogonality relations). Consider a set $I$ and suppose that $l^2(I)$ has an
orthonormal basis $(e_j)_{j \in J}$. For every $i \in I$, we define a function $f_i$ on $J$ by $f_i(j) = e_j(i)$. Then
$(f_i)_{i \in I}$ is an orthonormal basis for $l^2(J)$. If we apply this principle to the line $\ell_\theta$, the
previous proposition implies the next one.

\begin{proposition} \label{comult.ortho2}
Consider $\theta  \in \qz$ and define $J = I_q^+$ if $\theta > 0$ and $J = I_q^-$ if $\theta < 0$. For every $(x,y) \in
\ell_\theta$ we define the function $e_{(x,y)} : J \rightarrow \R$ such that $e_{(x,y)}(p) = a_p(x,y)$ for all $p \in
J$. Then the family $(\,e_{(x,y)} \mid (x,y) \in \ell_\theta\,)$ forms an orthonormal basis for $l^2(J)$.
\end{proposition}

\smallskip

The first of these two propositions will be shortly used to define the comultiplication on $M_q$, both of them pop up
in the proof of the left invariance of the Haar weight.

\smallskip

For the proofs of Propositions \ref{comult.ortho1} and \ref{comult.ortho2} we refer to the literature. For $1\geq
\theta>0$ Proposition \ref{comult.ortho2} is \cite[Thm.~4.1]{Cic} and Proposition \ref{comult.ortho1} follows from
Corollary 4.2 and the following remark of \cite{Cic} in base $q^2$ and with $(c,q^{2\alpha})$ replaced by
$(1,\kappa(\theta))$. For $\theta\geq 1$ we derive Propositions \ref{comult.ortho1} and \ref{comult.ortho2} similarly
from \cite[Thm.~4.1 and Cor.~4.2]{Cic} in base $q^2$ and with $(c,q^{2\alpha})$ replaced by $(1,\kappa(\theta)^{-1})$.
The case $\theta\geq 1$  can also be reduced to the case $0<\theta\leq 1$ using the second symmetry property of
Proposition 3.5. The functions $a_p(x,y)\!\!\restriction_{\ell^\theta}$ are eigenfunctions of the operator $L_\theta$
as considered in (\ref{comult.eq1}). The general theory used in the first half of \cite{Cic} is explained in
\cite{Erik}. For $\theta<0$ the statements in Propositions \ref{comult.ortho1} and \ref{comult.ortho2} reduce to
statements on specific classes of orthogonal polynomials. For $\theta<0$ the orthogonality relations in Proposition
\ref{comult.ortho1} are directly obtainable from the orthogonality relations for the Al-Salam and Carlitz polynomials
$U_n^{(a)}(x;q)$ in base $q^2$ with $a$ replaced by $\kappa(\theta)$, see Al-Salam's and Carlitz's original paper
\cite{Salam} or references in \cite{KS}. This time we have that $a_p(x,y)\!\!\restriction_{\ell^\theta}$ are
eigenfunctions of the operator $L_\theta$, due to the second order $q$-difference equation for the Al-Salam and Carlitz
polynomials, see \cite{KS}. The dual result in Proposition \ref{comult.ortho2}, which follows immediately from
Proposition \ref{comult.ortho1} for $\theta<0$ since the corresponding moment problem is determinate, can also be
matched to orthogonal polynomials after applying elementary transformation formulas for basic hypergeometric series.
The dual orthogonality relations split up in three summation results. Two of these summations can be matched to the
orthogonality relations for the $q$-Charlier polynomials given in \cite[Eq.~(3.23.3)]{KS}, but for different
parameters. The remaining summation is an easy consequence of \cite[Eq.~(1.3.16)]{Gas}.

The functions $a_p(x,y)$ can be thought of as Clebsch-Gordan coefficients for the tensor product decomposition of the
representation described in (2.1) with itself. Note that the corresponding Clebsch-Gordan coefficients for the quantum
$SU(2)$ group are given in terms of Wall polynomials, see Koornwinder \cite[Rem.~4.2]{Koorn}, and we can write the Wall
polynomials in terms of the functions $a_p(x,y)$, precisely for $y\in -q^{-\NN}$, using \cite[Eq.~(III.3)]{Gas}.

\medskip

There is a nice symmetry in $a_p(x,y)$ with respect to interchanging $x$, $y$ and $p$. One of this symmetries is easy
to see while the proof of the other one requires an extra lemma.

\medskip

\begin{lemma}
Consider $x,y,p \in I_{q^2}$ and $m \in \Z$ such that $y/x = \sgn(p)\,q^{2m}$. Then
$$
\sqrt{\frac{(-y;q^2)_\infty}{(-x;q^2)_\infty}}\,\,\Psis{-q^2/y}{q^2x/y}{q^2x/p} =
\sgn(p)^{\log_{q^2}(|p|)+1}\,(-q^2/p)^m\,\,\sqrt{\frac{(-x;q^2)_\infty}{(-y;q^2)_\infty}}\, \,
\Psis{-q^2/x}{q^2y/x}{q^2y/p} \ .
$$
\end{lemma}
\begin{proof}
If $p > 0$, this follows easily from Proposition \ref{app.prop1}(1). Now assume that $p < 0$. By Result \ref{app.res2},
the above equation is equivalent to
$$\Psis{-q^2/p}{q^2x/p}{q^2x/y} =
\sgn(p)^{\log_{q^2}(|p|)+1}\,(-q^2/p)^m\,\,\frac{(-x;q^2)_\infty}{(-y;q^2)_\infty}\, \, \Psis{-q^2/p}{q^2y/p}{q^2y/x} \
.$$ This equality follows easily from Proposition \ref{app.prop1}(2).
\end{proof}

\medskip

With this lemma in hand, one verifies the second of the next symmetries. Use Result \ref{app.res2} to prove the first
symmetry, the third one is then a consequence of the previous two symmetries. For the first two symmetries you will
also need two of the manipulation rules for $s(x,y)$ discussed before Definition \ref{comult.def1}, namely the last one
of (1) and the first one of (3).

\smallskip

\begin{proposition} \label{comult.prop2}
If $x,y,p \in I_q$, then
\begin{eqnarray}
a_p(x,y) & = & (-1)^{\chi(yp)}\, \sgn(x)^{\chi(x)}\,\, |y/p|\,\, a_y(x,p) \nonumber
\\ a_p(x,y) & = &  \sgn(p)^{\chi(p)}\,\sgn(x)^{\chi(x)}\,\sgn(y)^{\chi(y)}\, a_p(y,x) \nonumber
\\ a_p(x,y) & = &  (-1)^{\chi(xp)} \, \,\sgn(y)^{\chi(y)}\,|x/p|\, a_x(p,y) \ . \nonumber
\end{eqnarray}
\end{proposition}

\bigskip

Now we produce the eigenvectors of our self-adjoint extension of $\de_0(\gamma_0^\dag \gamma_0)$ (see the remarks after
the proof of proposition \ref{comult.prop1} ). We will use these eigenvectors to define a unitary operator that will
induce the comultiplication. The presence of $\lambda^{\chi(y)}$ and $\mu^{\chi(x)}$ in the formulas for
$\digamma_{r,s,m,p}$ can be traced back to Eq.~(\ref{comult.eq3}). The dependence of $\digamma_{r,s,m,p}$ on $r$,$s$
and $p$ is chosen in such a way that Proposition \ref{comult.prop1} is true.

\smallskip

\begin{definition} \label{comult.def2}
Consider $r,s \in \Z$, $m \in \Z$ and $p \in I_q$. We define the element $\digamma_{r,s,m,p} \in H \ot H$ such that
$$\digamma_{r,s,m,p}(\lambda,x,\mu,y) = \begin{cases} a_p(x,y)\, \lambda^{r+\chi(y/p)}\,\mu^{s-\chi(x/p)} & \text{ if } y =
\sgn(p)\,q^m\,x \\ \hspace{13ex} 0 & \text{ otherwise}\end{cases}$$ for all $x,y \in I_q$ and $\lambda,\mu \in \T$.
\end{definition}

\medskip

We know from Proposition \ref{comult.ortho1} that the family $(\,\digamma_{r,s,m,p} \mid r,s \in \Z, m \in \Z, p \in
I_q\,)$ forms an orthonormal basis for $H \ot H$. As a consequence, we get the following result.

\smallskip

\begin{proposition} \label{comult.prop3}
If $x,y \in I_q$ and $r,s \in \Z$, then, using $L^2$-convergence,
\begin{eqnarray*}
\ze^r \ot \sde_x \ot \ze^s \ot \sde_y = \sum_{p \in I_q} \, a_p(x,y)\,\,\digamma_{r-\chi(y/p),s+\chi(x/p),\chi(y/x),p}
 \ .
\end{eqnarray*}
\end{proposition}

\medskip

Now we are ready to introduce the comultiplication of quantum $\widetilde{SU}_q(1,1)$.

\begin{proposition} \label{comult.prop7}
Define the unitary transformation $V : H \ot H \rightarrow L^2(\T) \ot L^2(\T) \ot H$ such that $V(\digamma_{r,s,m,p})
= \ze^r \ot \ze^s \ot \ze^m \ot \sde_p$ for all $r,s \in \Z$, $m \in \Z$ and $p \in I_q$. Then there exists a unique
injective normal $^*$-homomorphism $\de : M_q \rightarrow M_q \ot M_q$ such that $\de(a) = V^* (1_{L^2(\T)} \ot
1_{L^2(\T)} \ot a) V$ for all $a \in M_q$.
\end{proposition}
\begin{proof}
Define the $^*$-homomorphism $\de : M_q \rightarrow B(H \ot H)$ such that $\de(a) = V^* (1_{L^2(\T)} \ot 1_{L^2(\T)}
\ot a) V$ for all $a \in M_q$. Fix $a \in M_q$. If $r,s,m \in \Z$, $p \in I_q$, then Definition \ref{comult.def2}
implies that $(1_H \ot M_\ze \ot 1_{L^2(I_q)})\digamma_{r,s,m,p} = \digamma_{r,s+1,m,p}$. It follows that $V (1_H \ot
M_\ze \ot 1_{L^2(I_q)}) = (1_{L^2(\T)} \ot M_\ze \ot 1_H)V$. As a consequence, $V (1_H \ot b \ot 1_{L^2(I_q)}) =
(1_{L^2(\T)} \ot b \ot 1_H)V$ for all $b \in L^\infty(\T)$. It follows that $\de(a)(1_H \ot b \ot 1_{L^2(I_q)}) = (1_H
\ot b \ot 1_{L^2(I_q)})\de(a)$ for all $b \in L^\infty(\T)$. Thus, $\de(a) \in (1_H \ot L^\infty(\T) \ot 1_{L^2(I_q)})'
= B(H) \ot L^\infty(\T) \ot B(L^2(I_q)) = B(H) \ot M_q$ by Lemma \ref{neumann.lem2}. In a similar way, $\de(a) \in  M_q
\ot B(H)$. Hence $\de(a) \in  M_q \ot M_q$ by the commutator theorem for tensor products.
\end{proof}

\medskip

The requirement that $\de(M_q) \subseteq M_q \ot M_q$  is the primary reason for introducing the extra generator $u$.
We cannot  work with the von Neumann algebra $M_q^\diamond$ that is generated only by $\al$, $\gamma$ and $e$, because
$\de(M_q^\diamond) \not\subseteq M_q^\diamond \ot M_q^\diamond$.

\medskip

This definition of $\de$ and Eq.~(\ref{neumann.eq1}) imply easily that $\langle\,\digamma_{r,s,m,p} \mid r,s \in \Z, m
\in \Z, p \in I_q\,\rangle$ is a core for $\de(\al)$, $\de(\gamma)$ and
\begin{equation} \label{comult.eq2}
\begin{array}{lcl}
\de(\al)\,\digamma_{r,s,m,p} & = & \sqrt{\sgn(p)+p^{-2}}\,\,\digamma_{r,s,m,pq} \vspace{1.5ex}
\\ \de(\gamma)\,\digamma_{r,s,m,p} & = & p^{-1}\,\,\digamma_{r,s,m+1,p} \ .
\end{array}
\end{equation} for $r,s \in \Z$, $m \in \Z$ and $p \in I_q$.

\medskip

Since $\digamma_{r,s,m,p}(\lambda,x,\mu,y) = 0$ if $\sgn(x)\,\sgn(y) \not= \sgn(p)$ it follows that $(e \ot
e)F_{r,s,m,p} = \sgn(p) \, F_{r,s,m,p}$. Hence $V(e \ot e) = (1_{L^2(\T^2)} \ot e) V$. As a consequence, $\de(e) = e
\ot e$.

\medskip

Recall the linear operators $\de_0(\al_0)$, $\de_0(\gamma_0)$ acting on $E \odot E$ (Eqs.~(\ref{hopf.comult2})). Also
recall the distinction between $^*$ and $\dag$. Since $\de(\al) = V^*(1_{L^2(\T)} \ot 1_{L^2(\T)} \ot \al)V$ and
$\de(\gamma) = V^*(1_{L^2(\T)} \ot 1_{L^2(\T)} \ot \gamma)V$, Lemma \ref{neumann.lem1} implies that $\de(\al)$ and
$\de(\gamma)$ are balanced. For the proof of the next proposition $q$-contiguous  relations for $q$-hypergeometric
functions are essential (see Lemma \ref{app.contiguous}).

\smallskip

\begin{proposition} \label{comult.prop1}
The following inclusions hold: $\de_0(\al_0) \subseteq \de(\al)$, $\de_0(\al_0)^\dag \subseteq \de(\al)^*$,
$\de_0(\gamma_0) \subseteq \de(\gamma)$ and $\de_0(\gamma_0)^\dag \subseteq \de(\gamma)^*$.
\end{proposition}
\smallskip
\begin{trivlist}
\item[\ \,{\it Proof.} (1)] The first step is to prove for $i,j,m \in \Z$, $p \in I_q$ and $r,s \in \Z$, $x,y \in I_p$,
\begin{equation}
\langle \de(\al) \digamma_{i,j,m,p} , \ze^r \ot \sde_x \ot \ze^s \ot \sde_y \rangle = \langle \digamma_{i,j,m,p} ,
\de_0(\al_0^\dag)(\ze^r \ot \sde_x \ot \ze^s \ot \sde_y) \rangle \ . \label{short.eq1}
\end{equation}
and this boils down to the $q$-contiguous relations in Lemma \ref{app.contiguous}.  First we can rewrite
Eq.~(\ref{app.eq1}) in terms of the functions $a_p$ of Definition \ref{comult.def1} by a straightforward verification
as
\begin{eqnarray}
& & \sqrt{\sgn(p) + p^{-2}}\, a_{qp}(x,y) \nonumber
\\ & & \spat = \sqrt{(\sgn(x)+q^2/x^2)(\sgn(y)+q^2/y^2)}\,a_p(q^{-1}x,q^{-1}y)   + \sgn(x)\,(q/xy)\, a_p(x,y) \ .
\label{comult.eq6}
\end{eqnarray}
Here we apply Eq.~(\ref{app.eq1}) in base $q^2$ with $a$,$b$,$z$ replaced by $-q^2/\kappa(y)$, $q^2 \kappa(x/y)$,
$\kappa(x/p)$ respectively (and take some care if $x=-q$ or $y=-q$).

\smallskip

Now consider the left hand side of Eq.~(\ref{short.eq1}). By Eq.~(\ref{comult.eq2}) and Definition \ref{comult.def2},
the inner product on the left hand side of Eq.~(\ref{short.eq1}) is zero unless $|y/x| = q^m$, $\sgn(p) = \sgn(xy)$,
$i+\chi(y/qp) = r$ and $j-\chi(x/qp) = s$ in which case, this inner product equals $\sqrt{\sgn(p) + p^{-2}}\,
a_{qp}(x,y)$. \inlabel{short.eq3}

For the right hand side of Eq.~(\ref{short.eq1}) we recall that $\de_0(\al_0^\dag) \in \cL^+(E \odot E)$ is given by
$\de_0(\al_0^\dag) = q\, e_0 \gamma_0 \odot \gamma_0^\dag + \al_0^\dag \odot \al_0^\dag$. From Eqs.~(\ref{hopf.eq1})
and (\ref{hopf.eq2}), it follows that the right hand side of Eq.~(\ref{short.eq1}) equals
\begin{eqnarray}
& & \sgn(x)\,(q/xy)\,\, \langle \digamma_{i,j,m,p} , \ze^{r+1} \ot \sde_x \ot \ze^{s-1} \ot \sde_y \rangle \nonumber
\\ & & \spat +
\sqrt{(\sgn(x)+q^2/x^2)(\sgn(y)+q^2/y^2)}\,\, \langle \digamma_{i,j,m,p} , \ze^r \ot \sde_{q^{-1} x} \ot \ze^s \ot
\sde_{q^{-1} y} \rangle\ . \label{short.eq2}
\end{eqnarray}
By Definition \ref{comult.def2}, this implies that the right hand side of Eq.~(\ref{short.eq1}) is zero unless $|y/x| =
q^m$, $\sgn(p) = \sgn(xy)$, $i+\chi(y/qp) = r$ and $j-\chi(x/qp) = s$, which proves Eq.~(\ref{short.eq1}) if these
conditions are violated. If, on the other hand, these conditions are satisfied, Definition \ref{comult.def2} guarantees
that (\ref{short.eq2}) equals
$$ \sgn(x)\,(q/xy)\, a_p(x,y) + \sqrt{(\sgn(x)+q^2/x^2)(\sgn(y)+q^2/y^2)}\,a_p(q^{-1}x,q^{-1}y)    \ .$$
Therefore (\ref{short.eq3}) and Eq.~(\ref{comult.eq6}) imply Eq.~(\ref{short.eq1})  also in this case.

\smallskip

So we have established Eq.~(\ref{short.eq1}) in all possible cases. Since the elements $\digamma_{i,j,m,p}$ form a core
for $\de(\al)$, we conclude that $\ze^r \ot \sde_x \ot \ze^s \ot \sde_y \in D(\de(\al)^*)$ and $\de(\al)^*(\ze^r \ot
\sde_x \ot \ze^s \ot \sde_y) = \de_0(\al_0^\dag)(\ze^r \ot \sde_x \ot \ze^s \ot \sde_y)$, thus proving that
$\de_0(\al_0^\dag) \subseteq \de(\al)^*$. Taking the adjoint of this equation, we see that $\de(\al) \subseteq
\de_0(\al_0^\dag)^*$. In general, $\de_0(\al_0) \subseteq \de_0(\al_0^\dag)^*$. Since $D(\de(\al)) = D(\de(\al)^*)
\supseteq E \odot E$, we conclude that the  inclusion $\de_0(\al_0) \subseteq \de(\al)$ also holds.

\medskip

\item[\ \,(2)] The inclusions regarding $\de(\gamma)$ and $\de(\gamma)^*$ are treated in the same way but this time one
applies Eq.~(\ref{app.eq2}) in base $q^2$ with $a$,$b$,$z$ replaced by $-q^2/\kappa(y)$, $q^2 \kappa(x/y)$, $q^2
\kappa(x/p)$ respectively, to obtain
\begin{eqnarray*}
& & p^{-1}\,a_p(x,y) = \sgn(x)\, y^{-1}\,\sqrt{\sgn(x)+1/x^2}\,\, a_p(q x,y)  + x^{-1}\,\sqrt{\sgn(y)+ q^2/y^2}\,\,
a_p(x,q^{-1}y)  \
\end{eqnarray*}
for all $x,y,p \in I_q$. \qed
\end{trivlist}

\medskip

This proposition implies also that $\de(\gamma^* \gamma)$ is an extension of $\de_0(\gamma_0^\dag \gamma_0)$. We also
know that $\langle\,\digamma_{r,s,m,p} \mid r,s \in \Z, m \in \Z, p \in I_q\,\rangle$ is a core for $\de(\gamma^*
\gamma)$ and $\de(\gamma^* \gamma)\,\digamma_{r,s,m,p}  =  p^{-2}\,\,\digamma_{r,s,m,p}$ for $r,s,m \in \Z$, $p \in
I_q$.  Using this information it is not so difficult to check that $U_\theta^* \de(\gamma^*
\gamma)\!\!\restriction_{\ell_\theta'} U_\theta = 1 \ot L^{\sgn(\theta)}_\theta$ for all $\theta \in \qz$, but we will
not make any use of this fact in  this paper.

\medskip

In the next step we  investigate the behavior of  $V$ with respect to the flip map $\Si : H \ot H \rightarrow H \ot H$.
Later on, this will guarantee that the unitary antipode (see the discussion after the proof of
\cite[Prop.~1.4]{VaKust2}) commutes with the comultiplication up to the flip map. It will also reduce some of the
calculations in the proof of the coassociativity. This behavior is directly related to the second symmetry property of
Proposition \ref{comult.prop2}.

\smallskip

Let us  define anti-unitary permutation operators $\Si' : L^2(\T^2) \rightarrow L^2(\T^2)$ and
 $\tilde{\Si} : L^2(\T^4) \rightarrow L^2(\T^4)$ such that $\Si'(v_1 \ot v_2) = \bar{v}_2 \ot \bar{v}_1$ and $\tilde{\Si}(v_1
\ot v_2 \ot v_3 \ot v_4) = \bar{v}_4 \ot \bar{v}_3 \ot \bar{v}_2 \ot \bar{v}_1$ for all $v_1,v_2,v_3,v_4 \in L^2(\T)$.

\smallskip

We also introduce the anti-unitary operator $\breve{J}$ on $H$ such that $\breve{J}\,[f] = [g]$, where $f,g \in
\cL^2(\T \times I_q)$ are such that $g(\lambda,x) = \sgn(x)^{\chi(x)}$ $\overline{f(\lambda,x)}$ for all  $x \in I_q$
and $\lambda \in \T$.

\begin{proposition} \label{coass.prop1}
We have that $V \Si = (\Si' \ot \breve{J}) V (\breve{J} \ot \breve{J})$ and
$$V_{13} (1_H \ot V) \Si_{13} = (\tilde{\Si} \ot  \breve{J}) (1_{L^2(\T^2)}  \ot V) (V \ot
1_H) (\breve{J} \ot \breve{J} \ot \breve{J})\ .$$
\end{proposition}
\begin{proof} 
Take $r,s,m \in \Z$ and $p \in I_q$. Proposition \ref{comult.prop2} implies that $a_p(x,y) =
\sgn(p)^{\chi(p)}\,\sgn(x)^{\chi(x)}\,\sgn(y)^{\chi(y)}$ $a_p(y,x)$ for all $x,y \in I_q$. If $x,y \in I_q$ and
$\lambda,\mu \in \T$, then Definition \ref{comult.def2} implies that
\begin{eqnarray*}
& & \hspace{-5ex} (\Si(\breve{J} \ot \breve{J})\digamma_{r,s,m,p})(\lambda,x,\mu,y) = ((\breve{J} \ot \breve{J})
\digamma_{r,s,m,p})(\mu,y,\lambda,x) = \sgn(x)^{\chi(x)}\,\sgn(y)^{\chi(y)}\,\,\overline{F_{r,s,m,p}(\mu,y,\lambda,x)}
\\ & &  \hspace{-5ex}  \spat  =  \sde_{x,\sgn(p)y q^m} \, \sgn(x)^{\chi(x)}\,\sgn(y)^{\chi(y)}\,\, \overline{\mu^{r+\chi(x/p)} \,\lambda^{s-\chi(y/p)}\,a_p(y,x)}
\\ & &  \hspace{-5ex} \spat = \sde_{y,\sgn(p)x q^{-m}}\, \sgn(p)^{\chi(p)}\, \lambda^{-s+\chi(y/p)}\,\mu^{-r-\chi(x/p)}\,a_p(x,y)
=  \sgn(p)^{\chi(p)}\, \digamma_{-s,-r,-m,p}(\lambda,x,\mu,y) \ .
\end{eqnarray*}
Thus,
\begin{eqnarray*}
V \Si (\breve{J} \ot \breve{J})\digamma_{r,s,m,p}  & = & \sgn(p)^{\chi(p)}\, V \digamma_{-s,-r,-m,p} =
\sgn(p)^{\chi(p)}\, \ze^{-s} \ot \ze^{-r} \ot \ze^{-m} \ot \sde_p
\\  & = &  (\Si' \ot \breve{J})(\ze^r \ot \ze^s \ot \ze^m \ot \sde_p)
= (\Si' \ot  \breve{J})V\digamma_{r,s,m,p}
\end{eqnarray*}
From this all we conclude that $V \Si = (\Si' \ot  \breve{J}) V (\breve{J} \ot \breve{J})$. Notice that $\Si_{13} =
\Si_{23} \Si_{12} \Si_{23}$. Thus
\begin{eqnarray*}
& & V_{13} (1_H \ot V) \Si_{13}  = V_{13} (1_H \ot V) \Si_{23} \Si_{12} \Si_{23} = V_{13}(1_H \ot  [(\Si' \ot
\breve{J}) V (\breve{J} \ot \breve{J})]) \Si_{12} \Si_{23} \ .
\end{eqnarray*}
Let $\hat{\Si} : L^2(\T^2) \ot L^2(\T^2) \rightarrow L^2(\T^2) \ot L^2(\T^2)$ denote the flip map. If we let $\Xi :
L^2(\T^2) \ot H \ot H \rightarrow H \ot L^2(\T^2) \ot H$ denote the permutation map defined by $\Xi(u \ot v \ot w) = w
\ot u \ot v$ for all $u \in L^2(\T^2)$ and $v,w \in H$, we get that
\begin{eqnarray*}
& & V_{13} (1_H \ot V) \Si_{13}  = V_{13} \, \Xi\, ([(\Si' \ot \breve{J}) V (\breve{J} \ot \breve{J})] \ot 1_H)
\\ & & \spat = (\hat{\Si} \ot 1_H)(1_{L^2(\T^2)} \ot V \Si) \,([(\Si' \ot \breve{J}) V (\breve{J} \ot \breve{J})] \ot 1_H)
\\ & & \spat = (\hat{\Si} \ot 1_H)(1_{L^2(\T^2)}
\ot [(\Si' \ot  \breve{J}) V (\breve{J} \ot \breve{J})]) \, ([(\Si' \ot  \breve{J}) V (\breve{J} \ot \breve{J})] \ot
1_H)
\\ & & \spat = (\hat{\Si}(\Si' \ot \Si') \ot \breve{J})(1_{L^2(\T^2)}
\ot  V ) \, ( V \ot 1_H) (\breve{J} \ot \breve{J} \ot \breve{J})
\\ & & \spat = (\tilde{\Si} \ot \breve{J}) (1_{L^2(\T^2)} \ot V) (V \ot 1_H) (\breve{J} \ot \breve{J} \ot \breve{J})\
.
\end{eqnarray*}
\end{proof}

\bigskip

We end this section with the statement that the comultiplication is coassociative, the result of this paper that is the
most technical to prove.

\smallskip

\begin{theorem} \label{coass.thm1}
The $^*$-homomorphism $\de : M_q \rightarrow M_q \ot M_q$ is coassociative, i.e.~$(\de \ot \io)\de = (\io \ot \de)\de$.
\end{theorem}

\smallskip

In order to enhance the readability of this paper, the proof will be given in section 5 but we stress that the results
of section 4 do not play any role therein.

\bigskip\medskip

\sectie{The Haar weight on $(M_q,\de)$.}

\bigskip

In this section we construct the left Haar weight on $(M_q,\de)$ and prove its left invariance through the construction
of a partial isometry that afterwards turns out to be the multiplicative unitary of quantum $\widetilde{SU}(1,1)$. We
also establish the unimodularity of $(M_q,\de)$.

\smallskip

Since $M_q = L^\infty(\T) \ot B(L^2(I_q))$ we can consider the trace $\Tr$ on $M_q$ given by $\Tr = \Tr_{L^\infty(\T)}
\ot \Tr_{B(L^2(I_q))}$, where $\Tr_{L^\infty(\T)}$ and $\Tr_{B(L^2(I_q))}$ are the canonical traces on $L^\infty(\T)$
and $B(L^2(I_q))$ which we choose to be normalized in such a way that $\Tr_{L^\infty(\T)}(1) = 1$ and
$\Tr_{B(L^2(I_q))}(P) = 1$ for every rank one projection $P$ in $B(L^2(I_q))$.

\smallskip

Given a weight $\eta$ on $M_q$, we use the following standard concepts  from  weight theory:
$$\cM_\eta^+ = \{\,x \in M_q^+ \mid \eta(x) < \infty\,\}, \hspace{1cm} \cM_\eta = \text{linear span of } \cM_\eta^+,
\hspace{1cm}  \cN_\eta = \{\,x \in M_q \mid \eta(x^* x) < \infty\,\}\ .$$

Next we introduce a GNS-construction for the trace $\Tr$. Define $K = H \ot L^2(I_q) = L^2(\T) \ot L^2(I_q) \ot
L^2(I_q)$. If $m \in \Z$ and $p,t \in I_q$, we set $f_{m,p,t} = \ze^m \ot \sde_p \ot \sde_t \in K$. Thus $(\,f_{m,p,t}
\mid m \in \Z, p,t \in I_q\,)$ is an orthonormal basis for $K$. Now define
\begin{trivlist}
\item[\ \,(1)] a linear map $\la_{\Tr} : \cN_{\Tr} \rightarrow K$ such that $\la_{\Tr}(a) = \sum_{p \in I_q}
(a \ot 1_{L^2(I_q)}) f_{0,p,p}$ for $a \in \cN_{\Tr}$.
\smallskip
\item[\ \,(2)] a unital $^*$-homomorphism $\pi : M_q \rightarrow B(K)$ such that $\pi(a) =
a \ot 1_{L^2(I_q)}$ for all $a \in M_q$.
\smallskip
\end{trivlist}
Then $(K,\pi,\la_{\Tr})$ is a GNS-construction for $\Tr$.

\medskip

Now we are ready to define the weight that will turn out to be left- and right invariant with respect to $\de$.  Use
the remarks before \cite[Prop.~1.15]{VaKust} to define a linear map $\la =
(\la_{\Tr})_{\gamma^* \gamma} : D(\la)
\subseteq M_q \rightarrow K$.

\smallskip

\begin{definition}
We define the faithful normal semi-finite weight $\vfi$ on $M_q$ as $\vfi = \text{\rm Tr}_{\gamma^* \gamma}$. By
definition, $(K,\pi,\la)$ is a GNS-construction for $\vfi$.
\end{definition}

\smallskip

See \cite{Va2} for more details about this definition. This definition of $\vfi$ is of course compatible with the usual
construction of absolutely continuous weights (see \cite{Ped}). So we already know that the modular automorphism group
$\si^{\vfi}$ of $\vfi$ is such that $\si^{\vfi}_s(x) = |\gamma|^{2is}\, x \,|\gamma|^{-2is}$ for all $x \in M_q$ and $s
\in  \R$.

The modular conjugation of $\vfi$ with respect to $(K,\pi,\la)$ will be denoted by $J$, the modular operator of $\vfi$
will be denoted by $\nab$. Since $\vfi = \Tr_{\gamma^* \gamma}$ and $\la = (\la_{\Tr})_{\gamma^* \gamma}$, $J$ equals
the modular conjugation of $\Tr$ with respect to $(K,\pi,\la_{\Tr})$.

Let us recall how these modular objects are related to the modular group. So fix $a \in \Nfi$. If $t \in \R$, then
$\si^\vfi_t(a) \in \Nfi$ and $\nab^{it} \la(a) = \la(\si^\vfi_t(a))$. Provided $a \in D(\si^\vfi_{\frac{i}{2}})$, the
element $\si^\vfi_{\frac{i}{2}}(a)^* \in \Nfi$ and $J \la(a) =  \la\bigl(\si^\vfi_{\frac{i}{2}}(a)^*\bigr)$. For $x \in
D(\si^\vfi_{\frac{i}{2}})$, we have that $a x \in \Nfi$ and $\la(ax) = J \pi(\si^\vfi_{\frac{i}{2}}(x)^*) J \la(a)$.

\medskip

Let us first establish some easy formulas that make working with the weight $\vfi$ pretty easy.

\begin{lemma} \label{mult.lem1}
Consider $m \in \Z$ and $p,t \in I_q$. Then
\begin{trivlist}
\item[\ \,(1)] $\Phi(m,p,t) \in \cN_{\vfi}$ and $\la(\Phi(m,p,t)) = |t|^{-1} \, f_{m,p,t}$,
\item[\ \,(2)] $\Phi(m,p,t) \in \cM_{\vfi}$ and $\vfi(\Phi(m,p,t)) = \sde_{m,0} \, \sde_{p,t}\, t^{-2}$,
\item[\ \,(3)] $\Phi(m,p,t)$ is analytic with respect to $\si^{\vfi}$ and $\si^{\vfi}_z(\Phi(m,p,t)) =
|p^{-1} t|^{2iz} \, \Phi(m,p,t)$ for all $z \in \C$.
\end{trivlist}
\end{lemma}
\smallskip
\begin{trivlist}
\item[\ \,{\it Proof.} (1)] By Eq.~(\ref{neumann.eq4}), $\Phi(m,p,t)^* \Phi(m,p,t) = \Phi(0,t,t) = M_{1 \ot \sde_t}$ which clearly belongs to $\cM_{\Tr}^+$. Thus $\Phi(m,p,t) \in \cN_{\Tr}$ and
\begin{equation} \label{mult.eq6}
\la_{\Tr}(\Phi(m,p,t)) = \sum_{p' \in I_q} (\Phi(m,p,t) \ot 1) f_{0,p',p'} = f_{m,p,t} \ .
\end{equation}
It is clear that $\Phi(m,p,t)\,|\gamma|$ is bounded and that the extension of this element to an an element of $B(H)$
is given by $|t|^{-1}\,\Phi(m,p,t) \in \cN_{\Tr}$. Therefore the remarks before  \cite[Prop.~1.15]{VaKust} imply that
$\Phi(m,p,t) \in \cN_{\vfi}$ and $\la(\Phi(p,m,t)) = |t|^{-1}\, \la_{\Tr}(\Phi(m,p,t)) = |t|^{-1}\,f_{m,p,t}$.
\smallskip
\item[\ \,(2)] By Eq.~(\ref{neumann.eq4}),  $\Phi(m,p,t) = \Phi(0,p,p)^* \Phi(m,p,t)$, thus
$\vfi(\Phi(m,p,t)) = \langle \la(\Phi(m,p,t)) , \la(\Phi(0,p,p))\rangle$. Hence (2) follows from (1).
\smallskip
\item[\ \,(3)] Choose $s \in \R$. By definition, $\Phi(m,p,t) = \rho_{p^{-1} t} \, M_{\ze^m \ot \sde_t}$. Thus Lemma
\ref{neumann.lem3} implies that
\begin{eqnarray*}
& & \hspace{-2ex} |\gamma|^{2is}\,\Phi(m,p,t) =  M_{1 \ot |\xi|^{-2is}}\,\rho_{p^{-1} t}\, M_{\ze^m \ot \sde_t} =
\rho_{p^{-1} t}\, M_{T_{p t^{-1}}(1 \ot |\xi|^{-2is})} \, M_{\ze^m \ot \sde_t}
\\ &  &  = |p^{-1}
t|^{2is}\, \rho_{p^{-1} t}\, M_{1 \ot |\xi|^{-2is}} \, M_{\ze^m \ot \sde_t} = |p^{-1} t|^{2is}\, \rho_{p^{-1} t}\,
M_{\ze^m \ot \sde_t}\,M_{1 \ot |\xi|^{-2is}} = |p^{-1} t|^{2is}\, \Phi(m,p,t) \, |\gamma|^{2is} \ ,
\end{eqnarray*}
from which it follows that $\si^{\vfi}_s(\Phi(m,p,t)) = |p^{-1} t|^{2is}\, \Phi(m,p,t)$ and the claim follows. \qed
\end{trivlist}

\smallskip

These results imply easily the next lemma which is needed to extend results (involving $\la$) that have been proven for
the elements $\Phi(m,p,t)$ to the whole of $\Nfi$.

\smallskip

\begin{corollary} \label{mult.lem2}
If $a \in \Nfi$, there exists a bounded net $(a_i)_{i \in I}$ in $\langle \, \Phi(m,p,t) \mid m \in \Z, p,t \in
I_q\,\rangle$ such that $(a_i)_{i \in I}$ converges strongly$^*$ to $a$ and $(\la(a_i))_{i \in I}$ converges to
$\la(a)$.
\end{corollary}
\begin{proof}
Define $C = \langle \, \Phi(m,p,t) \mid m \in \Z, p,t \in I_q\,\rangle$. Set $B = \{\,b \in M_q \mid \|b\| \leq
\|a\|\,\}$ and equip $B$ with the strong$^*$ topology. In $B \times K$, we consider the set $G$ that is the closure of
$\{\,(b,\la(b)) \mid b \in B \cap C\,\}$.

Let $F(I_q)$ denote the set of all finite subsets of $I_q$ and turn $F(I_q)$ into a directed set by inclusion. For $L
\in F(I_q)$, we define the projection $P_L = \sum_{v \in L} \Phi(0,v,v)$. Thus, $(P_L)_{L \in F(I_q)}$ converges
strongly to 1.

Fix $L \in F(I_q)$ for the moment. By Lemma \ref{mult.lem1}(1) we see that $P_L \in \Nfi$. Since $\la(b P_L) =
\pi(b)\,\la(P_L)$ for all $b \in M_q$, Kaplansky's density Theorem (see the proof of \cite[Thm.~5.3.5]{Kad}), applied
to the $^*$-algebra $C$, implies that $(a P_L,\la(a P_L)) \in G$.

It is clear from Lemma $\ref{mult.lem1}(3)$ that $P_L \in D(\si_{\frac{i}{2}}^\vfi)$ and $\si_{\frac{i}{2}}^\vfi(P_L) =
P_L$. Thus, $\la(a P_L) = J \pi(P_L) J \la(a)$. It follows that $\bigl(\,(a P_L,\la(a P_L))\,\bigr)_{L \in F(I_q)}
\rightarrow (a,\la(a))$, implying that $(a,\la(a)) \in G$ and the lemma follows.
\end{proof}

\smallskip

\begin{corollary} \label{mult.cor1}
Consider $m \in \Z$ and $p,t \in I_q$. Then
\begin{trivlist}
\item[\ \,1)] $J f_{m,p,t} = f_{-m,t,p}$,
\item[\ \,2)] $f_{m,p,t} \in D(\nab)$ and $\nab f_{m,p,t} = |p^{-1} t|^2\, f_{m,p,t}$.
\end{trivlist}
\end{corollary}
\begin{proof}
Recall that $J$ is also the modular conjugation of the trace $\Tr$. So Eq.~(\ref{mult.eq6}) in the previous proof
implies that $J f_{m,p,t} = J \la_{\Tr}(\Phi(m,p,t)) = \la_{\Tr}(\Phi(m,p,t)^*) = \la_{\Tr}(\Phi(-m,t,p)) =
f_{-m,t,p}$.

Let $s \in \R$. By Lemma \ref{mult.lem1}, $\Phi(m,p,t) \in \cN_{\vfi}$ and $\si^\vfi_s(\Phi(m,p,t)) = |p^{-1}
t|^{2is}\, \Phi(m,p,t)$. So we get that $\nab^{is} \la(\Phi(m,p,t)) = \la\bigr(\si^\vfi_s(\Phi(m,p,t))\bigl) = |p^{-1}
t|^{2is} \, \la(\Phi(m,p,t))$. Thus, since $f_{m,p,t} = |t|\, \la(\Phi(m,p,t))$, $\nab^{is} f_{m,p,t} = |p^{-1}
t|^{2is} \, f_{m,p,t}$ and also the second claim follows.
\end{proof}

\medskip

Now we know enough about the weight $\vfi$ to proceed to the proof of the left invariance of $\vfi$ with respect to
$\de$. For this purpose we introduce a partial isometry that will later turn out to be the multiplicative unitary of
$(M_q,\de)$. The formula defining $W^*$ in the next proposition is obtained by formally calculating $W^*
\bigl(\la(\Phi(m_1,p_1,t_1)) \ot \la(\Phi(m_2,p_2,t_2))\bigr) = (\la \ot \la)(\de(\Phi(m_2,p_2,t_2))(\Phi(m_1,p_1,t_1)
\ot 1))$, but as this stage we do not know whether $\vfi$ is left invariant so we cannot use this  formula, as is done
in the general theory, to define  $W^*$.

\smallskip

\begin{proposition} \label{mult.prop1}
There exists a unique surjective partial isometry $W$ on $K \ot K$ such that
\begin{eqnarray*}
& & W^* (f_{m_1,p_1,t_1} \ot f_{m_2,p_2,t_2}) =   \sum_{\scriptsize \begin{array}{c} y,z \in I_q \\ \sgn(p_2
t_2)(yz/p_1)q^{m_2} \in  I_q
\end{array}} \ |t_2/y|\,a_{t_2}(p_1,y)\,a_{p_2}(z,\sgn(p_2 t_2) (yz/p_1) q^{m_2})
\\ & & \hspace{39ex} \times \ \ f_{m_1+m_2- \chi(p_1 p_2/t_2 z),z,t_1} \ot f_{\chi(p_1 p_2/t_2 z), \sgn(p_2 t_2)
(yz/p_1) q^{m_2} , y} \nonumber
\end{eqnarray*}
for all $m_1,m_2 \in \Z$ and $p_1,p_2,t_1,t_2 \in I_q$.
\end{proposition}
\begin{proof}
The proof of this proposition is based on the properties of the functions $a_p$, in particular on Propositions
\ref{comult.ortho1}, \ref{comult.ortho2} and \ref{comult.prop2}. We now give details.

Fix $m_1,m_2 \in \Z$ and $p_1,p_2,t_1,t_2 \in I_q$. If $y,z,y',z' \in I_q$ such that $\sgn(p_2 t_2)(yz/p_1)q^{m_2}$,
\newline $\sgn(p_2 t_2)(y'z'/p_1)q^{m_2} \in I_q$ and $(y,z) \not= (y',z')$, then it is clear that the element $f_{m_1+m_2- \chi(p_1
p_2/t_2 z),z,t_1}$ $\ot f_{\chi(p_1 p_2/t_2 z), \sgn(p_2 t_2) (yz/p_1) q^{m_2} , y}$ is orthogonal to the element
$f_{m_1+m_2- \chi(p_1 p_2/t_2 z'),z',t_1} \ot$ \newline $f_{\chi(p_1 p_2/t_2 z'), \sgn(p_2 t_2) (y'z'/p_1) q^{m_2} ,
y'}$. By Propositions \ref{comult.ortho1}, \ref{comult.ortho2} and \ref{comult.prop2} we get moreover
\begin{eqnarray}
& & \hspace{-10ex} \sum_{\scriptsize \begin{array}{c} y,z \in I_q \\ \sgn(p_2 t_2)(yz/p_1)q^{m_2} \in  I_q
\end{array}}\ \ |t_2/y|^2\,|a_{t_2}(p_1,y)|^2\,|a_{p_2}(z,\sgn(p_2 t_2) (yz/p_1) q^{m_2})|^2\, \nonumber
\\ & & \spat = \sum_{y \in I_q} \, |t_2/y|^2\,|a_{t_2}(p_1,y)|^2\,\,\sum_{\scriptsize \begin{array}{c} z \in I_q \\
\sgn(p_2 t_2)(yz/p_1)q^{m_2} \in  I_q \end{array}}\ \ |a_{p_2}(z,\sgn(p_2 t_2) (yz/p_1) q^{m_2})|^2\, \nonumber
\\ & & \spat \leq \sum_{y \in I_q} \, |t_2/y|^2\,|a_{t_2}(p_1,y)|^2 = \sum_{y \in I_q} \, |a_{y}(p_2,t_1)|^2 = 1 \ .
\label{mult.eq1}
\end{eqnarray}
From this, it follows that we can define the element $v(m_1,p_1,t_1;m_2,p_2,t_2) \in K \ot K$ as
\begin{eqnarray*}
& &  \hspace{-10ex} \sum_{\scriptsize \begin{array}{c} y,z \in I_q \\ \sgn(p_2 t_2)(yz/p_1)q^{m_2} \in  I_q
\end{array}} \ |t_2/y|\,a_{t_2}(p_1,y)\,a_{p_2}(z,\sgn(p_2 t_2) (yz/p_1) q^{m_2})
\\ & & \hspace{10ex} \times \ \ f_{m_1+m_2- \chi(p_1 p_2/t_2 z),z,t_1} \ot f_{\chi(p_1 p_2/t_2 z), \sgn(p_2 t_2) (yz/p_1) q^{m_2} , y} \ .\nonumber
\end{eqnarray*}
Note that $v(m_1,p_1,t_1;m_2,p_2,t_2)$ is just the right hand side of the formula we want to define $W^*$ by. Now
choose also $m_1',m_2' \in \Z$ and $p_1',p_2',t_1',t_2' \in I_q$. Notice that inequality (\ref{mult.eq1}) together with
the Cauchy-Schwarz equality implies that
\begin{eqnarray*}
& &  \hspace{-6ex} \sum_{\scriptsize \begin{array}{c} y,z \in I_q \\ \sgn(p_2 t_2)(yz/p_1)q^{m_2} \in  I_q
\\ \sgn(p_2' t_2')(yz/p_1')q^{m_2'} \in  I_q \end{array}}\ \ |\,(t_2/y)\, (t_2'/y) \,a_{t_2}(p_1,y)\,a_{t_2'}(p_1',y) \, \begin{tabular}[t]{l} $a_{p_2}(z,\sgn(p_2 t_2) (yz/p_1)
q^{m_2})$ \\ \hspace{15ex} $a_{p_2'}(z,\sgn(p_2' t_2') (yz/p_1') q^{m_2'})\,|$ \end{tabular}
\end{eqnarray*}
is finite, which allows us to compute the next sum in any order we want:
\begin{eqnarray*}
& & \langle v(m_1,p_1,t_1;m_2,p_2,t_2) , v(m_1',p_1',t_1';m_2',p_2',t_2') \rangle
\\ & &
\\ & & \spat = \sum_{\scriptsize \begin{array}{c} y,z \in I_q \\ \sgn(p_2 t_2)(yz/p_1)q^{m_2} \in  I_q
\\ \sgn(p_2' t_2')(yz/p_1')q^{m_2'} \in  I_q \end{array}}\ \ |t_2/y|\, |t_2'/y| \,a_{t_2}(p_1,y)\,a_{t_2'}(p_1',y) \, a_{p_2}(z,\sgn(p_2 t_2) (yz/p_1) q^{m_2})
\,
\\ & & \spat\spat\spat   \times \ \ a_{p_2'}(z,\sgn(p_2' t_2') (yz/p_1') q^{m_2'})\
\langle f_{m_1+m_2- \chi(p_1 p_2/t_2 z),z,t_1} , f_{m_1'+m_2'- \chi(p_1' p_2'/t_2' z),z,t_1'} \rangle\,\,
\\ & & \spat\spat\spat   \times \ \ \langle f_{\chi(p_1 p_2/t_2 z), \sgn(p_2 t_2) (yz/p_1) q^{m_2} , y} , f_{\chi(p_1' p_2'/t_2' z), \sgn(p_2' t_2') (yz/p_1')
q^{m_2'} , y} \rangle
\\ & &
\\ & & \spat = \sum_{\scriptsize \begin{array}{c} y,z \in I_q \\ \sgn(p_2 t_2)(yz/p_1)q^{m_2} \in  I_q
\\ \sgn(p_2' t_2')(yz/p_1')q^{m_2'} \in  I_q \end{array}}\ \ |t_2/y|\, |t_2'/y| \,a_{t_2}(p_1,y)\,a_{t_2'}(p_1',y) \,
a_{p_2}(z,\sgn(p_2 t_2) (yz/p_1) q^{m_2}) \,
\\ & & \spat\spat\spat   \times \ \ a_{p_2'}(z,\sgn(p_2' t_2') (yz/p_1') q^{m_2'})\
\sde_{m_1+m_2,m_1'+m_2'}\,\sde_{t_1,t_1'}\, \sde_{|p_1 p_2/t_2|,|p_1' p_2'/t_2'|}
\\ & & \spat\spat\spat   \times \ \ \,\sde_{\sgn(p_2
t_2)q^{m_2}/p_1,\sgn(p_2' t_2')q^{m_2'}/p_1'}
\\ & &
\\ & & \spat = \sde_{m_1+m_2,m_1'+m_2'}\,\sde_{t_1,t_1'}\, \sde_{p_1 p_2/t_2,p_1' p_2'/t_2'}\,\sde_{\sgn(p_2 t_2)q^{m_2}/p_1,\sgn(p_2' t_2')q^{m_2'}/p_1'}
\ \sum_{y \in I_q} \ |t_2/y|\, |t_2'/y| \,a_{t_2}(p_1,y)
\\ & & \spat\ \ \  a_{t_2'}(p_1',y)\ \sum_{\scriptsize \begin{array}{c} z \in I_q \\ \sgn(p_2 t_2)(yz/p_1)q^{m_2} \in  I_q
\end{array}}\ \ \, a_{p_2}(z,\sgn(p_2 t_2) (yz/p_1) q^{m_2}) \, a_{p_2'}(z,\sgn(p_2 t_2) (yz/p_1) q^{m_2}) \ .
\end{eqnarray*}
Thus, using the orthogonality relation in Proposition \ref{comult.ortho1} together with Proposition \ref{comult.prop2},
we get that
\begin{eqnarray*}
& & \langle v(m_1,p_1,t_1;m_2,p_2,t_2) , v(m_1',p_1',t_1';m_2',p_2',t_2') \rangle
 = \sde_{m_1+m_2,m_1'+m_2'}\,\sde_{t_1,t_1'}\, \sde_{p_1 p_2/t_2,p_1' p_2'/t_2'}
\\ & & \spat\spat \sde_{\sgn(p_2 t_2)q^{m_2}/p_1,\sgn(p_2'
t_2')q^{m_2'}/p_1'} \, \sum_{y \in I_q} \ |t_2/y|\, |t_2'/y| \,a_{t_2}(p_1,y)\,a_{t_2'}(p_1',y)\,
\sde_{p_2,p_2'}\,\sde_{\sgn(y),\sgn(p_1 t_2)}
\\ & & \spat = \sde_{m_1+m_2,m_1'+m_2'}\,\sde_{t_1,t_1'}\, \sde_{p_1/t_2,p_1'/t_2'}\,\sde_{\sgn(t_2)q^{m_2}/p_1,\sgn(
t_2')q^{m_2'}/p_1'}\, \sde_{p_2,p_2'}
\\ & & \spat\ \ \ \sum_{y \in I_q} \ (-1)^{\chi(y t_2)}\,\sgn(p_1)^{\chi(p_1)}\,(-1)^{\chi(y t_2')}\,\sgn(p_1')^{\chi(p_1')}\,a_{y}(p_1,t_2)\,a_{y}(p_1',t_2')
\\ & & \spat = \sde_{m_1+m_2,m_1'+m_2'}\,\sde_{t_1,t_1'}\, \sde_{p_1/t_2,p_1'/t_2'}\,\sde_{\sgn(t_2)q^{m_2}/p_1,\sgn(
t_2')q^{m_2'}/p_1'}\, \sde_{p_2,p_2'}
\\ & & \spat\ \ \ (-1)^{\chi(t_2 t_2')} \, \sgn(p_1)^{\chi(p_1)} \, \sgn(p_1')^{\chi(p_1')} \sum_{y \in I_q} \,
a_{y}(p_1,t_2)\,a_{y}(p_1',t_2') \ .
\end{eqnarray*}
Because of the factor $\sde_{p_1/t_2,p_1'/t_2'}$ in the above expression, Proposition \ref{comult.ortho2} implies that
\begin{eqnarray*}
& & \langle v(m_1,p_1,t_1;m_2,p_2,t_2) , v(m_1',p_1',t_1';m_2',p_2',t_2') \rangle
 = \sde_{m_1+m_2,m_1'+m_2'}\,\sde_{t_1,t_1'}
\\ & & \spat\ \ \ \sde_{p_1/t_2,p_1'/t_2'}\,\sde_{\sgn(t_2)q^{m_2}/p_1,\sgn(
t_2')q^{m_2'}/p_1'}\, \sde_{p_2,p_2'}\,(-1)^{\chi(t_2 t_2')} \, \sgn(p_1)^{\chi(p_1)} \, \sgn(p_1')^{\chi(p_1')} \,
\sde_{p_1,p_1'}\,\sde_{t_2,t_2'}
\\ & & \spat = \sde_{m_1,m_1'}\, \sde_{m_2,m_2'}\,\sde_{p_1,p_1'}\,\sde_{p_2,p_2'}\,\sde_{t_1,t_1'}\,\sde_{t_2,t_2'}
= \langle f_{m_1,p_1,t_1} \ot f_{m_2,p_2,t_2} , f_{m_1',p_1',t_1'} \ot f_{m_2',p_2',t_2'} \rangle \ .
\end{eqnarray*}
Now the lemma follows easily.
\end{proof}

\medskip

The next proposition deals with the essential step towards the left invariance of the weight $\vfi$ and the realization
that $W$ is the multiplicative unitary associated to quantum $\widetilde{SU}(1,1)$.

\smallskip

\begin{proposition} \label{mult.prop2}
Consider $\om \in B(K)_*$ and $a \in \cN_{\vfi}$. Then $(\om\pi \ot \io)\de(a) \in \cN_{\vfi}$ and $\la((\om\pi \ot
\io)\de(a)) = (\om \ot \io)(W^*)\,\la(a)$.
\end{proposition}
\begin{proof}
Choose $m,m_1,m_2 \in \Z$, $p,t,p_1,t_1,p_2,t_2 \in I_q$. Take also $v \in I_q$ and let us calculate the element
\newline $(\Phi(m_2,p_2,t_2)^* \ot 1)\de(\Phi(m,p,t))(\Phi(m_1,p_1,t_1) \ot \Phi(0,v,v))$.
Choose $r,s \in \Z$ and $x,y \in I_q$. Then Proposition \ref{comult.prop3},  Definition \ref{comult.def2} and
Eqs.~(\ref{neumann.eq3}), (\ref{neumann.eq4})  tell us that
\begin{eqnarray*}
& & (\Phi(m_2,p_2,t_2)^* \ot 1)\de(\Phi(m,p,t)) \, (\ze^r \ot \sde_x \ot \ze^s \ot \sde_y)
\\ & & \spat = \sum_{z \in I_q} \, a_z(x,y)\,(\Phi(-m_2,t_2,p_2) \ot
1)V^*(1 \ot \Phi(m,p,t))V\,\digamma_{r-\chi(y/z),s+\chi(x/z),\chi(y/x),z}
\\ & & \spat = \sum_{z \in I_q} a_z(x,y) \, (\Phi(-m_2,t_2,p_2) \ot
1)V^*(1 \ot \Phi(m,p,t))(\ze^{r-\chi(y/z)} \ot \ze^{s+\chi(x/z)} \ot \ze^{\chi(y/x)} \ot \sde_z)
\\ & & \spat = a_t(x,y)\,(\Phi(-m_2,t_2,p_2) \ot
1)V^* (\ze^{r-\chi(y/t)} \ot \ze^{s+\chi(x/t)} \ot \ze^{\chi(y/x)+m} \ot \sde_p)
\\ & & \spat = a_t(x,y)\,(\Phi(-m_2,t_2,p_2) \ot 1)\,\digamma_{r-\chi(y/t),s+\chi(x/t),\chi(y/x)+m,p}
\\ & & \spat = \sum_{\scriptsize \begin{array}{c} z \in I_q \\ \sgn(p) |y/x| q^m z \in I_q\end{array}}
a_t(x,y)\, a_p(z,\sgn(p) |y/x| q^m z )\,(\Phi(-m_2,t_2,p_2) \ot 1)\,
\\ & & \hspace{26ex} \ze^{r-\chi(y/t)+\chi(q^m z y/p x)} \ot \sde_z \ot \ze^{s+\chi(x/t)-\chi(z/p)} \ot \sde_{\sgn(p) |y/x| q^m z}
\\ & & \spat = \sum_{\scriptsize \begin{array}{c} z \in I_q \\ \sgn(pt) (y/x) q^m z \in I_q\end{array}}
a_t(x,y)\, a_p(z,\sgn(pt) (y/x) q^m z )\,(\Phi(-m_2,t_2,p_2) \ot 1)
\\ & & \hspace{27ex} \ze^{r+m-\chi(px/zt)} \ot \sde_z \ot \ze^{s+\chi(px/zt)} \ot \sde_{\sgn(pt) (y/x) q^m z}
\end{eqnarray*}
by Eq.~(\ref{neumann.eq3}). So we see that $(\Phi(m_2,p_2,t_2)^* \ot 1)\de(\Phi(m,p,t)) \, (\ze^r \ot \sde_x \ot \ze^s
\ot \sde_y) = 0$ if $\sgn(pt) (y/x) q^m p_2 \not\in I_q$. If, on the other hand, $\sgn(pt) (y/x) q^m p_2 \in I_q$, then
\begin{eqnarray*}
& & (\Phi(m_2,p_2,t_2)^* \ot 1)\de(\Phi(m,p,t)) \, (\ze^r \ot \sde_x \ot \ze^s \ot \sde_y)
\\ & & \spat = a_t(x,y)\, a_p(p_2,\sgn(pt) (y/x) q^m p_2 )\,(\ze^{r+m-m_2-\chi(px/p_2 t)} \ot \sde_{t_2} \ot \ze^{s+\chi(px/p_2t)} \ot
\sde_{\sgn(pt) (y/x) q^m p_2})\ .
\end{eqnarray*}
Now we have for $r,s \in \Z$ and $x,y \in I_q$ that
\begin{eqnarray*}
& & (\Phi(m_2,p_2,t_2)^* \ot 1)\de(\Phi(m,p,t))(\Phi(m_1,p_1,t_1) \ot \Phi(0,v,v))\, (\ze^r \ot \sde_x \ot \ze^s \ot
\sde_y)
\\ & & \spat = \sde_{t_1,x}\,\sde_{v,y} \, (\Phi(m_2,p_2,t_2)^* \ot 1)\de(\Phi(m,p,t)) (\ze^{m_1+r} \ot \sde_{p_1} \ot
\ze^s \ot \sde_v) \ .
\end{eqnarray*}
So we see that $(\Phi(m_2,p_2,t_2)^* \ot 1)\de(\Phi(m,p,t))(\Phi(m_1,p_1,t_1) \ot \Phi(0,v,v)) = 0$ if $\sgn(pt)
(v/p_1) q^m p_2 \not\in I_q$. If, on the other hand, $\sgn(pt) (v/p_1) q^m p_2 \in I_q$, then
\begin{eqnarray*}
& & (\Phi(m_2,p_2,t_2)^* \ot 1)\de(\Phi(m,p,t))(\Phi(m_1,p_1,t_1) \ot \Phi(0,v,v))\, (\ze^r \ot \sde_x \ot \ze^s \ot
\sde_y)
\\ & & \spat = \sde_{t_1,x}\,\sde_{v,y} \, a_t(p_1,v)\, a_p(p_2,\sgn(pt) (v/p_1) q^m p_2 )
\\ & & \spat\ \ \ \ze^{r+m_1+m-m_2-\chi(p p_1/p_2 t)} \ot \sde_{t_2} \ot \ze^{s+\chi(pp_1/p_2t)} \ot \sde_{\sgn(pt) (v/p_1) q^m
p_2} \ .
\end{eqnarray*}
and thus
\begin{eqnarray*}
& & (\Phi(m_2,p_2,t_2)^* \ot 1)\de(\Phi(m,p,t))(\Phi(m_1,p_1,t_1) \ot \Phi(0,v,v))
 = a_t(p_1,v)\, a_p(p_2,\sgn(pt) (v/p_1) q^m p_2 )
\\ & & \hspace{20ex} \times \, \Phi(m_1+m-m_2-\chi(pp_1/p_2 t),t_2,t_1) \ot \Phi(\chi(pp_1/p_2t),\sgn(pt) (v/p_1) q^m
p_2,v)\ .
\end{eqnarray*}
Applying $\vfi \ot \io$ to this equation and using Lemma \ref{mult.lem1}, we see that
\begin{eqnarray*}
& & [\,(\om_{\la(\Phi(m_1,p_1,t_1)),\la(\Phi(m_2,p_2,t_2))}\pi \ot \io)(\de(\Phi(m,p,t)))\,]\,\Phi(0,v,v) =
\sde_{m_1+m-m_2,\chi(pp_1/p_2 t)}\,\sde_{t_1,t_2}
\\ & & \spat\spat \times  \, (t_1)^{-2}\, a_t(p_1,v)\, a_p(p_2,\sgn(pt) (v/p_1) q^m p_2 )
\,\Phi(m+m_1-m_2,\sgn(pt) (v/p_1) q^m p_2,v) \ .
\end{eqnarray*}
Or, again by Lemma \ref{mult.lem1},
\begin{eqnarray*}
& & (\om_{f_{m_1,p_1,t_1},f_{m_2,p_2,t_2}}\pi \ot \io)(\de(\Phi(m,p,t)))\,\Phi(0,v,v) = \sde_{m_1+m-m_2,\chi(pp_1/p_2
t)}\,\sde_{t_1,t_2}
\\ & & \spat\spat \times  \, |t_2/t_1|\, a_t(p_1,v)\, a_p(p_2,\sgn(pt) (v/p_1) q^m p_2 )
\, \Phi(m+m_1-m_2,\sgn(pt) (v/p_1) q^m p_2,v) \ ,
\end{eqnarray*}
so we conclude that if $\sgn(pt) (v/p_1) q^m p_2 \in I_q$, the element $(\om_{f_{m_1,p_1,t_1},f_{m_2,p_2,t_2}}\pi \ot
\io)(\de(\Phi(m,p,t)))$ $\Phi(0,v,v)$ belongs to $\cN_{\vfi}$ and
\begin{eqnarray}
& & \la\bigl(\,(\om_{f_{m_1,p_1,t_1},f_{m_2,p_2,t_2}}\pi \ot \io)(\de(\Phi(m,p,t)))\,\Phi(0,v,v)\,\bigr) =
\sde_{m_1+m-m_2,\chi(pp_1/p_2 t)}\,\sde_{t_1,t_2} \nonumber
\\ & & \spat\spat \times  \, |v|^{-1} \, a_t(p_1,v)\, a_p(p_2,\sgn(pt) (v/p_1) q^m p_2 )
f_{m+m_1-m_2,\sgn(pt) (v/p_1) q^m p_2,v} \ . \label{mult.eq2}
\end{eqnarray}
Also recall that $(\om_{f_{m_1,p_1,t_1},f_{m_2,p_2,t_2}}\pi \ot \io)(\de(\Phi(m,p,t)))\,\Phi(0,v,v) = 0$ if $\sgn(pt)
(v/p_1) q^m p_2 \not\in I_q$.

\smallskip

Notice that Corollary \ref{mult.cor1} implies that $J \pi(\Phi(0,v,v)) J f_{n,x,y} = J \pi(\Phi(0,v,v)) f_{-n,y,x} =
\sde_{v,y}\,J f_{-n,y,x} = \sde_{v,y} f_{n,x,y}$. So we get by Proposition \ref{mult.prop1} that
\begin{eqnarray*}
& & J \pi(\Phi(0,v,v)) J \, (\om_{f_{m_1,p_1,t_1},f_{m_2,p_2,t_2}} \ot \io)(W^*)\,f_{m,p,t}  =
\\ & & \spat = \sde_{t_1,t_2}\, \sde_{m_1+m-\chi(p_1 p/t p_2),m_2}\,  \sum_{\scriptsize \begin{array}{c} y \in I_q \\ \sgn(p t)(yp_2/p_1)q^{m} \in  I_q
\end{array}} \ |t/y|\,a_{t}(p_1,y)\,a_{p}(p_2,\sgn(p t) (yp_2/p_1) q^{m})
\\ & & \hspace{50ex} \times \ \  J \pi(\Phi(0,v,v)) J \, f_{\chi(p_1 p/t p_2), \sgn(p t)
(y p_2/p_1) q^{m} , y}
\\ & & \spat = \sde_{t_1,t_2}\, \sde_{m_1+m-\chi(p_1 p/t p_2),m_2}\,  \sum_{\scriptsize \begin{array}{c} y \in I_q \\ \sgn(p t)(yp_2/p_1)q^{m} \in  I_q
\end{array}} \ |t/y|\,a_{t}(p_1,y)\,a_{p}(p_2,\sgn(p t) (yp_2/p_1) q^{m})
\\ & & \hspace{50ex} \times \ \  \sde_{v,y} \, f_{m+m_1-m_2, \sgn(p t)
(y p_2/p_1) q^{m} , y} \ .
\end{eqnarray*}
So we see that $J \pi(\Phi(0,v,v)) J \, (\om_{f_{m_1,p_1,t_1},f_{m_2,p_2,t_2}} \ot \io)(W^*)\,f_{m,p,t} = 0$ if $\sgn(p
t) (vp_2/p_1) q^{m} \not\in I_q$. If, on the other hand, $\sgn(p t) (v p_2/p_1) q^{m} \in I_q$,
\begin{eqnarray*}
& & J \pi(\Phi(0,v,v)) J \, (\om_{f_{m_1,p_1,t_1},f_{m_2,p_2,t_2}} \ot \io)(W^*)\,f_{m,p,t}
\\ & & \hspace{-2ex} \spat = \sde_{t_1,t_2}\, \sde_{m_1+m-\chi(p_1 p/t p_2),m_2}\, |t/v|\,a_{t}(p_1,v)\,a_{p}(p_2,\sgn(p t) (vp_2/p_1) q^{m})\, f_{m+m_1-m_2,
\sgn(p t) (v p_2/p_1) q^{m} , v} \ .
\end{eqnarray*}
From Eq.~(\ref{mult.eq2}) and Lemma \ref{mult.lem1}, we conclude that
\begin{eqnarray*}
& & \la\bigl(\,(\om_{f_{m_1,p_1,t_1},f_{m_2,p_2,t_2}}\pi \ot \io)(\de(\Phi(m,p,t)))\,\Phi(0,v,v)\,\bigr)
\\ & & \spat = J \pi(\Phi(0,v,v)) J \, (\om_{f_{m_1,p_1,t_1},f_{m_2,p_2,t_2}} \ot \io)(W^*)\,\la(\Phi(m,p,t)) \ .
\end{eqnarray*}
We use the net of projections $(P_L)_{L \in F(I_q)}$ introduced in the proof of Corollary \ref{mult.lem2}. If $L \in
F(I_q)$, we know by the previous equation  that  $(\om_{f_{m_1,p_1,t_1},f_{m_2,p_2,t_2}}\pi \ot
\io)(\de(\Phi(m,p,t)))\,P_L \in \cN_{\vfi}$ and
$$
\la\bigl(\,(\om_{f_{m_1,p_1,t_1},f_{m_2,p_2,t_2}}\pi \ot \io)(\de(\Phi(m,p,t)))\,P_L\,\bigr) = J \pi(P_L) J \,
(\om_{f_{m_1,p_1,t_1},f_{m_2,p_2,t_2}} \ot \io)(W^*)\,\la(\Phi(m,p,t)) \ .
$$
Therefore the $\si$-strong$^*$ closedness of $\la$ implies that $(\om_{f_{m_1,p_1,t_1},f_{m_2,p_2,t_2}}\pi \ot
\io)(\de(\Phi(m,p,t))) \in \cN_{\vfi}$ and
$$
\la\bigl(\,(\om_{f_{m_1,p_1,t_1},f_{m_2,p_2,t_2}}\pi \ot \io)(\de(\Phi(m,p,t)))\,\bigr)
 = (\om_{f_{m_1,p_1,t_1},f_{m_2,p_2,t_2}} \ot \io)(W^*)\,\la(\Phi(m,p,t)) \ .
$$
Since the linear span of such linear functionals $\om_{f_{m_1,p_1,t_1},f_{m_2,p_2,t_2}}\pi$ is norm dense in $(M_q)_*$,
the proposition now follows easily from Corollary \ref{mult.lem2} and the $\si$-strong$^*$ closedness $\la$.
\end{proof}

\medskip

Now, the proof of the left invariance of $\vfi$ can be dealt with  as in the proof of \cite[Prop.~8.15]{VaKust}.
\smallskip
\begin{proposition}
The weight $\vfi$ is left invariant with respect to $\de$.
\end{proposition}
\begin{proof}
Take $a \in \cN_{\vfi}$ and $\om \in (M_q)_*^+$. Since $\pi(M_q)$ is in  standard form with respect to $K$, there
exists a vector $v \in K$ such that $\om = \om_{v,v} \pi$. Take also an orthonormal basis $(e_i)_{i \in I}$ for $K$. It
follows from  the proof of \cite[Lem.~A.5]{VaKust} that $\bigl(\,\sum_{i \in L} \, (\om_{v,e_i}\pi \ot \io)(\de(a))^*
(\om_{v,e_i}\pi \ot \io)(\de(a))\,\bigr)_{L \in F(I)}$ is an increasing net in $M_q^+$ that converges strongly to $(\om
\ot \io)\de(a^* a)$ \ (as before, $F(I)$ denotes the set of finite subsets of $I$, directed by inclusion). Therefore
the $\si$-weak lower semi-continuity of $\vfi$ implies that
$$\vfi\bigl(\,(\om \ot \io)\de(a^* a)\,\bigr) = \sum_{i \in I} \, \vfi\bigl(\,(\om_{v,e_i}\pi \ot \io)(\de(a))^* (\om_{v,e_i}\pi \ot
\io)(\de(a))\,\bigr) \ .$$ By the previous proposition, this implies that
\begin{eqnarray*}
& & \hspace{-4ex} \vfi\bigl(\,(\om \ot \io)\de(a^* a)\,\bigr)  =  \sum_{i \in I} \, \langle \la\bigl((\om_{v,e_i}\pi
\ot \io)(\de(a))\bigr) , \la\bigl((\om_{v,e_i}\pi \ot \io)(\de(a))\bigr) \rangle
\\ & & \spat \hspace{-5ex}  =  \sum_{i \in I} \, \langle (\om_{v,e_i} \ot \io)(W^*) \la(a),(\om_{v,e_i} \ot \io)(W^*) \la(a) \rangle
=  \sum_{i \in I} \, \langle (\om_{v,e_i} \ot \io)(W^*)^* (\om_{v,e_i} \ot \io)(W^*) \la(a), \la(a) \rangle \ .
\end{eqnarray*}
So if we use first \cite[Lem.~A.5]{VaKust} and afterwards the fact that $W^*$ is an isometry (Proposition
\ref{mult.prop1}), we see that
$$ \vfi\bigl(\,(\om \ot \io)\de(a^* a)\,\bigr) =  \langle (\om_{v,v} \ot \io)(W W^*) \la(a), \la(a) \rangle
= \om_{v,v}(1)\, \langle \la(a) , \la(a) \rangle = \om(1) \, \vfi(a^* a)\ $$ and we have proven the proposition.
\end{proof}

\medskip

We have introduced the anti-unitary operator $\breve{J}$ on $H$ before Proposition \ref{coass.prop1}. Thus,
$\breve{J}(\ze^n \ot \sde_x) = \sgn(x)^{\chi(x)} \, \ze^{-n} \ot \sde_x$ for all $n \in \Z$ and $x \in I_q$. This
implies easily that $\breve{J} \Phi(m,p,t) \breve{J} = \sgn(p)^{\chi(p)}\, \sgn(t)^{\chi(t)}$ $\Phi(-m,p,t)$ for $m \in
\Z$ and $p,t \in I_q$. It shows first of all $\breve{J} M_q \breve{J} = M_q$. This allows us to define an
anti-$^*$-isomorphism $\breve{R}$ on $M_q$ such that $\breve{R}(a) = \breve{J} a^* \breve{J}$ for all $a \in M_q$.
Proposition \ref{coass.prop1} guarantees that $\flip(\breve{R} \ot \breve{R})\de = \de \breve{R}$, where $\flip$
denotes the flip map $\flip : M_q \ot M_q \rightarrow M_q \ot M_q$. Later on,  $\breve{R}$ turns out to be intimately
connected with the unitary antipode.

\smallskip

\begin{proposition}
We have that $\vfi = \vfi \breve{R}$. Thus $\vfi$ is also right invariant.
\end{proposition}
\begin{proof}
If $m \in \Z$ and $p,t \in I_q$, then $\breve{R}(\Phi(m,p,t))^* = \breve{J} \Phi(m,p,t) \breve{J} = \sgn(p)^{\chi(p)}\,
\sgn(t)^{\chi(t)}\, \Phi(-m,p,t)$. Hence, Lemma \ref{mult.lem1} implies that $\breve{R}(\Phi(m,p,t))^* \in \cN_{\vfi}$.
This same lemma implies moreover for $m,m' \in \Z$ and $p,t,p',t' \in I_q$:
\begin{eqnarray*}
& & \langle \la(\breve{R}(\Phi(m,p,t))^*) , \la(\breve{R}(\Phi(m',p',t'))^*) \rangle
\\ & & \spat = \sgn(p)^{\chi(p)}\,\sgn(t)^{\chi(t)}\, \sgn(t')^{\chi(t')}\, \sgn(p')^{\chi(p')} \,
\langle \la(\Phi(-m,p,t)) , \la(\Phi(-m',p',t')) \rangle
\\ & & \spat = \sgn(p)^{\chi(p)}\,\sgn(t)^{\chi(t)}\, \sgn(t')^{\chi(t')}\, \sgn(p')^{\chi(p')} \, |t t'|^{-1} \sde_{m,m'}\,\sde_{p,p'}\,\sde_{t,t'}
\\ & & \spat = |t t'|^{-1}\,\sde_{m,m'}\,\sde_{p,p'}\,\sde_{t,t'} = \langle \la(\Phi(m',p',t')) , \la(\Phi(m,p,t))
\rangle\ .
\end{eqnarray*}
If $a \in \langle \, \Phi(m,p,t) \mid m \in \Z, p,t \in I_q \, \rangle$, the above result implies that $\breve{R}(a)^*
\in \Nfi$ and $\|\la(\breve{R}(a)^*)\| = \|\la(a)\|$.  By Corollary \ref{mult.lem2} and the $\si$-strong$^*$ closedness
of $\la$ this implies that for all $a \in \cN_{\vfi}$, the element $\breve{R}(a)^* \in \cN_{\vfi}$ and
$\|\la(\breve{R}(a)^*)\| = \|\la(a)\|$, thus $\vfi(\breve{R}(a^* a)) = \vfi(\breve{R}(a) \breve{R}(a)^*) =
\|\la(\breve{R}(a)^*)\| = \|\la(a)\| = \vfi(a^* a)$. Since $\breve{R}^2 = \io$, it follows that $\vfi = \vfi
\breve{R}$. Because $\flip(\breve{R} \ot \breve{R})\de = \de \breve{R}$ and $\vfi$ is left invariant, we get that
$\vfi$ is right invariant.
\end{proof}

\medskip

Combined with Theorem \ref{coass.thm1},  we have proven the main result of this paper (see Definition 1 of the
Introduction for used terminology):

\begin{theorem}
The pair $(M_q,\de)$ is a von Neumann algebraic quantum group which we denote by $\widetilde{SU}_q(1,1)$.
\end{theorem}

\medskip

By Proposition \ref{mult.prop2} and the remark after  \cite[Thm.1.2]{VaKust2}, we also conclude that

\begin{proposition}
The element $W$ is the multiplicative unitary of $(M_q,\de)$ with respect to the GNS-construction  $(K,\pi,\la)$.
\end{proposition}

\medskip

The unitarity can also be proved directly from Propositions \ref{mult.prop1}, \ref{comult.ortho1}, \ref{comult.ortho2}
and \ref{comult.prop2}. We denote the antipode, unitary antipode and scaling group of $(M_q,\de)$ by $S$, $R$ and
$\tau$ respectively. All these objects are defined after the proof of \cite[Prop.~1.4]{VaKust2}. The modular element
and scaling constant of a quantum group are introduced at the same place. Because $(M_q,\de)$ is unimodular, that is,
possesses a weight that is  left and right invariant, we get the following result.

\smallskip

\begin{proposition} \label{mult.prop4}
We have that $\vfi R = \vfi$. Furthermore, the modular element of $(M_q,\de)$ is the identity operator  and the scaling
constant of $(M_q,\de)$ equals 1.
\end{proposition}
\begin{proof}
We denote the modular element and scaling constant of $(M_q,\de)$ by $\sde$ and $\nu$ respectively, see the discussion
after the proof of  \cite[Prop.~1.4]{VaKust2}. We know that $\si^\vfi_s(\sde) = \nu^s\,\sde$ for all $s \in \R$ (*) and
$\vfi R = \vfi_\sde$. Because $\vfi$ is right invariant, the uniqueness of right Haar weights implies the existence of
$\lambda > 0$ such that $\vfi R = \lambda \, \vfi$.  Thus, the uniqueness of the Radon-Nykodim derivative guarantees
that $\sde = \lambda\,1$. But, since $\de(\sde) = \sde \ot \sde$, it follows that $\sde = 1$. Therefore, (*) implies
that $\nu = 1$.
\end{proof}

\medskip

In the next part of this section we calculate explicit expressions for $S$, $R$ and $\tau$ acting on elements
$\Phi(m,p,t)$. For this purpose we calculate some explicit slices $(\io \ot \om)(W^*)$, where $\om \in B(K)_*$.

\medskip

Consider $m,m_1,m_2 \in \Z$, $p,p_1,p_2,t,t_1,t_2 \in I_q$. It follows easily from Proposition \ref{mult.prop1} that
\begin{eqnarray}
& & (\io \ot \om_{f_{m_1,p_1,t_1},f_{m_2,p_2,t_2}})(W^*)\,f_{m,p,t} \nonumber
\\ & & \spat  = |t_1/t_2|\,a_{t_1}(p,t_2)\,a_{p_1}(\sgn(p_1 t_1) (p p_2/t_2) q^{-m_1}, p_2)\, f_{m_1-m_2+m,\sgn(p_1 t_1) (p p_2/t_2)
q^{-m_1},t} \label{mult.eq3}
\end{eqnarray}
if $\sgn(p_1 t_1) (p p_2/t_2) q^{-m_1} \in I_q$ and $\chi(p_1 t_2/p_2 t_1) = m_2 - m_1$. If, on the other hand,
$\sgn(p_1 t_1) (p p_2/t_2) q^{-m_1}$  $\not\in I_q$ or $\chi(p_1 t_2/p_2 t_1) \not= m_2 - m_1$, then $(\io \ot
\om_{f_{m_1,p_1,t_1},f_{m_2,p_2,t_2}})(W^*)\,f_{m,p,t} = 0$.

\medskip

\begin{proposition} \label{mult.prop3}
Consider $m_1,m_2,m_3,m_4 \in \Z$ and $p,t \in I_q$. Then we define $\om \in B(K)_*$ as the absolutely convergent sum
\begin{eqnarray}
& & \om = \sum_{\scriptsize \begin{array}{c} x,y \in I_q \\ \sgn(p) q^{m_3} x \in I_q \\ \sgn(t) q^{m_4} y \in I_q
\end{array}} \, \sgn(x)^{\chi(x)} \, \sgn(y)^{\chi(y)}\, (-1)^{\chi(xy)}\,|xy| a_p(\sgn(p) q^{m_3} x, x)\,a_t(\sgn(t) q^{m_4} y,
y)\nonumber
\\ & & \hspace{21ex} \times \  \ \om_{f_{m_1+\chi(x/y),\sgn(p)q^{m_3}x,\sgn(t)q^{m_4}y},f_{m_2+\chi(x/y),x,y}} \
.\label{mult.eq5}
\end{eqnarray}
Then $(\io \ot \om)(W^*) = t^2\,q^{-m_1-m_3}\, (-1)^{m_2}\, \sde_{\chi(t/p),m_1}
\,\sde_{m_3-m_4,m_2-m_1}\,\pi(\Phi(m_1-m_2,p,t))$.
\end{proposition}
\begin{proof} Let $s \in I_q$ and $m \in \Z$. By Definition \ref{comult.def1},  we have for $z \in I_q$
satisfying $\sgn(s)\,q^m z \in I_q$,
\begin{eqnarray*}
\\ & & |a_s(\sgn(s) q^m z,z)|
\\ & & \hspace{2ex} \leq c_q\, \nu(sq^{-m}) \,\sqrt{\frac{(-\kappa(s);q^2)_\infty}{(q^2;q^2)_\infty}}\,\,
|z|\,\sqrt{(-\kappa(z);q^2)_\infty}\,\, |\Psi(-q^2/\kappa(s);(q^{2(1+m)}/s^2) \kappa(z);\sgn(s) q^{2(m+1)})| \ .
\end{eqnarray*}
which by Lemma \ref{alt.lem1} implies that $\sum_{z \in I_q, \sgn(s) q^m z \in I_q} \, |a_s(\sgn(s) q^m z,z)| <
\infty$. It follows that the sum in the statement of this proposition is absolutely convergent. So, as a norm
convergent limit of finite sums of vector functionals, we obtain $\om \in B(K)_*$. Take $n \in \Z$ and $c,d \in I_q$.
Now
\begin{eqnarray*}
& & (\io \ot \om)(W^*)\,f_{n,c,d}
 = \sum_{\scriptsize \begin{array}{c} x,y \in I_q \\ \sgn(p) q^{m_3} x \in I_q \\ \sgn(t) q^{m_4} y \in I_q
\end{array}} \, \sgn(x)^{\chi(x)} \, \sgn(y)^{\chi(y)}\, (-1)^{\chi(xy)}\,|xy| a_p(\sgn(p) q^{m_3} x, x)
\\ & & \hspace{16ex} \times \  \,a_t(\sgn(t) q^{m_4} y, y)\, (\io \ot \om_{f_{m_1+\chi(x/y),\sgn(p)q^{m_3}x,\sgn(t)q^{m_4}y},f_{m_2+\chi(x/y),x,y}})(W^*)\,f_{n,c,d} \ .
\end{eqnarray*}
By Eq.~(\ref{mult.eq3}), we get immediately that $(\io \ot \om)(W^*)\,f_{n,c,d} = 0$ if $\sgn(pt) c q^{-m_1} \not\in
I_q$ or $m_3 - m_4 \not= m_2 - m_1$. Now suppose that $\sgn(pt) c q^{-m_1} \in I_q$ and $m_3 - m_4 = m_2 - m_1$. Then
Eq.~(\ref{mult.eq3}) and Proposition \ref{comult.prop2} imply that
\begin{eqnarray*}
& & (\io \ot \om)(W^*)\,f_{n,c,d}
\\ & &
\\ & & \spat = \sum_{\scriptsize \begin{array}{c} x,y \in I_q \\ \sgn(p) q^{m_3} x \in I_q \\ \sgn(t) q^{m_4} y \in I_q
\end{array}} \, \sgn(x)^{\chi(x)} \, \sgn(y)^{\chi(y)}\, (-1)^{\chi(xy)}\,|xy| a_p(\sgn(p) q^{m_3} x, x)\,a_t(\sgn(t) q^{m_4} y, y)
\\ & & \hspace{19ex} \times \  \ q^{m_4}\,a_{\sgn(t)q^{m_4}y}(c,y)\,
a_{\sgn(p)q^{m_3}x}(\sgn(pt)cq^{-m_1},x)\,f_{m_1-m_2+n,\sgn(pt) c q^{-m_1},d}
\\ & &
\\ & & \spat = \sum_{\scriptsize \begin{array}{c} x,y \in I_q \\ \sgn(p) q^{m_3} x \in I_q \\ \sgn(t) q^{m_4} y \in I_q
\end{array}} \, \sgn(x)^{\chi(x)} \, \sgn(y)^{\chi(y)}\, (-1)^{\chi(xy)}\,|xy| a_p(\sgn(p) q^{m_3} x, x)\,a_t(\sgn(t) q^{m_4} y, y)
\\ & & \hspace{19ex} \times \  \ q^{m_4}\,(-1)^{\chi(cy)+m_4}\,\sgn(y)^{\chi(y)}\, |c/q^{m_4}y|\,a_c(\sgn(t)q^{m_4}y,y)\,(-1)^{\chi(cx)+m_3+m_1}
\\ & & \hspace{19ex} \times \ \  \sgn(x)^{\chi(x)}\,|c q^{-m_1}/q^{m_3} x|\,    a_{\sgn(pt)cq^{-m_1}}(\sgn(p)q^{m_3}x,x)\,f_{m_1-m_2+n,\sgn(pt) c q^{-m_1},d}
\end{eqnarray*}
\begin{eqnarray*}
\\ & & \spat = \sum_{\scriptsize \begin{array}{c} x,y \in I_q \\ \sgn(p) q^{m_3} x \in I_q \\ \sgn(t) q^{m_4} y \in I_q
\end{array}} \, c^2\,q^{-m_1-m_3}\, (-1)^{m_1+m_3+m_4}\, a_p(\sgn(p) q^{m_3} x, x)\,a_{\sgn(pt)cq^{-m_1}}(\sgn(p)q^{m_3}x,x)
\\ & & \hspace{20ex} \times \  \ a_t(\sgn(t) q^{m_4} y, y)\,a_c(\sgn(t)q^{m_4}y,y)\,f_{m_1-m_2+n,\sgn(pt) c
q^{-m_1},d}\ .
\end{eqnarray*}
\smallskip
This sum is absolutely convergent so we can compute it in any order we deem useful. Therefore the orthogonality
relations in Proposition \ref{comult.ortho1}  imply that
\begin{eqnarray}
& & (\io \ot \om)(W^*)\,f_{n,c,d} =  c^2\,q^{-m_1-m_3}\, (-1)^{m_1+m_3+m_4}\, \sde_{p,\sgn(pt)cq^{-m_1}} \,
\sde_{t,c}\,
 \,f_{n+m_1-m_2,\sgn(pt) c q^{-m_1},d} \nonumber
\\ & & \spat = t^2\,q^{-m_1-m_3}\, (-1)^{m_1+m_3+m_4}\, \sde_{|p|,|t|q^{-m_1}} \, \sde_{t,c} \, \,f_{n+m_1-m_2,p,d}
\nonumber
\\ & & \spat = t^2\,q^{-m_1-m_3}\, (-1)^{m_1+m_3+m_4}\,
\sde_{|p|,|t|q^{-m_1}} \, \pi(\Phi(m_1-m_2,p,t))\,f_{n,c,d} \ ,\label{mult.eq4}
\end{eqnarray}
where $\pi$ is the GNS-representation of $\vfi$. If $\sgn(pt) c q^{-m_1} \not\in I_q$, then
$$\sde_{|p|,|t|q^{-m_1}} \, \pi(\Phi(m_1-m_2,p,t))\,f_{n,c,d} = \sde_{|p|,|t|q^{-m_1}}\,\sde_{t,c}\,f_{m_1-m_2+n,p,d} = 0\ ,$$
implying that Eq.~(\ref{mult.eq4}) also holds if $\sgn(pt) c q^{-m_1} \not\in I_q$ and $m_3-m_4 = m_2 - m_1$. Therefore
the lemma follows.
\end{proof}

\medskip

Now we can easily calculate the action of $S$ on elements of the form $\Phi(m,p,t)$.

\smallskip

\begin{proposition} \label{mult.prop6}
Consider $m \in \Z$ and $p,t \in I_q$. Then $\Phi(m,p,t) \in D(S)$ and $S(\Phi(m,p,t)) =
\sgn(p)^{\chi(p)}\,\sgn(t)^{\chi(t)}\,(-1)^m\,q^m\, \Phi(m,t,p)$.
\end{proposition}

\begin{proof} Define $\om \in B(K)_*$ as in Eq.~(\ref{mult.eq5}) for $m_1 = \chi(t/p)$, $m_2 = \chi(t/p) - m$, $m_3 = -
\chi(t/p)$ and $m_4 = - \chi(t/p) +m$. Then Proposition \ref{mult.prop3} implies that $\pi(\Phi(m,p,t)) =
t^{-2}\,(-1)^{m+\chi(pt)}\, (\io \ot \om)(W^*)$.  Thus \cite[Prop.~8.3]{VaKust} implies that $\Phi(m,p,t) \in
D(S^{-1})$ and
\begin{equation}  \label{mult.eq7}
\pi\bigl(S^{-1}(\Phi(m,p,t))\bigr) = t^{-2}\,(-1)^{m+\chi(pt)}\, (\io \ot \om)(W) = t^{-2}\,(-1)^{m+\chi(pt)}\, (\io
\ot \bar{\om})(W^*)^* \ ,
\end{equation}
where, as is customary, $\bar{\om} \in B(K)_*$ is given by $\bar{\om}(a) = \overline{\om(a^*)}$ for all $a \in B(K)$.
By definition,
\begin{eqnarray*}
& & \om = \sum_{\scriptsize \begin{array}{c} x,y \in I_q \\ (p/|t|) x \in I_q \\ (|p|/t) q^m y \in I_q
\end{array}} \, \sgn(x)^{\chi(x)} \, \sgn(y)^{\chi(y)}\, (-1)^{\chi(xy)}\,|xy| a_p((p/|t|)x, x)\,a_t((|p|/t)q^m  y,
y)
\\ & & \hspace{21ex} \times \  \ \om_{f_{\chi(t/p)+\chi(x/y),(p/|t|)x,(|p|/t)q^my},f_{-m+\chi(t/p)+\chi(x/y),x,y}} \ .
\end{eqnarray*}
Thus,
\begin{eqnarray*}
\bar{\om} & = & \sum_{\scriptsize \begin{array}{c} x,y \in I_q \\ (p/|t|) x \in I_q \\ (|p|/t) q^m y \in I_q
\end{array}} \, \sgn(x)^{\chi(x)} \, \sgn(y)^{\chi(y)}\, (-1)^{\chi(xy)}\,|xy|\, a_p((p/|t|)x, x)\,a_t((|p|/t)q^m  y,
y)\rm
\\ & & \hspace{21ex} \times \  \ \om_{f_{-m+\chi(t/p)+\chi(x/y),x,y},f_{\chi(t/p)+\chi(x/y),(p/|t|)x,(|p|/t)q^my}} \
\\ & &
\\ & = & \sum_{\scriptsize \begin{array}{c} x,y \in I_q \\ (p/|t|) x \in I_q \\ (|p|/t) q^m y \in I_q
\end{array}} \, (-1)^{\chi(xy)}\,|xy|\, \sgn(p)^{\chi(p)} \, \sgn((p/|t|)x)^{\chi((p/|t|)x)}\, a_p(x,(p/|t|)x)
\\ & & \hspace{23ex} \times \  \ \sgn(t)^{\chi(t)}\,\sgn((|p|/t)q^m y)^{\chi((|p|/t)q^m y)} \, a_t(y,(|p|/t)q^m  y)
\\ & & \hspace{23ex} \times \  \  \om_{f_{-m+\chi(t/p)+\chi(x/y),x,y},f_{\chi(t/p)+\chi(x/y),(p/|t|)x,(|p|/t)q^my}} \ ,
\end{eqnarray*}
where we used Proposition \ref{comult.prop2} in the last equation. A simple change in summation variables $x' =
(p/|t|)x$ and $y' = (|p|/t) y q^{m}$ then reveals that
\begin{eqnarray*}
& & \bar{\om} = \sum_{\scriptsize \begin{array}{c} x',y' \in I_q \\ (|t|/p) x' \in I_q \\ (t/|p|) q^{-m} y' \in I_q
\end{array}} \, \sgn(p)^{\chi(p)}\,\sgn(t)^{\chi(t)}\,(-1)^m\,q^{-m} (t^2/p^2)\,\sgn(x')^{\chi(x')} \, \sgn(y')^{\chi(y')}\, (-1)^{\chi(x'y')}\,|x'y'|
\\ & & \hspace{23ex} \times \  \  a_p((|t|/p)x',x')\,    a_t((t/|p|)q^{-m}y',y')
\\ & & \hspace{23ex} \times \  \
\om_{f_{\chi(t/p)+\chi(x'/y'),(|t|/p)x',(t/|p|)q^{-m}y'},f_{m+\chi(t/p)+\chi(x'/y'),x',y'}}\ ,
\end{eqnarray*}
which upon close inspection shows that $\bar{\om}$ is $\sgn(p)^{\chi(p)}\,\sgn(t)^{\chi(t)}\,(-1)^m\,q^{-m} (t^2/p^2)$
times the functional in Eq.~(\ref{mult.eq5}) where $m_1 = \chi(t/p)$, $m_2 = \chi(t/p)+m$, $m_3 = \chi(t/p)$ and $m_4 =
\chi(t/p) - m$. By Proposition \ref{mult.prop3} this implies that
\begin{eqnarray*}
& & (\io \ot \bar{\om})(W^*) = \sgn(p)^{\chi(p)}\,\sgn(t)^{\chi(t)}\,(-1)^m\,q^{-m} (t^2/p^2)\, t^2\,(p^2/t^2)\,
(-1)^{m+\chi(pt)} \, \pi(\Phi(-m,p,t))
\\ & & \spat = \sgn(p)^{\chi(p)}\,\sgn(t)^{\chi(t)}\,(-1)^m\,q^{-m}\,t^2\,(-1)^{m+\chi(pt)} \, \pi(\Phi(-m,p,t)) \ .
\end{eqnarray*}
Therefore Eq.~(\ref{mult.eq7}) implies that
$$S^{-1}(\Phi(m,p,t)) =  \sgn(p)^{\chi(p)}\,\sgn(t)^{\chi(t)}\,(-1)^m\,q^{-m}\,
\Phi(m,t,p) \ ,$$
proving the result. 
\end{proof}

\smallskip

It is now possible to recognize the polar decomposition of the antipode. For this purpose, define the anti-unitary
operator $I$ on $H$ and the strictly positive operator $Q$ in $H$ such that $\langle \, \ze^n \ot \sde_x \mid n \in \Z,
x \in I_q \,\rangle$ is a core for $Q$ and $I(\ze^n \ot \sde_x) = (-1)^n\,\sgn(x)^{\chi(x)}\, \ze^{-n} \ot \sde_x$ and
$Q(\ze^n \ot \sde_x) = q^{2n}\,\ze^n \ot \sde_x$ for all $n \in \Z$ and $x \in I_q$. Recall that the unitary antipode
$R$ and scaling group $\tau$ are introduced after the proof of \cite[Prop.~1.4]{VaKust2}.

\smallskip

\begin{proposition}
The unitary antipode $R$ and the scaling group $\tau$ for $(M_q,\de)$ are such that $R(a) = I a^* I$ and $\tau_s(a) =
Q^{is} a Q^{-is}$ for all $a \in M_q$ and $s \in \R$. Let $m \in \Z$ and $p,t \in I_q$. Then $R(\Phi(m,p,t)) =
\sgn(p)^{\chi(p)}\,\sgn(t)^{\chi(t)}\,(-1)^m\, \Phi(m,t,p)$. If $z \in \C$, then the element $\Phi(m,p,t)$ belongs to
$D(\tau_z)$ and \newline $\tau_z(\Phi(m,p,t)) = q^{2miz}\, \Phi(m,p,t)$.
\end{proposition}
\begin{proof}
By Proposition \ref{mult.prop4} and the discussion after the proof of \cite[Prop.~1.4]{VaKust2} we know that there
exists a strictly positive operator $P$ on $K$ such that $P^{is} \la(a) = \la(\tau_s(a))$ for all $a \in \Nfi$ and $s
\in \R$. Note that $\pi(\tau_s(a)) = P^{is} \pi(a) P^{-is}$ for all $s \in \R$.

Choose $m \in \Z$ and $p,t \in \R$. Since $S = R \tau_{-\frac{i}{2}}$ and $R$ and $\tau$ commute, we see that
$\tau_{-i} = S^2$. Thus, Proposition \ref{mult.prop6} implies that $\Phi(m,p,t) \in D(\tau_{-i})$ and
$\tau_{-i}(\Phi(m,p,t)) = q^{2m}\,\Phi(m,p,t)$. Thus, Lemma \ref{mult.lem1} and the von Neumann algebraic version of
\cite[Prop.~4.4]{Kust1} guarantee that $\la(\Phi(m,p,t)) \in D(P)$ and $P\,\la(\Phi(m,p,t)) =
q^{2m}\,\la(\Phi(m,p,t))$. Therefore, Lemma \ref{mult.lem1} tells us that $f_{m,p,t} \in D(P)$ and $P f_{m,p,t} =
q^{2m}\,f_{m,p,t} = (Q \ot 1)f_{m,p,t}$. Because such elements $f_{m,p,t}$ form a core for $Q \ot 1$, we conclude that
$Q \ot 1 \subseteq P$. Taking the adjoint of this inclusion, we see that $P \subseteq Q \ot 1$. As a consequence, $P =
Q \ot 1$. So if $a \in M_q$ and $s \in \R$, we see that $\tau_s(a) \ot 1 = \pi(\tau_s(a)) = P^{is} \pi(a) P^{-is} =
Q^{is} a Q^{-is} \ot 1$, implying that $\tau_s(a) = Q^{is} a Q^{-is}$.

\smallskip

Note that $Q^{is}(\ze^n \ot \sde_x) = q^{2nis}\,\ze^n \ot \sde_x$ for all $n \in \Z$, $x \in I_q$ and $s \in \R$. Now a
straightforward calculation reveals that $\tau_s(\Phi(m,p,t)) =   Q^{is} \Phi(m,p,t) Q^{-is} = q^{2mis}\, \Phi(m,p,t)$
for all $m \in \Z$,  $p,t \in I_q$ and $s \in \R$.

Fix $m \in \Z$ and $p,t \in I_q$ for the moment. By definition $S = R \tau_{-\frac{i}{2}}$. The previous paragraph
implies that $\tau_{-\frac{i}{2}}(\Phi(m,p,t)) = q^m\,\Phi(m,p,t)$. Therefore Proposition \ref{mult.prop6} guarantees
that $R(\Phi(m,p,t)) = q^{-m}\, S(\Phi(m,p,t)) = \sgn(p)^{\chi(p)}\,\sgn(t)^{\chi(t)}\,(-1)^m\,\Phi(m,t,p)$. Using the
definition of $I$, one checks easily that also $I \Phi(m,p,t)^* I = I \Phi(-m,t,p) I =
\sgn(p)^{\chi(p)}\,\sgn(t)^{\chi(t)}\,(-1)^m\, \Phi(m,t,p)$ and thus by the previous calculation, $R(\Phi(m,p,t)) = I
\Phi(m,p,t)^* I$. Since the von Neumann algebra $M_q$ is generated by such elements $\Phi(m,p,t)$, we conclude that
$R(a) = I a^* I$ for all $a \in M_q$.
\end{proof}

\medskip

One easily checks that $R$ and $\breve{R}$ are connected through the formula $\breve{R}(\Phi(m,p,t)) = (-1)^m\,
R(\Phi(m,p,t))$ if $m \in \Z$ and $p,t \in I_q$.

\smallskip

Consider $z \in \C$. Since the $^*$-algebra $\langle\, \Phi(m,p,t) \mid m \in \Z, p,t \in I_q \, \rangle$ is
$\si$-strongly$^*$ dense in $M_q$ and this same $^*$-algebra is clearly invariant under each $\tau_s$, where $s \in
\R$, it follows from the von Neumann algebraic version of  \cite[Cor.~1.22]{Kust2} that $\langle\, \Phi(m,p,t) \mid m
\in \Z, p,t \in I_q \, \rangle$ is a $\si$-strong$^*$ core for $\tau_z$.

If $n \in \Z$, the fact that  $S^n = R \, \tau_{-n\frac{i}{2}}$ implies that $\langle\, \Phi(m,p,t) \mid m \in \Z, p,t
\in I_q \, \rangle$ is a $\si$-strong$^*$ core for $S^n$.

\medskip

Following \cite[Prop.~1.6]{VaKust2} we associate to the von Neumann algebraic quantum group $(M_q,\de)$ a reduced
\cst-algebraic quantum group in the sense of \cite{VaKust}. Therefore we define the \cst-subalgebra $A_q$ of $M_q$ such
that $\pi(A_q)$ is the norm closure, in $B(K)$, of the set $\{\,(\io \ot \om)(W^*) \mid \om \in B(H)_*\,\}$. Then
\cite[Prop.~1.6]{VaKust2} guarantees that $(A_q,\de\!\!\restriction_{A_q})$ is a reduced \cst-algebraic quantum group
in the sense of \cite{VaKust}. In the next proposition we give an explicit description for $A_q$.

\smallskip

We will use the following notation. For $f \in C(\T \times I_q)$ and $x \in I_q$ we define $f_x \in C(\T)$ such that
$f_x(\lambda) = f(\lambda,x)$ for all $\lambda \in \T$.

\smallskip

\begin{proposition} \label{mult.prop5}
Denote by $\cC$ the \cst-algebra of all functions $f \in C(\T \times I_q)$ such that $\text{$(1)$}$\ $f_x$ converges
uniformly to $0$ as $x \rightarrow 0$ and \text{$(2)$}\ $f_x$ converges uniformly to a constant function as $x
\rightarrow \infty$. Then $A_q$ is the norm closed linear span, in $B(H)$, of the set $\{\,\rho_p \,M_f \mid f \in \cC,
p \in \qz\,\}$.
\end{proposition}
\begin{proof}
Let $p,t \in I_q$ and define $F_{p,t} \in \cF(I_q)$ such that $F_{p,t}(x) = a_p(x,t)$ for all $x \in I_q$. The reason
for introducing the function $F_{p,t}$ stems from the fact that
\begin{eqnarray}
& & (\io \ot \om_{f_{m_1,p_1,t_1},f_{m_2,p_2,t_2}})(W^*) \nonumber
\\ & & \spat = \sde_{\chi(p_1 t_2/p_2 t_1),m_2 - m_1} \,\,
|t_1/t_2| \,\, \pi\left(\ \rho_{\sgn(p_1 t_1) (t_2/p_2) q^{m_1}}\,\, M_{\ze^{m_1-m_2} \ot T_{\sgn(p_1 t_1) (p_2/t_2)
q^{-m_1}}(F_{p_1,p_2}) F_{t_1,t_2}}\ \right) \nonumber
\\ & & \label{cst.eq4}
\end{eqnarray}
for $p_1,p_2,t_1,t_2 \in I_q$ and $m_1,m_2 \in \Z$\, (and which follows from Eq.~(\ref{mult.eq3})\,). Since $B(K)_*$ is
the closed linear span of vector functionals, $A_q$ is the closed linear span of elements of the form (\ref{cst.eq4}).

\smallskip

Let $p,t \in I_q$. Take  $m \in \Z$ such that $|p/t| = q^m$. We consider two different cases:
\begin{trivlist}
\item[\ \,\boldmath\bf $\sgn(p) = \sgn(t)$ :] Note first of all that in this case, $F_{p,t}(x) = 0$ for all $x \in I_q^-$. Definition \ref{comult.def1} implies the existence of a bounded function $h :
I_q^+ \rightarrow \R$ such that
$$F_{p,t}(x) = h(x) \, \frac{\nu(pt/x)}{\sqrt{(-\kappa(x);q^2)_\infty}} \,\,\Psi(-q^2/\kappa(t);q^2\kappa(x/t);q^2
\kappa(x/p))$$ for all $x \in I_q^+$. If $x \rightarrow 0$, we have that $\nu(pt/x) \rightarrow 0$,
$(-\kappa(x);q^2)_\infty \rightarrow 1$ and $\Psi(-q^2/\kappa(t);q^2\kappa(x/t);q^2 \kappa(x/p))$ $\rightarrow
\Psi(-q^2/\kappa(t);0;0)$. It follows that $F_{p,t}(x) \rightarrow 0$ as $x \rightarrow 0$.

\smallskip

Proposition \ref{comult.prop2} and Definition \ref{comult.def1} imply for $x \in I_q^+$,
\begin{eqnarray}
& & F_{p,t}(x) = a_p(x,t) = \sgn(p)^{\chi(pt)} \, a_p(t,x) = (-1)^{\chi(pt)} \,
\sgn(p)^{\chi(pt)}\,c_q\,(-\kappa(p);q^2)_\infty^{\frac{1}{2}} \nonumber
\\ & & \spat\spat (-\kappa(t);q^2)_\infty^{-\frac{1}{2}} \,\,x\,\nu(px/t)\,(-\kappa(x);q^2)_\infty^{\frac{1}{2}}\,
\Psi(-q^2/\kappa(x);q^2 \kappa(t/x);q^2 \kappa(t/p)) \label{cst.eq1}
\end{eqnarray}
If $x \rightarrow \infty$, then $\Psi(-q^2/\kappa(x);q^2 \kappa(t/x);q^2 \kappa(t/p)) \rightarrow \Psi(0;0;q^2
\kappa(t/p)) = (q^2 \kappa(t/p);q^2)_\infty$, where we used \cite[Eq.~(1.3.16)]{Gas}. \!\! \inlabel{cst.eq2}

Let $x \in I_q^+$ and choose $s \in \Z$ such that $x=q^s$. Then Result \ref{app.res1} implies that
\begin{eqnarray*}
& & x\,\nu(px/t)\,(-\kappa(x);q^2)_\infty^{\frac{1}{2}} = q^s \, q^{\frac{1}{2}(s+m-1)(s+m-2)}\,
\sqrt{(-q^{2s};q^2)_\infty}
\\ & & \spat = q^s \, q^{\frac{1}{2}s(s-1)}\,q^{\frac{1}{2}(m-1)(m-2)}\,q^{(m-1)s}\,
q^{-\frac{1}{2}s(s-1)}\,(-1,-q^2;q^2)_\infty^{\frac{1}{2}}\,(-q^{2(1-s)};q^2)_\infty^{-\frac{1}{2}}
\\ & & \spat =  q^{\frac{1}{2}(m-1)(m-2)}\,
\sqrt{2}\,(-q^2;q^2)_\infty \,x^m\,(-q^2/x^2;q^2)_\infty^{-\frac{1}{2}} \ .
\end{eqnarray*}
So if $p=t$ (and thus $m=0$),  $x\,\nu(px/t)\,(-\kappa(x);q^2)_\infty^{\frac{1}{2}} \rightarrow
\sqrt{2}\,q\,(-q^2;q^2)_\infty$ as $x \rightarrow \infty$. If, on the other hand, $p > t$ (and thus $m < 0$), we see
that $x\,\nu(px/t)\,(-\kappa(x);q^2)_\infty^{\frac{1}{2}} \rightarrow 0$ as $x \rightarrow \infty$. If we combine this
with Eq.~(\ref{cst.eq1}) and the convergence in (\ref{cst.eq2}), we conclude, since $c_q^{-1} =
\sqrt{2}\,q\,(q^2,-q^2;q^2)_\infty$, that
\begin{trivlist}
\item[\ \,(1)] $F_{p,p}(x) \rightarrow 1$ as $x \rightarrow \infty$,
\item[\ \,(2)] if $p > t$, then $F_{p,t}(x) \rightarrow 0$ as $x \rightarrow \infty$ \ .
\end{trivlist}
But Proposition \ref{comult.prop2} implies that $|F_{p,t}| =  |t/p|\,|F_{t,p}|$. By the convergence in (2), this
implies that if $p < t$, also  $F_{p,t}(x) \rightarrow 0$ as $x \rightarrow \infty$.

\smallskip

\item[\ \,\boldmath\bf $\sgn(p) \not= \sgn(t)$ :] Note that in this case, $F_{p,t}(x) = 0$ for all $x \in I_q^+$.
Completely similar as in the beginning of the first part of this proof, one shows that $F_{p,t}(x) \rightarrow 0$ as $x
\rightarrow 0$.
\end{trivlist}
\smallskip
From the previous discussion we only have to remember that $F_{p,t}(x) \rightarrow 0$ as $x \rightarrow 0$. Provided $p
\not= t$, we have that $F_{p,t}(x) \rightarrow 0$ as $x \rightarrow \infty$. And $F_{p,p}(x) \rightarrow 1$ as $x
\rightarrow \infty$. \!\!\!\inlabel{cst.eq3}

\smallskip

Let $B$ denote the norm closed linear span of $\{\,\rho_p \,M_f \mid f \in \cC, p \in \qz\,\}$ in $B(H)$.

Use the notation of Eq.~(\ref{cst.eq4}). If $T_{\sgn(p_1 t_1) (p_2/t_2) q^{-m_1}}(F_{p_1,p_2}) F_{t_1,t_2} \in
C_0(I_q)$, then Eqs.~(\ref{cst.eq4}) and (\ref{cst.eq3}) immediately guarantees that $(\io \ot
\om_{f_{m_1,p_1,t_1},f_{m_2,p_2,t_2}})(W^*) \in \pi(B)$. If, on the other hand, \newline $T_{\sgn(p_1 t_1) (p_2/t_2)
q^{-m_1}}(F_{p_1,p_2}) F_{t_1,t_2}$ $\not\in C_0(I_q)$, we must by (\ref{cst.eq3}) necessarily have that $t_1=t_2$ and
$p_1 = p_2$. Thus, Eq.~(\ref{cst.eq4}) implies in this case that
\begin{eqnarray*}
 (\io \ot \om_{f_{m_1,p_1,t_1},f_{m_2,p_1,t_1}})(W^*)
 &  = & \sde_{0,m_2 - m_1} \,\,
 \pi\left(\ \rho_{q^{m_1}|t_1/p_1|}\,\, M_{\ze^{m_1-m_2} \ot T_{q^{-m_1}|p_1/t_1|}(F_{p_1,p_1}) F_{t_1,t_1}}\ \right)
\\ &  = & \sde_{0,m_2 - m_1} \,\,
 \pi\left(\ \rho_{q^{m_1}|t_1/p_1|}\,\, M_{1 \ot T_{q^{-m_1}|p_1/t_1|}(F_{p_1,p_1}) F_{t_1,t_1}}\ \right)
\in \pi(B) \ ,
\end{eqnarray*}
where we used (\ref{cst.eq3}) in the last relation. Hence $\pi(A_q) \subseteq \pi(B)$,  thus $A_q \subseteq B$.

\medskip

Proposition \ref{mult.prop3} implies that $\rho_{t^{-1}p} \,M_{\ze^m \ot \sde_p} = \Phi(m,t,p) \in A_q$ for all $p,t
\in I_q$ and $m \in \Z$. It follows that $\rho_x \,M_{\ze^m \ot \sde_p} \in A_q$ for all $m \in \Z$, $x \in \qz$ and $p
\in I_q$. Thus,  $\rho_x \,M_f \in A_q$ for all $x \in \qz$ and $f \in C_0(\T \times I_q)$.

If $x \in - q^{\Z}$ and $f \in \cC$, it is easy to see that there exists $g \in C_0(\T \times I_q)$ such that
$\rho_x\,M_f = \rho_x\,M_g$. Hence,  $\rho_x\,M_f \in A_q$ in this case.

Now fix $x \in q^\Z$ and $f \in \cC$. Take $m \in \Z$ such that $x = q^m$. Because $f \in \cC$,  there exists  $c \in
\C$ such that $f_y \rightarrow c\,1_\T$ uniformly as $y \rightarrow \infty$. Consequently, $f - c\,(1_\T \ot
T_{x^{-1}}(F_{1,1}) F_{1,1})$ belongs to $C_0(\T \times I_q)$ by (\ref{cst.eq3}). As such, $\rho_x\,M_{f - c\,(1_\T \ot
T_{x^{-1}}(F_{1,1}) F_{1,1})} \in A_q$. Moreover, Eq.~(\ref{cst.eq4}) guarantees that $\pi(\rho_x\,M_{1 \ot
T_{x^{-1}}(F_{1,1}) F_{1,1}}) = (\io \ot \om_{f_{m,1,1},f_{m,1,1}})(W^*) \in \pi(A_q)$. Hence $\rho_x\,M_{1 \ot
T_{x^{-1}}(F_{1,1}) F_{1,1}} \in A_q$. It follows that $\rho_x \, M_f \in A_q$. Hence $B \subseteq A_q$, proving that
$B = A_q$.
\end{proof}

\bigskip\medskip

\sectie{Proof of the coassociativity of the comultiplication}

\bigskip

This section is devoted to the proof of Theorem \ref{coass.thm1}. First one proves that $(\de \ot \io)\de(\gamma) =
(\io \ot \de)\de(\gamma)$. It is easy to prove that both operators agree on the intersection of their domain, but it
takes a lot of  work to establish the equality of their domains. In order to enhance the continuity of the exposition,
we give the proof in the second half of this section. Once we have proven that $(\de \ot \io)\de(\gamma) = (\io \ot
\de)\de(\gamma)$, it is straightforward to show that $(\de \ot \io)\de(x) = (\io \ot \de)\de(x)$ for all $x \in
M_q^\diamond$. In the final step of the proof of the coassociativity it is shown that $(\de \ot \io)\de(u) = (\io \ot
\de)\de(u)$.

\medskip

Since $\de(x) = V^*(1_{L^2(\T^2)} \ot x) V$ for all $x \in M_q$,  we get for any closed densely defined operator $T$ in
$H$, affiliated with $M_q$,
\begin{equation}\label{coass.eq2}
(\de \ot\io)\de(T) = (V^* \ot 1_H)(1_{L^2(\T^2)} \ot V^*) (1_{L^2(\T^2)} \ot 1_{L^2(\T^2)} \ot T) (1_{L^2(\T^2)} \ot V)
(V \ot 1_H)
\end{equation}
and
\begin{equation} \label{coass.eq3}
(\io \ot\de)\de(T) = (1_H \ot V^*) V_{13}^*(1_{L^2(\T^2)} \ot 1_{L^2(\T^2)} \ot T) V_{13} (1_H \ot V)
\end{equation}
where $V_{13} : H \ot L^2(\T^2) \ot H \rightarrow L^2(\T^2) \ot L^2(\T^2) \ot H$ is defined as usual.

\medskip

Using the map $\de_0^{(2)} : \cA_q \rightarrow \cA_q \odot \cA_q \odot \cA_q$ defined by $\de_0^{(2)} = (\de_0 \odot
\io)\de_0 = (\io \odot \de_0)\de_0$, we get adjointable operators $\de_0^{(2)}(\al_0)$, $\de_0^{(2)}(\al_0^\dag)$,
$\de_0^{(2)}(\gamma_0)$, $\de_0^{(2)}(\gamma_0^\dag)$  in $\cL^+(E \odot E \odot E)$. Proposition \ref{comult.prop1}
will imply the next result.

\begin{proposition} \label{coass.prop2}
If $T=\al$ or $T = \gamma$, then
$$\begin{array}{lclclcl}  \de_0^{(2)}(T_0) & \subseteq & (\de \ot\io)\de(T)   & \hspace{12ex} &  \de_0^{(2)}(T_0) & \subseteq & (\io \ot\de)\de(T)
\\ \de_0^{(2)}(T_0^\dag) & \subseteq & (\de \ot\io)\de(T^*)   & \hspace{12ex} &  \de_0^{(2)}(T_0^\dag) & \subseteq & (\io
\ot\de)\de(T^*) \ .
\end{array}$$
\end{proposition}
\begin{proof} 
We only proof the inclusion for $T=\al$. The other ones are dealt with in the same way. Proposition
\ref{comult.prop1} tells us that $\al_0 \odot \al_0 + q\, e_0 \gamma_0^\dag \odot \gamma_0 \subseteq V^*(1_{L^2(\T^2)}
\ot \al)V$ and thus,
\begin{equation}
 1_{L^2(\T^2)}  \odot  \al_0 \odot \al_0 + q \, 1_{L^2(\T^2)} \odot e_0 \gamma_0^\dag \odot \gamma_0
 \subseteq (1_{L^2(\T^2)} \ot V^*) (1_{L^2(\T^2)} \ot 1_{L^2(\T^2)}   \ot \al) (1_{L^2(\T^2)} \ot V) \ .
\label{coass.eq1}
\end{equation}
By Lemma \ref{neumann.lem1}, $D(1_{L^2(\T^2)} \ot \al) = D(1_{L^2(\T^2)} \ot e \gamma^*)$. Let $v \in D(1_{L^2(\T^2)}
\ot \al)$ and $w \in E$. There exists a sequence $(v_n)_{n=1}^\infty$ in $L^2(\T^2) \odot E$ such that
$(v_n)_{n=1}^\infty \rightarrow v$ and $\bigl(\,(1_{L^2(\T^2)} \odot \al_0)(v_n)\,\bigr)_{n=1}^\infty \rightarrow
(1_{L^2(\T^2)} \ot \al)(v)$. Since $(1_{L^2(\T^2)}  \odot  \al_0^\dag)(1_{L^2(\T^2)}  \odot  \al_0) = (1_{L^2(\T^2)}
\odot e_0 \gamma_0)(1_{L^2(\T^2)} \odot e_0 \gamma_0^\dag) + 1_{L^2(\T^2)} \odot  e_0$, $1_{L^2(\T^2)} \odot  e_0$ is
bounded  and $1_{L^2(\T^2)} \ot e \gamma^*$ is closed, it follows that $\bigl(\,(1_{L^2(\T^2)} \odot e_0
\gamma_0^\dag)(v_n)\,\bigr)_{n=1}^\infty \rightarrow (1_{L^2(\T^2)} \ot e \gamma^*)(v)$.

Therefore  the net
$$\bigl(\,(1_{L^2(\T^2)} \odot
\al_0 \odot \al_0 + q \, 1_{L^2(\T^2)} \odot e_0 \gamma_0^\dag \odot \gamma_0)(v_n \ot w)\,\bigr)_{n=1}^\infty$$
converges to $(\,(1_{L^2(\T^2)} \ot \al) \odot \al_0 + q \,(1_{L^2(\T^2)} \ot e \gamma^*)\odot \gamma_0\,)(v \ot w)$.
Hence, inclusion (\ref{coass.eq1}) and the fact that the operator on the right hand side of this inclusion is closed
imply that $v \ot w \in D\bigl((1_{L^2(\T^2)} \ot V^*) (1_{L^2(\T^2)} \ot 1_{L^2(\T^2)} \ot  \al) (1_{L^2(\T^2)} \ot
V)\bigr)$ and
$$
(\,(1_{L^2(\T^2)} \ot \al) \odot \al_0 + q \, (1_{L^2(\T^2)} \ot e \gamma^*)\odot \gamma_0)(v \ot w)  = (1_{L^2(\T^2)}
\ot V^*) (1_{L^2(\T^2)} \ot 1_{L^2(\T^2)}   \ot  \al) (1_{L^2(\T^2)} \ot V)(v \ot w) \ .
$$
In other words,
$$
(1_{L^2(\T^2)} \ot \al) \odot \al_0 + q \, (1_{L^2(\T^2)} \ot e \gamma^*)\odot \gamma_0 \subseteq (1_{L^2(\T^2)} \ot
V^*) (1_{L^2(\T^2)} \ot 1_{L^2(\T^2)}   \ot  \al) (1_{L^2(\T^2)} \ot V) \ .
$$
As a consequence, proposition \ref{comult.prop1} implies that
\begin{eqnarray*}
& & \de_0^{(2)}(\al_0) = \de_0(\al_0) \odot \al_0 + q \, \de_0(e_0 \gamma_0^\dag) \odot \gamma_0
 \subseteq V^* (1_{L^2(\T^2)} \ot \al)V \odot \al_0 + q \, V^*(1_{L^2(\T^2)} \ot e \gamma^*)V \odot \gamma_0
\\ & & \spat \subseteq (V^* \ot 1_H)(1_{L^2(\T^2)} \ot V^*) (1_{L^2(\T^2)} \ot 1_{L^2(\T^2)}   \ot  \al) (1_{L^2(\T^2)} \ot V)(V \ot 1_H)
 = (\de \ot \io)\de(\al) \ .
\end{eqnarray*}
\end{proof}

\medskip

For later purposes we introduce some extra notation: if $z,w \in \C$, we set
\begin{equation}
\phi(w;z) = \sum_{r=0}^\infty \frac{(w q^{2r};q^2)_\infty}{(q^2;q^2)_r}\,q^{2r(r-1)}\,z^r \ . \label{coass.eq32}
\end{equation}
Note that this $\phi$-function is closely related to the $_0\vfi_2$-function (see \cite[Def.~(1.2.22)]{Gas}).

\medskip

The next step of the proof consists of showing that

\smallskip

\begin{proposition} \label{coass.prop5}
$(\de \ot \io)\de(\gamma) = (\io \ot \de)\de(\gamma)$.
\end{proposition}

\smallskip

The proof of this result is a tedious but not so difficult task and  is given in the second half of this section, after
Proposition \ref{coass.prop6}. It is however important to remember that this result is needed for the last two parts
(Corollary \ref{coass.cor1} and Proposition \ref{coass.prop6}) of  the proof of the coassociativity of $\de$.

\medskip

Recall that we introduced $M_q^\diamond$ after Definition \ref{neumann.def1}.

\smallskip

\begin{corollary} \label{coass.cor1}
We have that $(\de \ot \io)\de(x) = (\io \ot \de)\de(x)$ for all $x \in M_q^\diamond$.
\end{corollary}
\begin{proof}
The previous proposition implies that $D\bigl((\de \ot \io)\de(\gamma)\bigr) = D\bigl((\io \ot \de)\de(\gamma)\bigr)$ \
. \inlabel{co.eq3}

By Eq.~(\ref{coass.eq2}), we know that $$D\bigl((\de \ot \io)\de(T)\bigr) = (V^* \ot 1_H)(1_{L^2(\T^2)} \ot
V^*)\,D(1_{L^2(\T^2)} \ot 1_{L^2(\T^2)} \ot T)$$ for $T = \gamma$ and $T = \al$. Lemma \ref{neumann.lem1} guarantees
that $D(1_{L^2(\T^2)} \ot 1_{L^2(\T^2)} \ot \gamma) = D(1_{L^2(\T^2)} \ot 1_{L^2(\T^2)} \ot \al)$. As a consequence,
$D\bigl((\de \ot \io)\de(\gamma)\bigr) = D\bigl((\de \ot \io)\de(\al)\bigr)$. In a similar way one shows that
$D\bigl((\io \ot \de)\de(\gamma)\bigr) = D\bigl((\io \ot \de)\de(\al)\bigr)$. By Eq.~(\ref{co.eq3}), this guarantees
that $D\bigl((\de \ot \io)\de(\al)\bigr) = D\bigl((\io \ot \de)\de(\al)\bigr)$.

Proposition \ref{coass.prop2} tells us that $\de_0^{(2)}(\al_0^\dag) \subseteq (\de \ot \io)\de(\al^*)$ and
$\de_0^{(2)}(\al_0^\dag) \subseteq (\io \ot \de)\de(\al^*)$. Taking the adjoints of this inclusions, we see that $(\de
\ot \io)\de(\al)$ and $(\io \ot \de)\de(\al)$ are both restrictions of $\de_0^{(2)}(\al_0^\dag)^*$. The previous
discussion learned us that also the domains of $(\de \ot \io)\de(\al)$ and $(\io \ot \de)\de(\al)$ agree. Thus, $(\de
\ot \io)\de(\al) = (\io \ot \de)\de(\al)$.

\smallskip

Since $M_\ze \ot 1$ is the unitary in the polar decomposition of $\gamma$ and $(\de \ot \io)\de(\gamma) = (\io \ot
\de)\de(\gamma)$, the uniqueness of the polar decomposition implies that $(\de \ot \io)\de(M_\ze \ot 1) = (\io \ot
\de)\de(M_\ze \ot 1)$ (*) and $(\de \ot \io)\de(|\gamma|) = (\io \ot \de)\de(|\gamma|)$. Since $\de(e) = e \ot e$,
$(\de \ot \io)\de(e) = (\io \ot \de)\de(e)$, so we get that $(\de \ot \io)\de(e\,|\gamma|) = (\io \ot
\de)\de(e\,|\gamma|)$. Applying functional calculus to this equation, we see that $(\de \ot \io)\de(M_{1 \ot g}) = (\io
\ot \de)\de(M_{1 \ot g})$ for all $g \in \cL^\infty(I_q)$. Combining this with (*), it follows that $(\de \ot
\io)\de(\pi(f)) = (\io \ot \de)\de(\pi(f))$ for all $f \in \cL^\infty(\T \times I_q)$ (cf. the proof of Lemma
\ref{neumann.lem2}).

Because $w$ is the isometry in the polar decomposition of $\al$ and $(\de \ot \io)\de(\al) = (\io \ot \de)\de(\al)$,
the uniqueness of the polar decomposition implies that $(\de \ot \io)\de(w) = (\io \ot \de)\de(w)$. The lemma follows
by Lemma \ref{neumann.lem2}.
\end{proof}

\medskip

In order to finalize the proof of the coassociativity we rely on formulas discovered by S.L. Woronowicz. In particular,
we borrow the use of the extended Hilbert space,  the definition of $\tilde{\Theta}$ as in (\ref{coass.eq30}) and
Eq.~(\ref{co.eq1}) in the proof below. These formulas can be found in \cite[Sec.~5]{Wor7} \ (Eq.~(\ref{co.eq1}) is a
slight variation of the one in \cite{Wor7}, where $e$ is replaced by 1). S.L. Woronowicz told us that he has a proof
for these identities, not included in the version of \cite{Wor7} yet. Our proof of the strong convergence in
(\ref{coass.eq30}) and Eq.~(\ref{co.eq1}) are probably different from his ones, due to the difference of both
approaches. The proof of the equality $(\de \ot \io)\de(u) = (\io \ot \de)\de(u)$ is not present in \cite{Wor7}, but it
is possible that S.L. Woronowicz also knows a proof for this. However, the property of tensor products of
$\widetilde{SU}_q(1,1)$-quadruples that corresponds to the equality $D((\de \ot \io)\de(\gamma)) = D((\io \ot
\de)\de(\gamma))$ is still not proven within the approach of \cite{Wor7}.

\smallskip

\begin{proposition} \label{coass.prop6}
The $^*$-homomorphism $\de : M_q \rightarrow M_q \ot M_q$ is coassociative.
\end{proposition}
\begin{proof}
Since $M_q$ is generated by $M_q^\diamond \cup \{u\}$, it only remains to show that $(\de \ot \io)\de(u) = (\io \ot
\de)\de(u)$ by the previous corollary. Define first of all the isometry $v \in M_q^\diamond$ as $v = e\,w\,(M_\ze^* \ot
1)$.

Set $K_q = - q^{\Z} \cup q^{\Z}$ and use the uniform measure on $K_q$. Define $\tilde{H} = L^2(\T) \ot L^2(K_q) $. Then
$\tilde{H}$ is obviously the orthogonal sum of $H$ and $L^2(\T) \ot L^2(K_q \setminus I_q)$.

Let $L$ be a Hilbert space. Then we define a normal $^*$-homomorphism $B(H \ot L) \rightarrow B(\tilde{H} \ot L) : a
\rightarrow \hat{a}$ such that $\hat{a}\!\!\restriction_{H \ot L} = a$ and $\hat{a}\!\!\restriction_{(H \ot L) ^\perp}
= 0$. If $L_1,L_2$ are two Hilbert spaces, $(b \ot c)\hoed  = \hat{b} \ot c$ for all $b \in B(H \ot L_1)$, $c \in
B(L_2)$. \inlabel{coass.eq33}

On $\tilde{H}$ we define unitary operators $\tilde{e}$, $\tilde{w}$ and $\tilde{u}$ such that for $f \in L^2(\T \times
K_q)$,
$$(\tilde{e} f)(\lambda,x) = \sgn(x)\, f(\lambda,x) \hspace{2cm} (\tilde{w} f)(\lambda,x) = f(\lambda,q^{-1} x)
\hspace{2cm} (\tilde{u} f)(\lambda,x) = f(\lambda,-x) $$ for almost all $(\lambda,x) \in \T \times K_q$. Now define the
unitary element $\tilde{v}$ in $B(\tilde{H})$ as $\tilde{v} = \tilde{e}\,\tilde{w}\,(M_\ze^* \ot 1)$.

\smallskip

Since $(L^2(\T) \ot L^2(K_q \setminus I_q)) \ot H$ is separable, we can extend the orthonormal basis for $H \ot H$
given by the functions of Definition \ref{comult.def2} to an orthonormal basis $(\,\digamma_{r,s,m,p} \mid r,s,m \in
\Z, p \in K_q \,)$ for $\tilde{H} \ot H$.  If $p \in K_q$, we define $\tilde{\sde}_p \in L^2(K_q)$ such that
$\tilde{\sde}_p(x) = \sde_{x,p}$ for all $x \in K_q$. Now define the unitary operator $\tilde{V} : \tilde{H} \ot H
\rightarrow L^2(\T^2) \ot \tilde{H}$ such that $\tilde{V} \digamma_{r,s,m,p} = \ze^r \ot \ze^s \ot \ze^m \ot
\tilde{\sde}_p$ for all $r,s,m \in \Z$ and $p \in K_q$. So the restriction of $\tilde{V}$ to $H \ot H$ equals $V$ and
$\tilde{V}(H^\perp \ot H) = L^2(\T^2) \ot H^\perp$.

Define the normal injective $^*$-homomorphism $\tde : B(\tilde{H}) \rightarrow B(\tilde{H} \ot H)$ such that $\tde(a) =
\tilde{V}^* (1_{L^2(\T^2)}  \ot a) \tilde{V}$ for all $a \in B(\tilde{H})$. Thus, if $a \in M_q$, then $\tde(\hat{a}) =
\de(a)\hoed$. By Eq.~(\ref{coass.eq33}), this implies that $(\tde \ot \io)(\hat{b}) = (\de \ot \io)(b)\hoed$ for all $b
\in M_q \ot B(H)$.

\smallskip

Take $r,s,m \in \Z$ and $p \in K_q$. Choose $n \in \NE$.  One easily checks that  $\tde(\tilde{v}) \digamma_{r,s,m,p} =
\sgn(p)\,\digamma_{r,s,m-1,pq}$. Thus, $\tde(\tilde{v})^n \digamma_{r,s,m,p} = \sgn(p)^n\, \digamma_{r,s,m-n,p q^n}$.
Notice that $p\,q^n \in I_q$ if $n \geq 1 - \chi(p)$.

Assume that $n \geq 1 - \chi(p)$. Choose $x \in K_q$, $y \in I_q$ and $\lambda, \mu \in \T$. If $q^n x \in I_q$ and
$y/x = \sgn(p)\,q^m$, then $y/(q^n x) = \sgn(p q^n)\,q^{m-n}$ which by Definition \ref{comult.def2} and Result
\ref{app.res2} guarantees that
\begin{eqnarray}
& & \bigl((\tilde{v}^{-n} \ot 1)\tde(\tilde{v})^n \digamma_{r,s,m,p}\bigr)(\lambda,x,\mu,y) = \sgn(p)^n \,
[(\tilde{v}^{-n} \ot 1) \digamma_{r,s,m-n,p q^n}](\lambda,x,\mu,y) \nonumber
\\ & & \spat = \sgn(p)^n \,\sgn(x)^n\,\lambda^n \, \digamma_{r,s,m-n,p q^n}(\lambda,q^n x,\mu,y) \nonumber
\\ & & \spat = \sgn(y)^n \, \lambda^n\, a_{q^n p}(q^n x,y) \, \lambda^{r+\chi(y/p q^n)}\, \mu^{s-\chi(x q^n/ p q^n)} \nonumber
\\ & & \spat = \sgn(y)^n \, a_{q^n p}(q^n x,y) \, \lambda^{r+\chi(y/p)}\, \mu^{s-\chi(x/p)} \nonumber
\\ & & \spat = \lambda^{r+\chi(y/p)}\, \mu^{s-\chi(x/p)}\, \sgn(y)^n \, c_q\,s(q^n x,y)\, (-1)^{\chi(q^n p)}\,(-\sgn(y))^{\chi(q^n x)}\,|y|\,\nu(q^n p y/q^n
x)\,\nonumber
\\ & & \spat \ \ \  \sqrt{\frac{(-\kappa(q^n p),-\kappa(y);q^2)_\infty}{(-\kappa(q^n x);q^2)_\infty}}\,\,\Psis{-q^2/\kappa(q^n p)}{
q^2 \kappa(q^n x/q^n p)}{q^2 \kappa(q^n x/y)} \nonumber
\\ & & \spat = \lambda^{r+\chi(y/p)}\,
\mu^{s-\chi(x/p)}\,c_q\,s(x,y)\, (-1)^{\chi(p)}\,(-\sgn(y))^{\chi(x)}\,|y|\,\nu(p y/x)\,\nonumber
\\ & & \spat \ \ \  \sqrt{\frac{(-\kappa(q^n p),-\kappa(y);q^2)_\infty}{(-\kappa(q^n x);q^2)_\infty}}\,\,\Psis{-q^2/\kappa(q^n p)}{
q^2 \kappa(x/p)}{q^2 \kappa(q^n x/y)} \nonumber
\end{eqnarray}
\begin{eqnarray}
& & \spat = \lambda^{r+\chi(y/p)}\, \mu^{s-\chi(x/p)}\,c_q\,s(x,y)\,
(-1)^{\chi(p)}\,(-\sgn(y))^{\chi(x)}\,|x|\,q^m\,\nu(q^m p)\,\nonumber
\\ & & \spat \ \ \  \sqrt{\frac{(-\kappa(q^n p),-\sgn(p)\,q^{2m}\,\kappa(x);q^2)_\infty}{(-\kappa(q^n x);q^2)_\infty}}\,\,\Psis{-q^2/\kappa(q^n p)}{
q^2 \kappa(x/p)}{\sgn(p)\,q^{2(1-m+n)}} \ . \nonumber
\\ & &  \label{coass.eq25}
\end{eqnarray}

If, on the other hand, $q^n x \not\in I_q$ or  $y/x \not= \sgn(p)\,q^m$, then
\begin{eqnarray*}
\bigl((\tilde{v}^{-n} \ot 1)\tde(\tilde{v})^n \digamma_{r,s,m,p}\bigr)(\lambda,x,\mu,y) & = & \sgn(p)^n \,
[(\tilde{v}^{-n} \ot 1) \digamma_{r,s,m-n,p q^n}](\lambda,x,\mu,y)
\\ & = & \sgn(p)^n \,\sgn(x)^n\,\lambda^n \, \digamma_{r,s,m-n,p q^n}(\lambda,q^n x,\mu,y) = 0 \ .
\end{eqnarray*}
Using notation (\ref{coass.eq32}), we define the function $f_{r,s,m,p} : \T \times K_q \times \T \times I_q \rightarrow
\C$ such that for $x \in K_q$, $y \in I_q$, $\lambda,\mu \in \T$, $f_{r,s,m,p}(\lambda,x,\mu,y)$ equals
\begin{equation} \label{coass.eq31}
\lambda^{r+\chi(y/p)}\, \mu^{s-\chi(x/p)}\,c_q\,s(x,y)\, (-1)^{\chi(p)}\,(-\sgn(y))^{\chi(x)}\,|y|\,\nu(q^m p)\,
\sqrt{(-\kappa(y);q^2)_\infty} \,\,\phi(q^2 \kappa(x/p);- q^{2(2-m)}/p^2)
\end{equation}
if $y/x = \sgn(p)\,q^m$ and $f_{r,s,m,p}(\lambda,x,\mu,y) = 0$ if $y/x \not= \sgn(p)\,q^m$. Then it is clear from the
above computations (and the proof of Lemma \ref{app.lem2}) that $\bigl(\,(\tilde{v}^{-n} \ot 1)\tde(\tilde{v})^n
\digamma_{r,s,m,p}\,\bigr)_{n=1}^\infty$ converges pointwise to $f_{r,s,m,p}$. \!\!\inlabel{coass.eq24}

\medskip

Define $J = \{\, x \in K_q \mid |x| \leq q^{1-m} \text{ or } \sgn(x) = \sgn(p)\,\}$. By Lemma \ref{alt.lem1} there
exists a function $H \in l^2(J)^+$  such that
$$|\,\,|x|\,\sqrt{(-\sgn(p)\,q^{2m}\,\kappa(x);q^2)_\infty}\,\,\Psi(-q^2/\kappa(q^n p);q^2\kappa(x/p);\sgn(p)
q^{2(1-m+n)})\,| \leq H(x)$$ for all $x \in J$ and $n \in \NE$.

\medskip

There exists clearly a number $C > 0$ such that $q^m\,\nu(q^m p) \, (-\kappa(q^n
p);q^2)_\infty^{\frac{1}{2}}\,(q^2;q^2)_\infty^{-\frac{1}{2}} \leq C$ for all $n \in \NE$.

Let $(x,y) \in K_q \times I_q$ such that $y/x = \sgn(p) q^m$. If $|x| > q^{1-m}$ then $|q^m x| \geq 1$ and since
$\sgn(p) q^m x = y \in I_q$, this implies that $\sgn(x) = \sgn(p)$. Hence, $x$ belongs to $J$. So we can define an
element $G \in \cL^2(\T \times K_q \times \T \times I_q)$ such that for $\lambda,\mu \in \T$ and $x \in K_q$, $y \in
I_q$, we have that $G(\lambda,x,\mu,y) = 0$ if $y/x \not= \sgn(p) q^m$ and $G(\lambda,x,\mu,y) = C\,H(x)$ if $y/x =
\sgn(p) q^m$. Then Eq.~(\ref{coass.eq25}) and the remarks thereafter imply that $|(\tilde{v}^{-n} \ot
1)\tde(\tilde{v})^n \digamma_{r,s,m,p}| \leq G$ for all $n \in \NE$ such that $n \geq 1-\chi(p)$. Therefore the
convergence in (\ref{coass.eq24}) and  the dominated convergence theorem imply that $f_{r,s,m,p}$ belongs to $\tilde{H}
\ot H$ and that the sequence  $\bigl(\,(\tilde{v}^{-n} \ot 1)\tde(\tilde{v})^n \digamma_{r,s,m,p}\,\bigr)_{n=1}^\infty$
converges to $f_{r,s,m,p}$ in $\tilde{H} \ot H$.

\medskip

Since the sequence $\bigl(\,(\tilde{v}^{-n} \ot 1)\tde(\tilde{v})^n\,\bigr)_{n=1}^\infty$ is a sequence of isometries
and the linear span of such elements $\digamma_{r,s,m,p}$ is dense in $\tilde{H} \ot H$, we conclude that there exists
an isometry $\tilde{\Theta} : \tilde{H} \ot H \rightarrow \tilde{H} \ot H$ such that the sequence
$\bigl(\,(\tilde{v}^{-n} \ot 1)\tde(\tilde{v})^n\,\bigr)_{n=1}^\infty$ converges strongly to $\tilde{\Theta}$.
\inlabel{coass.eq30}

Thus, $\tilde{\Theta} \digamma_{r,s,m,p} = f_{r,s,m,p}$ for all $r,s,m \in \Z$ and $p \in K_q$. Using
Eq.~(\ref{coass.eq31}) and the definition of $\tde$ and $\tilde{u}$, one easily checks that $(\tilde{u} \ot e)
f_{r,s,m,p} = f_{r,s,m,-p}$ and $\tde(\tilde{u}) \digamma_{r,s,m,p} = \digamma_{r,s,m,-p}$ for all $r,s,m \in \Z$ and
$p \in K_q$ (note that $s(-x,y) = \sgn(y)\,s(x,y)$ for all $x,y \in \R \setminus \{0\}$). \newline As a consequence,
$\tilde{\Theta} \tde(\tilde{u}) = (\tilde{u} \ot e) \tilde{\Theta}$. \hspace{-5ex} \inlabel{co.eq1}

\medskip

Define $\hat{\Theta} \in B(\tilde{H} \ot H)$ such that $\hat{\Theta}\!\!\restriction_{H \ot H} =
\tilde{\Theta}\!\!\restriction_{H \ot H}$ and $\hat{\Theta}\!\!\restriction_{H^\perp \ot H} = 0$. One sees that $\de(v)
\digamma_{r,s,m,p} = \sgn(p) \digamma_{r,s,m-1,qp} = \tde(\tilde{v}) \digamma_{r,s,m,p}$ for all $r,s,m \in \Z$, $p \in
I_q$, implying that $\de(v)$ is the restriction of $\tde(\tilde{v})$ to $H \ot H$. So we get for every $n \in \N$ that
$\tilde{\de}(\hat{v})^n\!\!\restriction_{H \ot H} = \bigl(\de(v)\hoed\,\bigr)^n\!\!\restriction_{H \ot H} =
\tde(\tilde{v})^n\!\!\restriction_{H \ot H}$, whereas $\tilde{\de}(\hat{v})^n\!\!\restriction_{H^\perp \ot H} = 0$. It
follows that $\hat{\Theta}$ is the strong limit of the sequence of isometries $\bigl(\,(\tilde{v}^{-n} \ot
1)\tilde{\de}(\hat{v})^n\bigl)_{n=1}^\infty$. As such, $\hat{\Theta}$ belongs to $B(\tilde{H}) \ot M_q$.

 We set
$P = u^* u$, which belongs to $M_q^\diamond$ because it is a multiplication operator. Multiplying Eq.~(\ref{co.eq1})
from the right by  $\tde(\hat{P})$ and using the fact that $\tilde{u}\,\hat{P} = \hat{u}$, we see that $\tilde{\Theta}
\tde(\hat{u}) = (\tilde{u} \ot e) \tilde{\Theta} \tde(\hat{P})$. Thus, $\hat{\Theta} \de(u)\hoed = (\tilde{u} \ot e)
\hat{\Theta} \de(P)\hoed$. \inlabel{co.eq2}

Restricting the previous equation to $H \ot H$, we get $\hat{\Theta} \de(u) = (\tilde{u} \ot e) \hat{\Theta} \de(P)$.
Applying $\io \ot \de$ to this equality, it follows that $(\io \ot \de)(\hat{\Theta})\, (\io \ot \de)(\de(u)) =
(\tilde{u} \ot e \ot e) (\io \ot \de)(\hat{\Theta}) (\io \ot \de)\de(P)$. \inlabel{coass.eq35}

If we apply $\tde \ot \io$ to Eq.~(\ref{co.eq2}), we get  $(\tde \ot \io)(\hat{\Theta}) (\de \ot \io)(\de(u))\hoed =
(\tde(\tilde{u}) \ot e) (\tde \ot \io)(\hat{\Theta}) (\de \ot \io)(\de(P))\hoed$. So the restriction of this equation
to $H \ot H \ot H$ gives $(\tde \ot \io)(\hat{\Theta}) (\de \ot \io)\de(u) = (\tde(\tilde{u}) \ot e) (\tde \ot
\io)(\hat{\Theta}) (\de \ot \io)\de(P)$.  Multiplying this equation from the left by $\tilde{\Theta} \ot 1$, we see
that
\begin{eqnarray}
& & (\tilde{\Theta} \ot 1)\,(\tde \ot \io)(\hat{\Theta}) (\de \ot \io)\de(u) = (\tilde{\Theta} \tde(\tilde{u}) \ot e)
(\tde \ot \io)(\hat{\Theta}) (\de \ot \io)\de(P) \nonumber
\\ & & \spat = (\tilde{u} \ot e \ot e)(\tilde{\Theta} \ot 1) (\tde \ot \io)(\hat{\Theta}) (\io \ot \de)\de(P)
\label{coass.eq34}
\end{eqnarray}
where in the last equality we used Eq.~(\ref{co.eq1}), the fact that $P \in M_q^\diamond$ and Corollary
\ref{coass.cor1}.

Now we calculate $(\tilde{\Theta} \ot 1)\,(\tde \ot \io)(\hat{\Theta})$. Let $n \in \N$. Since $v$ belongs to
$M_q^\diamond$, Corollary \ref{coass.cor1} guarantees that $(\de \ot \io)\de(v) = (\io \ot \de)\de(v)$ which implies
that $(\tde \ot \io)\tde(\hat{v}) = (\io \ot \de)\tilde{\de}(\hat{v})$. Hence,
\begin{eqnarray*}
(\,[(\tilde{v}^{-n} \ot 1) \tde(\tilde{v}^n)] \ot 1) (\tde \ot \io)\bigl(\, (\tilde{v}^{-n} \ot 1)
\tde(\hat{v}^n)\,\bigr) & = & (\tilde{v}^{-n} \ot 1 \ot 1)\, (\tde \ot \io)\tde(\hat{v}^n)
\\ & = & (\tilde{v}^{-n} \ot 1 \ot 1)\, (\io \ot \de)\tde(\hat{v}^n) \ .
\end{eqnarray*}
If we let $n$ tend to infinity in this equality, the expression we started with converges strongly to
\newline $(\tilde{\Theta} \ot 1)(\tde \ot \io)(\hat{\Theta})$ whereas the expression we end up converges strongly to $(\io
\ot \de)(\hat{\Theta})$. Thus, \newline $(\tilde{\Theta} \ot 1)(\tde \ot \io)(\hat{\Theta})$ $= (\io \ot
\de)(\hat{\Theta})$. Inserting this in Eq.~(\ref{coass.eq34}), we conclude that $(\io \ot \de)(\hat{\Theta})\,(\de \ot
\io)\de(u) = (\tilde{u} \ot e \ot e)(\io \ot \de) (\hat{\Theta})$ $(\io \ot \de)\de(P)$. Comparing this with
Eq.~(\ref{coass.eq35}), we conclude that
\begin{equation} \label{coass.eq8}
(\io \ot \de)(\hat{\Theta})\,(\io \ot \de)\de(u) =  (\io \ot \de)(\hat{\Theta})\,(\de \ot \io)\de(u) \ .
\end{equation}
The initial projection of $\hat{\Theta}$ is given by $X \ot 1$ where $X \in B(\tilde{H})$ is the orthogonal projection
onto $H$. In other words, $\hat{\Theta}^* \hat{\Theta} = X \ot 1$. Thus, $(\io \ot \de)(\hat{\Theta})^* (\io \ot
\de)(\hat{\Theta}) = X \ot 1 \ot 1$, implying that $(\io \ot \de)(\hat{\Theta})$ is a partial isometry with initial
space $H \ot H \ot H$. Therefore Eq.~(\ref{coass.eq8}) guarantees that $(\io \ot \de)\de(u) = (\de \ot \io)\de(u)$.
\end{proof}

\medskip

It is clear from the end of the proof above that we need $\hat{\Theta}$ to be an isometry on $H \ot H \ot H$. This is
precisely the reason for using the extended Hilbert space $\tilde{H}$, a definition of $\hat{\Theta}$ as the strong
limit of $\bigl(\,((v^*)^n \ot 1)\de(v)^n\,\bigr)_{n=1}^\infty$ would not result in an isometry. The need for using
$\tilde{\Theta}$ should be clear from the fact that the range of the operator $(\tde \ot \io)(\hat{\Theta}) (\de \ot
\io)\de(u)$ in Eq.~(\ref{coass.eq34}) is not contained in $H \ot H \ot H$.

\bigskip\medskip

In the last part of this section we provide the proof of Proposition \ref{coass.prop5}. The most intricate part of this
proof is contained in Proposition \ref{coass.prop3}. The results before this proposition reduce the problem to a
simpler form. We stress that this exposition does not depend on Corollary \ref{coass.cor1} or Proposition
\ref{coass.prop6}.

\smallskip

In order to filter out some distracting features we will change our Hilbert space $H^{\ot 3}$. For this purpose, define
a measure $\Omega$ on $I_{q^2}^3$ such that $$\Omega(\{(x,z,y)\}) = |x|\,|y|\,\frac{(-x,-y;q^2)_\infty
}{(-z;q^2)_\infty}$$ for all $x,y,z \in I_{q^2}$. In this section, all $L^2$-spaces that involve $I_{q^2}^3$ are
defined by this measure.

\smallskip

Define the unitary transformation $U : H^{\ot 3} \rightarrow L^2(\T^3) \ot L^2(I_{q^2}^3,\Omega)$ such that $U\,[f] =
[g]$, where the functions $f \in \cL^2((\T \times I_q)^3)$ and $g \in \cL^2(\T^3 \times I_{q^2}^3)$ are such that
\begin{eqnarray*}   g(\lambda,\mu,\eta,x,z,y) & = &
s(xz,y)\,s(z,x)\,\lambda^{-\log_{q^2}\!|x|}\,\mu^{\log_{q^2}\!|x/y|}\,\eta
^{\log_{q^2}\!|y|}\,(-1)^{\log_{q^2}\!|z|}\,\sgn(z)^{\log_{q^2}\!|x|} \, \sgn(y)^{\log_{q^2}\!|xz|}\,
\\ &  & |x|^{-1/2} \,  |y|^{-1/2}\,
\sqrt{\frac{(-z;q^2)_\infty}{(-x,-y;q^2)_\infty }}\,\,f(\lambda,\kappa^{-1}(x),\mu,\kappa^{-1}(z),\eta,\kappa^{-1}(y))
\end{eqnarray*}
for all $\lambda,\mu,\eta \in \T$ and $x,z,y \in I_{q^2}$.

\medskip\medskip

If $\theta \in -q^{2\Z} \cup q^{2\Z}$,  we define the sets (see Fig.~1 and Fig.~2 in the proof of proposition
\ref{coass.prop3})

$$K(\theta) = \{\,(x,y) \mid x,y \in I_{q^2} \text{ such that } \theta x y \in I_{q^2}\,\}\hspace{3ex} \text{and}
\hspace{3ex} L(\theta) = \{\,(x,\theta x y,y)\, \mid (x,y) \in K(\theta)\,\} \ .$$
\smallskip

Consider $m,k \in \Z$ and $p \in I_{q^2}$. Define the function $F^k_{m,p} \in  \cF(K(\sgn(p)q^{2k}))$ such that for
$(x,y) \in K(\sgn(p)q^{2k})$,
\begin{eqnarray}
F^k_{m,p}(x,y) & = & c_q^2\,q^{\frac{1}{2}(m+k+1)(m+k+2)}\, \nu(\kappa^{-1}(p)q^m)\,\,\sqrt{(-p;q^2)_\infty} \nonumber
\\ & & \Psis{-q^2/p}{q^{2(1-m)}\,y/|p|}{\sgn(p)\,q^{2(1-m)}}\,\,\Psis{-q^2/x}{\sgn(p)\,q^{2(1+k)}
y}{q^{2(1+m+k)}x} \label{coass.eq16}
\end{eqnarray}
if $\sgn(p)q^{-2m} y \in I_{q^2}$ and $F^k_{m,p}(x,y) = 0$ if $\sgn(p) q^{-2m} y \not\in I_{q^2}$. We also define the
function $G^k_{m,p} \in \cF(K(\sgn(p)q^{2k}))$ by $G^k_{m,p}(x,y) = F^k_{-m,p}(y,x)$ for all $(x,y) \in
K(\sgn(p)q^{2k})$.

\medskip

Next we define functions $\tilde{F}^k_{m,p}, \tilde{G}^k_{m,p} \in  \cF(I_{q^2}^3)$ both of which have support in
$L(\sgn(p)q^{2k})$ and satisfy
$$\tilde{F}^k_{m,p}(x,\sgn(p)q^{2k} x y,y) = F^k_{m,p}(x,y) \hspace{10ex} \text{and} \hspace{10ex} \tilde{G}^k_{m,p}(x,\sgn(p)q^{2k} x y,y) =
G^k_{m,p}(x,y)$$ for all $(x,y) \in K(\sgn(p)q^{2k})$.

\medskip

Define the unitary element $\vsi \in C(\T^3)$ such that $\vsi(\lambda,\eta,\mu) = \lambda \mu \bar{\eta}$ for all
$\lambda,\eta,\mu \in \T$.

\smallskip

\begin{lemma} \label{coass.lem1}
The families $\{\,\tilde{F}^k_{m,p} \mid p \in I_{q^2}, m,k \in \Z \,\}$ and $\{\,\tilde{G}^k_{m,p} \mid p \in I_{q^2},
m,k \in \Z \,\}$ are both orthonormal bases of $L^2(I_{q^2}^3,\Om)$. Moreover,
\begin{trivlist}
\item[\ \,1)] the space $\langle \, v \ot \tilde{F}^k_{m,p} \mid v \in L^2(\T^3), p \in I_{q^2}, m,k \in \Z \,\rangle$  is a core
for $U\,(\de \ot \io)(\de(\gamma))\,U^*$,
\smallskip
\item[\ \,2)] the space $\langle\,v \ot  \tilde{G}^k_{m,p} \mid v \in L^2(\T^3), p \in I_{q^2}, m,k \in \Z \,\rangle$  is a core
for $U\,(\io \ot \de)(\de(\gamma))\,U^*$,
\smallskip
\item[\ \,3)] if $v \in L^2(\T^3)$, $p \in I_{q^2}$ and $m,k \in \Z$, then
\begin{eqnarray*}
(U(\de \ot \io)(\de(\gamma))U^*) \,(v \ot \tilde{F}^k_{m,p})  & = & |p|^{-\frac{1}{2}}\, M_\vsi v \ot
\tilde{F}^k_{m+1,p}
\\ (U(\io \ot \de)(\de(\gamma))U^*) \,(v \ot \tilde{G}^k_{m,p})  & = &  |p|^{-\frac{1}{2}}\, M_\vsi v \ot
\tilde{G}^k_{m+1,p}\ .
\end{eqnarray*}
\end{trivlist}
\end{lemma}
\begin{proof}
The proof of this fact is pretty straightforward. For $r,s,t,n,m \in \Z$ and $p \in I_q$, we define $F^{r,s,t,n}_{m,p},
G^{r,s,t,n}_{m,p} \in H^{\ot 3}$ as
\begin{eqnarray*}
F^{r,s,t,n}_{m,p} & = & (V^* \ot 1_H)(1_{L^2(\T^2)} \ot V^*)(\ze^r \ot \ze^s \ot \ze^t \ot \ze^n \ot \ze^m \ot \sde_p)
\\ G^{r,s,t,n}_{m,p} & = & (1_H \ot V^*) V_{13}^*(\ze^r \ot \ze^s \ot \ze^t \ot \ze^n \ot \ze^m
\ot \sde_p) \ .
\end{eqnarray*}
From Proposoition \ref{comult.prop7}, we get immediately  that the families $(\,F^{r,s,t,n}_{m,p} \mid r,s,t,n,m \in
\Z, p \in I_q \,)$ and $(\,G^{r,s,t,n}_{m,p} \mid r,s,t,n,m \in \Z, p \in I_q \,)$ are both orthonormal bases of
$H^{\ot 3}$. Moreover, Eqs.~(\ref{neumann.eq1}), (\ref{coass.eq2}) and (\ref{coass.eq3}) imply that
\begin{trivlist}
\item[\ \,1)] the space $\langle\,F^{r,s,t,n}_{m,p} \mid r,s,t,n,m \in \Z, p \in I_q \rangle$  is a core
for the operator $(\de \ot \io)(\de(\gamma))$,
\smallskip
\item[\ \,2)] the space $\langle\,G^{r,s,t,n}_{m,p} \mid r,s,t,n,m \in \Z, p \in I_q \rangle$  is a core
for the operator $(\io \ot \de)(\de(\gamma))$,
\smallskip
\item[\ \,3)] if $r,s,t,m,n \in \Z$ and $p \in I_q$, then
\begin{eqnarray*}
(\de \ot \io)(\de(\gamma)) \,F^{r,s,t,n}_{m,p}  =  p^{-1}\, F^{r,s,t,n}_{m+1,p} \hspace{5ex} \text{and} \hspace{5ex}
(\io \ot \de)(\de(\gamma)) \,G^{r,s,t,n}_{m,p}  =  p^{-1}\, G^{r,s,t,n}_{m+1,p} \ .
\end{eqnarray*}
\hspace{96ex} \inlabel{coass.eq132}
\end{trivlist}
By Definition \ref{comult.def2}, we get the following equalities (in which the infinite sums are $L^2$-convergent),
\begin{eqnarray}
F^{r,s,t,n}_{m,p}  & = & (V^* \ot 1_H)(1_{L^2(\T^2)} \ot V^*) (\ze^r \ot \ze^s \ot \ze^t \ot \ze^n \ot \ze^m \ot
\sde_p) \nonumber
 =  (V^* \ot 1_H)(\ze^r \ot \ze^s \ot \digamma_{t,n,m,p}) \nonumber
\\ & = & \sum_{\scriptsize \begin{array}{c}  y \in I_q \\  \sgn(p)q^{-m} y \in I_q\end{array}} \ a_p(\sgn(p)q^{-m}y,y)\ (V^* \ot
1_H)(\ze^r \ot \ze^s \ot  \begin{tabular}[t]{l} $\ze^{t+\chi(y/p)} \ot$ \\ $ \hspace{-1ex} \sde_{\sgn(p)q^{-m}y} \ot
\ze^{n-\chi(q^{-m}y/p)} \ot \sde_y)$\end{tabular} \nonumber
\\ & = & \sum_{\scriptsize \begin{array}{c}  y \in I_q \\  \sgn(p)q^{-m} y \in I_q\end{array}}
\ a_p(\sgn(p)q^{-m}y,y)\ \digamma_{r,s,t+\chi(y/p),\sgn(p)q^{-m}y} \ot \ze^{n+m+\chi(p/y)} \ot \sde_y \nonumber
\\ & = & \sum_{\scriptsize \begin{array}{c}  y \in I_q \\  \sgn(p)q^{-m} y \in I_q\end{array}}
\sum_{\scriptsize \begin{array}{c}  x \in I_q \\  \sgn(p)\sgn(y) q^{t+\chi(y/p)} x \in I_q\end{array}} \
a_{\sgn(p)q^{-m}y}(x,\sgn(p)\sgn(y) q^{t+\chi(y/p)} x) \nonumber
\\ & & \nonumber
\\ & & \hspace{5ex} a_p(\sgn(p)q^{-m}y,y)\,\ze^{r+t+\chi(y/p)+\chi(x)+m-\chi(y)} \nonumber
\ot \sde_x \ot \ze^{s-\chi(x) -m + \chi(y)} \ot \nonumber
\\ & & \hspace{5ex} \sde_{\sgn(p)\sgn(y) q^{t+\chi(y/p)} x} \ot \ze^{n+m+\chi(p/y)} \ot
\sde_y \nonumber
\\ & & \nonumber
\\ & = & \sum_{\scriptsize \begin{array}{c}  y \in I_q \\  \sgn(p)q^{-m} y \in I_q\end{array}}
\sum_{\scriptsize \begin{array}{c}  x \in I_q \\   (q^t/p)\,x y  \in I_q\end{array}} \
a_p(\sgn(p)q^{-m}y,y)\,a_{\sgn(p)q^{-m}y}(x,(q^t/p)\,x y ) \nonumber
\\ & & \nonumber
\\ & & \hspace{5ex}\ze^{r+t+m+\chi(x/p)}
\ot \sde_x \ot \ze^{s-m+\chi(y/x)} \ot \sde_{(q^t/p)\,x y } \ot \ze^{n+m+\chi(p/y)} \ot \sde_y \ . \label{coass.eq4}
\end{eqnarray}

Choose $\lambda,\mu,\eta \in \T$ and $x,z,y \in I_{q^2}$. Look first at the case where $\sgn(p)\,q^{-2m}\,y \in
I_{q^2}$ and $z = \sgn(p)\,q^{2(t-\chi(p))}\,x y$. Set $x' = \kappa^{-1}(x)$, $y' = \kappa^{-1}(y)$ and $z' =
\kappa^{-1}(z)$. Thus, $\sgn(p) q^{-m}\,y'  \in I_q$ and $(q^t/p) x' y' = z' \in I_q$. Set $d =
\lambda^{r+t+m+\chi(x'/p)}\,\mu^{s-m+\chi(y'/x')}\,\eta^{n+m+\chi(p/y')}$. The above equation,
 Proposition \ref{comult.prop2} and Definition \ref{comult.def1} imply that
\begin{eqnarray*}
& & F^{r,s,t,n}_{m,p}(\lambda,\kappa^{-1}(x),\mu,\kappa^{-1}(z),\eta,\kappa^{-1}(y)) \, =
d\,\,a_p(\sgn(p)q^{-m}y',y')\,a_{\sgn(p)q^{-m}y'}(x',z')\,
\\ & & \spat = d\, \,(\sgn(p) \sgn(y'))^{m+\chi(y')}\,\sgn(z')^{\chi(z')}\,\sgn(x')^{\chi(x')}
\, a_p(\sgn(p)q^{-m}y',y')\, a_{\sgn(p)q^{-m}y'}(z',x')
\\ & & \spat = d\,\,(\sgn(p)\,\sgn(y'))^{m+\chi(y')} \,\sgn(z')^{\chi(x')}\,(\sgn(p)\,\sgn(x')\,\sgn(y'))^{\chi(p)+\chi(y')+t} \,
\,\sgn(x')^{\chi(x')}
\\ & & \spat \ \ \ \ c_q \, s(\sgn(p)q^{-m}y',y')\,  (-1)^{\chi(p)+m+\chi(y')}\, \sgn(y')^{m+\chi(y')} \, |y'|
\\ & & \spat \ \ \ \ \nu(p q^m)\,\sqrt{\frac{(-\kappa(p),-y;q^2)_\infty}{(-\sgn(p)q^{-2m}y;q^2)_\infty}}
 \ \Psis{-q^2/\kappa(p)}{q^{2(1-m)}\,y/|\kappa(p)|}{\sgn(p)\,q^{2(1-m)}}
\\ & & \spat \ \ \ \ c_q \, s(z',x')\,  (-1)^{m+\chi(y')+\chi(z')}\, \sgn(x')^{\chi(z')} \, |x'|
\\ & & \spat \ \ \ \ \nu(q^{-m}y'x'/z')\,\sqrt{\frac{(-\sgn(p)q^{-2m}y,-x;q^2)_\infty}{(-z;q^2)_\infty}}
 \ \Psis{-q^2/x}{q^2 z/x}{\sgn(p)\,q^{2(1+m)}z/y}
\\ & & \spat = d\,c_q^2 \, \nu(p q^m)\,\nu(q^{-m-t+\chi(p)})\, s(xz,y)\,s(z,x)\, (-1)^{\chi(p)+\chi(z')}
\\ & & \spat \ \ \ \ \sgn(p)^{m+t+\chi(p)} \, \sgn(z')^{\chi(x')}\, \sgn(y')^{\chi(x'z')} \,
\, |x|^{\frac{1}{2}} \, |y|^{\frac{1}{2}}
\\ & & \spat \ \ \ \ \sqrt{\frac{(-\kappa(p),-x,-y;q^2)_\infty}{(-z;q^2)_\infty}}\,
\ \Psis{-q^2/\kappa(p)}{q^{2(1-m)}\,y/|\kappa(p)|}{\sgn(p)\,q^{2(1-m)}}
\\ & & \spat \ \ \ \ \Psis{-q^2/x}{\sgn(p)\,q^{2(1+t-\chi(p))} y}{q^{2(1+m+t-\chi(p))}x}
\\ & & \spat = \lambda^{\log_{q^2}\!|x|} \, \mu^{\log_{q^2}\!|y/x|}\,\eta^{-\log_{q^2}\!|y|}\,\,s(xz,y)\,s(z,x)\,(-1)^{\log_{q^2}\!|z|}
\,\sgn(z)^{\log_{q^2}\!|x|}\,\sgn(y)^{\log_{q^2}\!|xz|}\,|x|^{\frac{1}{2}}\,|y|^{\frac{1}{2}}
\\ & & \spat \ \ \ \ \sqrt{\frac{(-x,-y;q^2)_\infty}{(-z;q^2)_\infty}}\, (-1)^{\chi(p)}\, \sgn(p)^{m+t+\chi(p)}
 \ \lambda^{r+t+m-\chi(p)}\,\mu^{s-m}\,\eta^{n+m+\chi(p)}\,\,\tilde{F}^{t-\chi(p)}_{m,\kappa(p)}(x,z,y) \ .
\end{eqnarray*}
If $\sgn(p)\,q^{-2m}\,y \not\in I_{q^2}$ or $z \not= \sgn(p)\,q^{2(t-\chi(p))}\,x y$, then Eq.~(\ref{coass.eq4})
implies immediately that \newline $F^{r,s,t,n}_{m,p}(\lambda,\kappa^{-1}(x),\mu,\kappa^{-1}(z),\eta,\kappa^{-1}(y)) =
0$. So we see that
$$U F^{r,s,t,n}_{m,p} = (-1)^{\chi(p)}\, \sgn(p)^{m+t+\chi(p)}\,\,\ze^{r+t+m-\chi(p)}\ot \ze^{s-m}\ot \ze^{n+m+\chi(p)}\,\ot
\tilde{F}^{t-\chi(p)}_{m,\kappa(p)}\ .$$ Proposition \ref{coass.prop1} implies that
\begin{eqnarray*}
U G^{r,s,t,n}_{m,p} & = & U \Sigma_{13} (\breve{J} \ot \breve{J} \ot \breve{J})(V^* \ot 1_H)(1_{L^2(\T^2)} \ot
V^*)(\tilde{\Si} \ot \breve{J})
 (\ze^r \ot \ze^s \ot \ze^t \ot \ze^n \ot \ze^m \ot \sde_p)
\\ & = & \sgn(p)^{\chi(p)}\,U \Sigma_{13} (\breve{J} \ot \breve{J} \ot \breve{J})(V^* \ot 1_H)(1_{L^2(\T^2)} \ot V^*)
(\ze^{-n} \ot \ze^{-t} \ot \ze^{-s} \ot \ze^{-r} \ot \ze^{-m} \ot \sde_p)
\\ & = & \sgn(p)^{\chi(p)}\,U \Sigma_{13} (\breve{J} \ot \breve{J} \ot \breve{J})F^{-n,-t,-s,-r}_{-m,p} \ .
\end{eqnarray*}
Define the anti-unitary operator $\breve{I}$ on $L^2(\T^3) \ot L^2(I_{q^2}^3,\Om)$ such that $\breve{I}\,[f] = [g]$
where $f,g \in \cL^2(\T^3 \times I_{q^2}^3)$ are such that $g(\lambda,\mu,\eta,x,z,y) = \sgn(xyz)^{\log_{q^2}\!|xyz|}
\,\, \overline{ f(\eta,\mu,\lambda,y,z,x)}$ for all  $x,y,z \in I_{q^2}$ and $\lambda,\mu,\eta \in \T$.

\smallskip

Let $x,y,z \in I_{q^2}$. If $z > 0$, it follows from the considerations before Definition \ref{comult.def2} that
$s(xz,y) = \sgn(xy)\,s(yz,x)$ and $s(z,x) = \sgn(xy)\,s(z,y)$. If on the other hand $z < 0$, we have $s(xz,y) =
s(yz,x)$ and $s(z,x) = s(z,y)$. Thus, in both cases $s(xz,y) \,s(z,x) = s(yz,x) \,s(z,y)$. Now it is straightforward to
check that $ U \Sigma_{13} (\breve{J} \ot \breve{J} \ot \breve{J}) = \breve{I} U$. Thus,
\begin{eqnarray*}
U G^{r,s,t,n}_{m,p} & = & \sgn(p)^{\chi(p)}\,\breve{I} U F^{-n,-t,-s,-r}_{-m,p} \\ & = &
(-1)^{\chi(p)}\,\sgn(p)^{m+s}\, \breve{I}\,(\ze^{-n-s-m-\chi(p)} \ot \ze^{-t+m} \ot \ze^{-r-m+\chi(p)} \ot
\tilde{F}^{-s-\chi(p)}_{-m,\kappa(p)}) \ .
\end{eqnarray*}
Since $\tilde{F}^{-s-\chi(p)}_{-m,\kappa(p)}$ is supported on the set $L(q^{-2s}/\kappa(p))$ and
$\sgn(xyz)^{\log_{q^2}\!|xyz|} = \sgn(p)^{s+\chi(p)}$ for all $(x,y,z) \in L(q^{-2s}/\kappa(p))$, we get that
\begin{eqnarray*}
U G^{r,s,t,n}_{m,p} & = & (-1)^{\chi(p)}\,\sgn(p)^{m+\chi(p)}\,\ze^{r+m-\chi(p)} \ot \ze^{t-m} \ot \ze^{n+s+m+\chi(p)}
 \ot \tilde{G}^{-s-\chi(p)}_{m,\kappa(p)} \ .
\end{eqnarray*}
This implies for $r',s',t' \in \Z$, $k,m \in \Z$ and $p' \in I_{q^2}$,
\begin{eqnarray*}
\ze^{r'} \ot \ze^{s'} \ot \ze^{t'} \ot \tilde{F}_{m,p'}^k & = & (-1)^{\log_{q^2}\!|p'|}\,\sgn(p')^{m+k}\,U
F^{r'-m-k,s'+m,k+\log_{q^2}\!|p'|,t'-m-\log_{q^2}\!|p'|}_{m,\kappa^{-1}(p')}
\\ \ze^{r'} \ot \ze^{s'} \ot \ze^{t'} \ot \tilde{G}_{m,p'}^k & = & (-1)^{\log_{q^2}\!|p'|}\,\sgn(p')^{m+\log_{q^2}\!|p'|}\,U
F^{r'-m+\log_{q^2}\!|p'|,-k-\log_{q^2}\!|p'|,s'+m,t'-m+k}_{m,\kappa^{-1}(p')} \ .
\end{eqnarray*}
From these two equations and the properties listed in (\ref{coass.eq132}) in the beginning of the proof, all the claims
in the statement of this lemma easily follow.
\end{proof}

\medskip

Set $D = \langle \, \ze^r \mid r \in \Z \, \rangle \subseteq L^2(\T)$.

\begin{lemma} \label{coass.lem2} The operator $U\de_0^{(2)}(\gamma_0)U^*$ is an adjointable operator on $D^{\odot 3} \odot \cK(I_{q^2})^{\odot 3}$ such that \newline $(U\de_0^{(2)}(\gamma_0)U^*)^\dag = U\de_0^{(2)}(\gamma_0^\dag)U^*$.
Moreover
\begin{eqnarray*}
& & (U\de_0^{(2)}(\gamma_0)U^* \,g)(\lambda,\mu,\eta,x,z,y) = q\,|xy/z|^{\frac{1}{2}}\, \frac{\lambda \bar{\mu}
\eta}{xy}\,\bigl[\, (1+q^{-2} y) \,\,f(\lambda,\mu,\eta,q^2x,z,q^{-2}y)
\\ & & \hspace{8ex} - (1+q^{-2}y) \,\,g(\lambda,\mu,\eta,x,q^{-2}z,\eta,q^{-2}y) - (1+z) \,\,g(\lambda,\mu,\eta,q^2 x,q^2 z,y) +
g(\lambda,\mu,\eta,x,z,y) \, \bigr]
\end{eqnarray*}
for $g \in D^{\odot 3} \odot \cK(I_{q^2})^{\odot 3}$, $\lambda,\mu,\eta \in \T$ and $x,z,y \in I_{q^2}$.
\end{lemma}
\begin{proof}
Since the domain of $\de_0^{(2)}(\gamma_0)$ and $\de_0^{(2)}(\gamma_0^\dag)$ is $E^{\odot 3}$ and clearly
$U(E^{\odot 3}) = D^{\odot 3} \odot \cK(I_{q^2})^{\odot 3}$, the statements about the adjointability of
$\de_0^{(2)}(\gamma_0)$ follow easily. Using Eqs.~(\ref{hopf.comult2}), one checks that
\begin{eqnarray*}
& & \de_0^{(2)}(\gamma_0)  =  \gamma_0 \odot \al_0 \odot \al_0 + e_0 \al_0^\dag \odot \gamma_0 \odot \al_0  +  e_0
\al_0^\dag \odot e_0 \al_0^\dag \odot \gamma_0 + q\, \gamma_0 \odot e_0 \gamma_0^\dag \odot \gamma_0 \ .
\end{eqnarray*}
So we get for $f \in E^{\odot 3}$, $x,y,z \in I_q$ and $\lambda,\mu,\eta \in \T$ that
\begin{eqnarray*}
& &  (\de_0^{(2)}(\gamma_0) \,f)(\lambda,x,\mu,z,\eta,y)
\\ & & \spat   = \frac{\lambda}{x}\,\sqrt{\sgn(z)+q^2/z^2}\,\sqrt{\sgn(y)+q^2/y^2} \,\,f(\lambda,x,\mu,q^{-1}z,\eta,q^{-1}y)
\\ &  & \spat \ \ + \sgn(x) \, \frac{\mu}{z}\,\sqrt{\sgn(x)+1/x^2}\,\sqrt{\sgn(y)+q^2/y^2} \,\,f(\lambda,qx,\mu,z,\eta,q^{-1}y)
\\ &  & \spat \ \ + \sgn(xz) \frac{\eta}{y}\,\sqrt{\sgn(x)+1/x^2}\,\sqrt{\sgn(z)+1/z^2} \,\,f(\lambda,qx,\mu,qz,\eta,y)
  + \sgn(z)\ \frac{q \lambda \bar{\mu} \eta}{xyz} \, f(\lambda,x,\mu,z,\eta,y)
\\ & & \spat  =  \sgn(yz)\,\frac{q^2 \lambda}{xyz}\,\sqrt{1+q^{-2}\,\kappa(z)}\,\sqrt{1+q^{-2}\kappa(y)} \,\,f(\lambda,x,\mu,q^{-1}z,\eta,q^{-1}y)
\\ &  & \spat \ \ + \sgn(y)
\, \frac{q \mu}{xyz}\,\sqrt{1+\kappa(x)}\,\sqrt{1+q^{-2} \kappa(y)} \,\,f(\lambda,qx,\mu,z,\eta,q^{-1}y)
\\ &  & \spat \ \ + \frac{\eta}{xyz}\,\sqrt{1+\kappa(x)}\,\sqrt{1+\kappa(z)} \,\,f(\lambda,qx,\mu,qz,\eta,y)
+ \sgn(z)\ \frac{q\lambda \bar{\mu} \eta}{xyz} \, f(\lambda,x,\mu,z,\eta,y) \ .
\end{eqnarray*}
The inverse $U^*$ of $U$ is such that
\begin{eqnarray*}
(U^* f)(\lambda,x,\mu,z,\eta,y) & = &
s(xz,y)\,s(z,x)\,\lambda^{\log_q\!|x|}\,\mu^{\log_q\!|y/x|}\,\eta^{-\log_q\!|y|}\,(-1)^{\log_q\!|z|}\,\sgn(z)^{\log_q\!|x|}
\, \sgn(y)^{\log_q\!|xz|}\,
\\ &  & |x| \,  |y|\, \sqrt{\frac{(-\kappa(x),-\kappa(y);q^2)_\infty}{(-\kappa(z);q^2)_\infty}}\,\,f(\lambda,\mu,\eta,\kappa(x),\kappa(z),\kappa(y))
\end{eqnarray*}
for $f \in D^{\odot 3}  \odot \cK(I_q)^{\odot 3}$, $\lambda,\mu,\eta \in \T$ and $x,z,y \in I_q$.

Then the following rules are easily established for $g \in D^{\odot 3}  \odot \cK(I_{q^2})^{\odot 3}$, $x,y,z \in
I_{q^2}$ and $\om = (\lambda,\eta,\mu) \in \T^3$.
\begin{trivlist}
\item[\ \,(1)] If $f \in D^{\odot 3}  \odot \cK(I_q)^{\odot 3}$, then $(U M_f U^* g)(\om,x,z,y) =
f(\lambda,\kappa^{-1}(x),\mu,\kappa^{-1}(z),\eta,\kappa^{-1}(y))$ $g(\om,x,z,y)$. \vspace{1ex}

\item[\ \,(2)] $(U  (1  \odot  T_{q^{-1}} \odot
 T_{q^{-1}}) U^* g)(\om,x,z,y) = - q^{-1}\,\sgn(y)\,\bar{\mu} \eta \,
\sqrt{(1+q^{-2}y)(1+q^{-2}z)^{-1}}\, g(\om,x,q^{-2} z,q^{-2} y)$. \vspace{1ex}

\item[\ \,(3)] $(U  (T_q \odot 1  \odot
 T_{q^{-1}}) U^* g)(\om,x,z,y) =  \sgn(yz)\,\lambda \bar{\mu}^2 \eta \,
\sqrt{(1+q^{-2}y)(1+x)^{-1}}\, g(\om,q^2x,z,q^{-2} y)$. \vspace{1ex}

\item[\ \,(4)] $(U  (T_q \odot  T_q \odot 1) U^* g)(\om,x,z,y) =  - q\,\sgn(z)\,\lambda \bar{\mu} \, \sqrt{(1+z)(1+x)^{-1}}\,
g(\om,q^2x,q^2z,y)$.
\end{trivlist}
This implies that for $g \in D^{\odot 3}  \odot \cK(I_{q^2})^{\odot 3}$, $x,y,z \in I_{q^2}$ and $\om =
(\lambda,\eta,\mu) \in \T^3$,
\begin{eqnarray*}
& & (U\de_0^{(2)}(\gamma_0)U^* \,g)(\om,x,z,y)
\\ & & \spat =  - \sgn(z)\,\frac{q \lambda\bar{\mu}\eta}{\kappa^{-1}(xyz)}\,(1+q^{-2}y) \,\,g(\om,x,q^{-2}z,q^{-2}y)
\\ & & \spat\ \ \  + \sgn(z)\,\frac{q \lambda\bar{\mu}\eta}{\kappa^{-1}(xyz)}\, \, (1+q^{-2} y) \,\,g(\om,q^2x,z,q^{-2}y)
\\ &  & \spat\ \ \ - \sgn(z)\,\frac{q \lambda\bar{\mu}\eta}{\kappa^{-1}(xyz)}\,(1+z) \,\,g(\om,q^2 x,q^2 z,y)
\\ & & \spat\ \ \  + \sgn(z)\ \frac{q\lambda \bar{\mu} \eta}{\kappa^{-1}(xyz)} \, g(\om,x,z,y) \ .
\\ & & \spat = q\,|xy/z|^{\frac{1}{2}}\, \frac{\lambda \bar{\mu} \eta}{xy}\,\bigl(\,(1+q^{-2} y) \,\,g(\om,q^2x,z,q^{-2}y)
- (1+q^{-2}y) \,\,g(\om,x,q^{-2}z,q^{-2}y)
\\ & & \spat\ \ - (1+z) \,\,g(\om,q^2 x,q^2 z,y) + g(\om,x,z,y) \, \bigr) \ .
\end{eqnarray*}
\end{proof}

Let us also get rid of part of the operators  acting on $L^2(\T^3)$. For this purpose we introduce normal operators
$\gamma_l$ on $\gamma_r$ in $L^2(I_q^3,\Om)$ such that
\begin{trivlist}
\item[\ \,(1)] the space $\langle \, \tilde{F}^k_{m,p} \mid p \in I_{q^2}, m,k \in \Z \rangle$ is a core for $\gamma_l$,
\smallskip
\item[\ \,(2)] the space $\langle \, \tilde{G}^k_{m,p} \mid p \in I_{q^2}, m,k \in \Z \rangle$ is a core for $\gamma_r$,
\smallskip
\item[\ \,(3)] if $p \in I_{q^2}$ and $m,k \in \Z$, then $\gamma_l \,\tilde{F}^k_{m,p}  =  |p|^{-\frac{1}{2}}\, \tilde{F}^k_{m+1,p}$ and $\gamma_r \,\tilde{G}^k_{m,p}  =  |p|^{-\frac{1}{2}}\,  \tilde{G}^k_{m+1,p}$.
\end{trivlist}

\medskip

By Lemma \ref{coass.lem1} we know that $U((\de \ot \io)\de(\gamma))U^* = M_\vsi \ot \gamma_l$ and $U((\io \ot
\de)\de(\gamma))U^* = M_\vsi \ot \gamma_r$. So in order to prove that $(\de \ot \io)\de(\gamma) = (\io \ot
\de)\de(\gamma)$, it is enough to show that $\gamma_l = \gamma_r$.

\smallskip

The description of $\gamma_l$ and $\gamma_r$ in terms of eigenvectors will not be sufficient to solve this problem. We
will also need to know the explicit action of $\gamma_l$ and $\gamma_r$ on functions in its domain. For this purpose we
introduce the linear operators $\tilde{\gamma}$, $\tilde{\gamma}^\dag$ belonging to  $\End(\cF(I_{q^2}^3))$ such that
\begin{eqnarray*}
(\tilde{\gamma}\,f)(x,z,y) &  = & q\,|xy/z|^{\frac{1}{2}}\, \frac{1}{xy}\,\bigl(\, (1+q^{-2} y) \,\,f(q^2x,z,q^{-2}y) -
(1+q^{-2}y) \,\,f(x,q^{-2}z,q^{-2}y)
\\ & & \hspace{22ex} - (1+z) \,\,f(q^2 x,q^2 z,y) + f(x,z,y) \, \bigr)
\end{eqnarray*}
and
\begin{eqnarray*}
(\tilde{\gamma}^\dag\,f)(x,z,y) &  = & q\,|xy/z|^{\frac{1}{2}}\, \frac{1}{xy}\,\bigl(\, (1+q^{-2} x)
\,\,f(q^{-2}x,z,q^2y) -  (1+z) \,\,f(x,q^2z,q^2y)
\\ & & \hspace{22ex} - (1+q^{-2} x) \,\,f(q^{-2} x,q^{-2} z,y) + f(x,z,y) \, \bigr)
\end{eqnarray*}
for all $f \in \cF(I_{q^2}^3)$ and $x,y,z \in I_{q^2}$. Note that $\tilde{\gamma}^\dag$ is obtained from
$\tilde{\gamma}$ by interchanging $x$ and $y$.

\medskip

\medskip

\begin{lemma} \label{coass.lem3}
If $f \in D(\gamma_l)$ and $g \in D(\gamma_r^*)$, then $\gamma_l(f) = \tilde{\gamma}(f)$ and $\gamma_r^*(g) =
\tilde{\gamma}^\dag(g)$.
\end{lemma}
\begin{proof}
Define the sesquilinear form $\langle .\,,\,.\rangle : \cK(I_{q^2}^3) \times \cF(I_{q^2}^3) \rightarrow \C$ such that
$$\langle f , g \rangle = \sum_{(x,z,y) \in I_{q^2}^3}\,f(x,z,y)\,\bar{g}(x,z,y)\, |x|\,|y|\,\frac{(-x,-y;q^2)_\infty}{(-z;q^2)_\infty}$$
for all $f \in \cK(I_{q^2}^3)$ and $g \in \cF(I_{q^2}^3)$.

\smallskip

Under this pairing, $\cF(I_{q^2}^3)$ is anti-linearly isomorphic to the algebraic dual of $\cK(I_{q^2}^3)$. So there
exist an anti-linear anti-homomorphism $^\circ : \End(\cK(I_{q^2}^3)) \rightarrow \End(\cF(I_{q^2}^3))$ such that
$\langle T(f) , g \rangle = \langle f , T^\circ(g) \rangle$ for all $f \in \cK(I_{q^2}^3)$ and $g \in \cF(I_{q^2}^3)$.
If $f \in \cK(I_{q^2}^3)$ and $g \in L^2(I_{q^2}^3,\Om)$ then $\langle f , g \rangle$ equals the inner product in
$L^2(I_{q^2}^3,\Om)$. So if $T \in \End(\cK(I_{q^2}^3))$, then $T^*(f) = T^\circ(f)$ for all $f \in D(T^*)$\   (and the
domain of $T^*$ consists of all $f \in L^2(I_{q^2}^3,\Om)$ such that $T^\circ(f) \in L^2(I_{q^2}^3,\Om)$\,). Moreover,
$T \in \End(\cK(I_{q^2}^3))$ is adjointable if and only if $T^\circ(\cK(I_{q^2}^3)) \subseteq \cK(I_{q^2}^3)$ in which
case $T^\dag$ is the restriction of $T^\circ$ to $\cK(I_{q^2}^3)$.

\smallskip

If $p = (p_1,p_2,p_3) \in (-q^{2\Z} \cup q^{2\Z})^3$ and $g \in \cF(I_{q^2}^3)$, we  define $\tilde{T}_p, \tilde{M}_g
\in \text{End}(\cF(I_{q^2}^3))$ such that
$$\tilde{T}_p(f)(x,z,y) = f(p_1 x, p_2 z, p_3 y) \hspace{10ex} \text{and} \hspace{10ex}
\tilde{M}_g(f)(x,z,y) = g(x,z,y) \, f(x,z,y)$$ for all $f \in \cF(I_{q^2}^3)$ and $x,z,y \in I_{q^2}$. We denote the
restrictions of $\tilde{M}_g, \tilde{T}_p$ to $\cK(I_{q^2}^3)$ by $\hat{M}_g, \hat{T}_p \in \End(\cK(I_{q^2}^3))$
respectively.

Let $\xi_1,\xi_2,\xi_3$ denote the 3 coordinate functions on $I_{q^2}^3$. Define also an auxiliary function $\om :
I_{q^2}^3 \rightarrow \R$ such that for $x,y,z \in I_{q^2}$, $\om(x,z,y) = (1+q^{-2} z)^{-1}$ if $z \not= -q^2$ and
$\om(x,z,y) = 0$ if $z = - q^2$. It is not so difficult to check that
$$\begin{array}{lcl}
\hat{T}_{q^2,1,1}^\circ = q^{-2}\,\tilde{M}_{1+q^{-2} \xi_1}\,\tilde{T}_{q^{-2},1,1} & \hspace{8ex} \vspace{1ex} &
\hat{T}_{1,1,q^{-2}}^\circ = q^2\,\tilde{M}_{(1+ \xi_3)^{-1}}\,\tilde{T}_{1,1,q^2}
\\  \hat{T}_{1,q^2,1}^\circ = \,\tilde{M}_\om\,\tilde{T}_{1,q^{-2},1} & \vspace{1ex} &
\hat{T}_{1,q^{-2},1}^\circ = \tilde{M}_{1+ \xi_2}\,\tilde{T}_{1,q^2,1}
\\ \hat{M}_g^\circ = \tilde{M}_{\bar{g}} \ \text{ if } g \in \cF(I_{q^2}^3)  \ . & &
\end{array}$$
Thus
$$\begin{array}{lcl} \hat{T}_{q^2,q^2,1}^\circ = q^{-2}\,\tilde{M}_{\om(1+q^{-2} \xi_1)}\,\tilde{T}_{q^{-2},q^{-2},1}
& \hspace{8ex} \vspace{1ex} & \hat{T}_{1,q^{-2},q^{-2}}^\circ = q^2\,\tilde{M}_{(1+ \xi_2)(1+
\xi_3)^{-1}}\,\tilde{T}_{1,q^2,q^2}
\\ \hat{T}_{q^2,1,q^{-2}}^\circ = \tilde{M}_{(1+q^{-2} \xi_1)(1+ \xi_3)^{-1}} \, \tilde{T}_{q^{-2},1,q^2} \ . & &
\end{array} $$

Let $\hat{\gamma}$ denote the restriction of $\tilde{\gamma}$ to $\cK(I_{q^2}^3)$. Clearly, $\hat{\gamma} =
\sum_{p}\,\, \hat{T}_p\, \hat{M}_{g_p}$ where $p$ ranges over $(q^2,q^2,1)$, $(1,q^{-2},q^{-2})$, $(q^2,1,q^{-2})$,
$(1,1,1)$ and $g_p$ belongs to $\cF(I_q^{3})$. Therefore $\hat{\gamma}$ is adjointable.

Moreover, $\hat{\gamma}^\circ = \sum_{p}\,\, \hat{M}_{g_p}\, \hat{T}_p^\circ$. Now an easy calculation reveals that
$\hat{\gamma}^\circ = \tilde{\gamma}^\dag$. Similarly,  $(\hat{\gamma}^\dag)^\circ = \tilde{\gamma}$.

\smallskip

Lemma \ref{coass.lem2} implies that $U\de_0^{(2)}(\gamma_0)U^* = M_\vsi\!\!\restriction_{D^{\odot 3}} \odot
\hat{\gamma}$ and $U\de_0^{(2)}(\gamma^\dag)U^* = M_\vsi^*\!\!\restriction_{D^{\odot 3}} \odot \hat{\gamma}^\dag$.

By Proposition \ref{coass.prop2} and the remarks before this lemma we get that
$$M_\vsi^*\!\!\restriction_{D^{\odot 3}} \odot
\hat{\gamma}^\dag = U\de_0^{(2)}(\gamma_0^\dag)U^* \subseteq U((\de \ot \io)\de(\gamma^*))U^* = M_\vsi^* \ot \gamma_l^*
\ .$$ Multiplying this inclusion from the right with $M_\vsi \ot 1_{L^2(I_{q^2}^3,\Om)}$, we get that $1_{D^{\odot 3}}
\odot \hat{\gamma}^\dag \subseteq 1_{L^2(\T^3)} \ot \gamma_l^*$.

Fix a unit vector $w \in D^{\odot 3}$. Define the bounded linear operator $T : L^2(\T^3) \rightarrow \C$ such that
$T(x) = \langle x , w \rangle$ for all $x \in L^2(\T^3)$. Then $(T \ot 1_{L^2(I_{q^2}^3,\Om)})\,(1_{L^2(\T^3)} \ot
\gamma_l^*) \subseteq \gamma_l^*\,(T \ot 1_{L^2(I_{q^2}^3,\Om)})$. So if $v \in \cK(I_{q^2}^3)$, then $w \ot v$ $\in$
$D(1_{L^2(\T^3)} \ot \gamma_l^*)$ and $(1_{L^2(\T^3)} \ot \gamma_l^*)(w \ot v) = w \ot \hat{\gamma}^\dag\,v$. Thus, $v
= (T \ot 1_{L^2(I_{q^2}^3,\Om)})(w \ot v) \in D(\gamma_l^*)$ and
$$\gamma_l^* v  =  \gamma_l^* (T \ot 1_{L^2(I_{q^2}^3,\Om)})(w \ot v) = (T \ot 1_{L^2(I_{q^2}^3,\Om)})(1_{L^2(\T^3)} \ot \gamma_l^*)(w \ot v)
 =  (T \ot 1_{L^2(I_{q^2}^3,\Om)})(w \ot \hat{\gamma}^\dag v) = \hat{\gamma}^\dag v\ .$$
So we have shown that $\hat{\gamma}^\dag \subseteq \gamma_l^*$. Taking the adjoint of this inclusion we see that
$\gamma_l \subseteq (\hat{\gamma}^\dag)^*$. By the remarks in the beginning of the proof, this implies for $f \in
D(\gamma_l)$, $\gamma_l(f) = (\hat{\gamma}^\dag)^*(f) = (\hat{\gamma}^\dag)^\circ(f) = \tilde{\gamma}(f)$. Starting
from the inclusion $\de_0^{(2)}(\gamma_0) \subseteq  (\io \ot \de)\de(\gamma)$, the statement concerning $\gamma_r^*$
is proven in the same way.
\end{proof}

\medskip

Let us further reduce the 3-dimensional problem  to a 2-dimensional one. If $\theta \in q^{2\Z} \cup - q^{-2\Z}$, we
define $L^2_\theta(I_{q^2}^3,\Om)$ to be the closed subspace of $L^2_\theta(I_{q^2}^3,\Om)$ consisting of all functions
that have support in $L(\theta)$. Then $L^2(I_{q^2}^3,\Om)$ is the orthogonal sum $\oplus_{\theta \in q^{2\Z} \cup -
q^{-2\Z}} L^2_\theta(I_{q^2}^3,\Om)$. We know that  $\tilde{F}^k_{m,p}, \tilde{G}^k_{m,p} \in
L^2_{\sgn(p)\,q^{2k}}(I_{q^2}^3,\Om)$ for all $p \in I_{q^2}$, $m,k \in \Z$.

Thus, lemma \ref{coass.lem1} implies that both families $(\,\tilde{F}^{\log_{q^2}\!|\theta|}_{m,p} \mid p \in I_{q^2}
\text{ such that } \sgn(p) = \sgn(\theta), m \in \Z\,)$ and $(\,\tilde{G}^{\log_{q^2}\!|\theta|}_{m,p} \mid p \in
I_{q^2} \text{ such that } \sgn(p) = \sgn(\theta), m \in \Z\,)$ are orthonormal bases of $L^2_\theta(I_{q^2}^3,\Om)$.

\medskip

It follows from the definition of $\gamma_l$ and $\gamma_r$ that $L^2_\theta(I_{q^2}^3,\Om)$ is invariant under
$\gamma_l$ and $\gamma_r$.

\medskip

Before starting the proof of the next result, let us first introduce some notation. Therefore suppose that $C \subseteq
(-q^{2\Z} \cup q^{2\Z})\times (-q^{2\Z} \cup q^{2\Z})$ and consider a function $f : C \rightarrow \C$. Then we define
new functions $\partial_1 f, \partial_2 f : C \rightarrow \C$ such that $(\partial_1 f)(x,y) = \frac{f(q^2 x,y) -
f(x,y)}{x}$ and $(\partial_2 f)(x,y) = \frac{f(x,q^2 y) - f(x,y)}{y}$ for all $(x,y) \in C$ (recall discussion (2) in
Notation and conventions).

\medskip

\begin{proposition} \label{coass.prop3}
We have that $\gamma_l = \gamma_r$.
\end{proposition}
\begin{proof}
Fix $\theta \in -q^{2\Z} \cup q^{2\Z}$. Set $k = \log_{q^2}\!|\theta|$. On $K(\theta)$ we define the measure
$\Om_\theta$ such that
$$\Omega_\theta(\{(x,y)\}) = |x|\,|y|\,\frac{(-x,-y;q^2)_\infty}{(-\theta xy;q^2)_\infty}$$ for all $(x,y)
\in K(\theta)$. Define a unitary transformation $\Xi : L_\theta(I_{q^2}^3,\Om) \rightarrow L^2(K(\theta),\Om_\theta)$
such that $(\Xi f)(x,y) = f(x,\theta xy,y)$ for all $f \in L_\theta(I_{q^2}^3,\Om)$ and $(x,y) \in K(\theta)$. Notice
that $\Xi \tilde{F}^k_{m,p} = F^k_{m,p}$ and $\Xi \tilde{G}^k_{m,p} = G^k_{m,p}$ for all $m \in \Z$ and $p \in I_{q^2}$
such that $\sgn(p) = \sgn(\theta)$. Set $\gamma_r^\theta = \Xi\, \gamma_r\!\!\restriction_{L_\theta(I_{q^2}^3,\Om)}
\Xi^*$ and $\gamma_l^\theta = \Xi\, \gamma_l\!\!\restriction_{L_\theta(I_{q^2}^3,\Om)} \Xi^*$. Then,
\begin{trivlist}
\item[\ \,(1)] the space $\langle\,F^k_{m,p} \mid p \in I_{q^2} \text{ such that } \sgn(p) = \sgn(\theta), m \in
\Z\,\rangle$ is a core of $\gamma_l^\theta$,
\smallskip
\item[\ \,(2)] the space $\langle\,G^k_{m,p} \mid p \in I_{q^2} \text{ such
that } \sgn(p) = \sgn(\theta), m \in \Z\,\rangle$ is a core of $(\gamma_r^\theta)^*$. \ \inlabel{coass.eq15}
\end{trivlist}

Choose $f \in D(\gamma^\theta_l)$ and $g \in D((\gamma^\theta_r)^*)$. By Lemma \ref{coass.lem3}, we know that for
$(x,y) \in K(\theta)$,
\begin{equation} \label{coass.eq18}
(\gamma_l^\theta\,f)(x,y)   =  q\,|\theta|^{-\frac{1}{2}}\, \frac{1}{xy}\,\bigl(\, (1+q^{-2} y) \,\,f(q^2x,q^{-2}y) -
(1+q^{-2}y) \,\,f(x,q^{-2}y)  - (1+\theta xy) \,\,f(q^2 x,y) + f(x,y) \, \bigr)
\end{equation}
and
\begin{equation} \label{coass.eq19}
((\gamma_r^\theta)^*\,g)(x,y) =  q\,|\theta|^{-\frac{1}{2}}\, \frac{1}{xy}\,\bigl(\,(1+q^{-2} x) \,\,g(q^{-2}x,q^2y) -
(1+\theta x y) \,\,g(x,q^2y) - (1+q^{-2} x) \,\,g(q^{-2} x,y) + g(x,y) \, \bigr) \ .
\end{equation}

Note the symmetry between the above formulas for $\gamma^\theta_l$ and $(\gamma^\theta_r)^*$ with respect to
interchanging $x$ and $y$.

\smallskip

Define the auxiliary function $h : K(\theta) \rightarrow \C$ such that $h(x,y) = \frac{(-x,-y;q^2)_\infty}{(-\theta x
y;q^2)_\infty}$ for all $(x,y) \in K(\theta)$.

\smallskip

We want to show that $\langle \gamma^\theta_l f , g \rangle = \langle f ,  (\gamma^\theta_r)^* g \rangle$. The
difference of these inner products are sums over the whole area $K(\theta)$ but we will start by approximating these
sums by sums over a finite subset of $K(\theta)$ obtained by cutting off $K(\theta)$ by 4 rectangles (see Fig.~1 and
Fig.~2 later in the proof). Then we let these rectangles get bigger and bigger. In this way, these finite sums converge
certainly to the sum over the whole of $K(\theta)$. In the second part of the proof, we will then prove that these
finite sums converge to 0, implying that $\langle \gamma^\theta_l f , g \rangle$ must be equal to $\langle f ,
(\gamma^\theta_r)^* g \rangle$.

\smallskip

So let us first calculate the contribution of one rectangle, lying in one of the quadrants, to the difference of the
above inner product. Therefore let  $a,b,c,d \in I_{q^2}$ such that $\sgn(a) = \sgn(b)$, $\sgn(c) = \sgn(d)$, $|a| <
|b|$ and $|c| < |d|$ and set
\begin{equation}
C(a,b;c,d)
 =  \sum_{\scriptsize \begin{array}{c} (x,y) \in K(\theta) \\ x \text{ between } a \text{ and } b \\
y \text{ between } c \text{ and } d \end{array}} \bigr(\,(\gamma_l^\theta f)(x,y)\, \bar{g}(x,y)\,- f(x,y)\,
\overline{((\gamma_r^\theta)^* g)(x,y)}\,\bigl)\,  |x|\,|y|\,h(x,y) \ . \label{coass.eq17}
\end{equation}
where the statement \lq$x \text{ between } a \text{ and } b$\rq\  allows for $x$ to be equal to $a$ or $b$ as well
(similarly for $y$).

Using Eqs.~(\ref{coass.eq18}) and (\ref{coass.eq19}) and recalling discussion (2) in Notations and conventions, we get
that
\begin{eqnarray}
& & \hspace{-6ex} q^{-1}\,|\theta|^{\frac{1}{2}}\,\sgn(ac)\,C(a,b;c,d) \nonumber
\\ & & \hspace{-6ex}  \spat =  \sum_{x=a}^b \sum_{y=c}^d  (1+q^{-2} y) \,\,f(q^2x,q^{-2}y)\,\bar{g}(x,y)\,h(x,y)
- \sum_{x=a}^b \sum_{y=c}^d (1+q^{-2} x) \,\,f(x,y)\,\bar{g}(q^{-2}x,q^2y)\,h(x,y) \nonumber
\\ & & \hspace{-6ex}  \spat\ \ \ - \sum_{x=a}^b \sum_{y=c}^d (1+q^{-2}y) \,\,f(x,q^{-2}y)\,\bar{g}(x,y)\,h(x,y)
+ \sum_{x=a}^b \sum_{y=c}^d (1+\theta x y) \,\,f(x,y) \bar{g}(x,q^2y)\,h(x,y) \nonumber
\\ & & \hspace{-6ex}  \spat\ \ \ - \sum_{x=a}^b \sum_{y=c}^d (1+\theta xy) \,\,f(q^2 x,y)\,\bar{g}(x,y)\,h(x,y)
+ \sum_{x=a}^b \sum_{y=c}^d (1+q^{-2} x) \,\,f(x,y)\,\bar{g}(q^{-2} x,y)\,h(x,y) \nonumber
\\ & & \hspace{-6ex}  \spat
= \sum_{x=a}^b \sum_{y=q^{-2} c}^{q^{-2} d} (1+ y) \,\,f(q^2x,y)\,\bar{g}(x,q^2 y)\,h(x,q^2 y) -
\sum_{x=q^{-2}a}^{q^{-2}b} \sum_{y=c}^d (1+ x) \,\,f(q^2 x,y)\,\bar{g}(x,q^2y)\,h(q^2 x,y) \nonumber
\\ & & \hspace{-6ex} \spat\ \ \ - \sum_{x=a}^b \sum_{y=q^{-2}c}^{q^{-2}d} (1+y) \,\,f(x,y)\,\bar{g}(x,q^2 y)\,h(x,q^2 y)
+ \sum_{x=a}^b \sum_{y=c}^{d} (1+\theta x y) \,\,f(x,y) \bar{g}(x,q^2y)\,h(x,y) \nonumber
\\ & & \hspace{-6ex} \spat\ \ \  - \sum_{x=a}^{b} \sum_{y=c}^d (1+\theta xy) \,\,f(q^2 x,y)\,\bar{g}(x,y)\,h(x,y)
+ \sum_{x=q^{-2}a}^{q^{-2}b} \sum_{y=c}^d (1+ x) \,\,f(q^2x,y)\,\bar{g}(x,y)\,h(q^2 x,y) \label{coass.eq5} \ .
\end{eqnarray}
Define the set $S = \{\,(x,y) \in (-q^{2\Z} \cup q^{2\Z})\times  (-q^{2\Z} \cup q^{2\Z}) \mid q^2 \theta x y \in
I_{q^2} \, \}$. Define a new auxiliary function $h' : S \rightarrow \C$ such that $h'(x,y) =
\frac{(-x,-y;q)_\infty}{(-q^2 \theta xy;q^2)_\infty}$ for all $(x,y) \in S$. Now observe the following basic facts:
$$
\begin{array}{lll}
(1+y) h(x,q^2y) = h'(x,y)  \ \ \text{ if } (x,q^2 y) \in K(\theta) & \hspace{0.8cm} &  (1+x) h(q^2 x,y) = h'(x,y) \ \
\text{ if } (q^2 x, y) \in K(\theta)
\\ \vspace{-2mm}
\\  (1+\theta xy) h(x,y) = h'(x,y)\ \  \text{ if } (x,y) \in K(\theta) \ .
\end{array}
$$
Remembering that $f$ and $g$ are defined on $K(\theta)$, these facts combined with Eq.~(\ref{coass.eq5}) imply that
\begin{eqnarray*}
& & q^{-1}\,|\theta|^{\frac{1}{2}}\,\sgn(ac)\,C(a,b;c,d) \nonumber
\\ & & \spat = \sum_{x=a}^b \sum_{y=q^{-2} c}^{q^{-2} d} f(q^2x,y)\,\bar{g}(x,q^2 y)\,h'(x,y)
- \sum_{x=q^{-2}a}^{q^{-2}b} \sum_{y=c}^d f(q^2 x,y)\,\bar{g}(x,q^2y)\,h'(x,y) \nonumber
\\ & & \spat\ \ \ - \sum_{x=a}^b \sum_{y=q^{-2}c}^{q^{-2}d} f(x,y)\,\bar{g}(x,q^2 y)\,h'(x,y)
+ \sum_{x=a}^b \sum_{y=c}^{d} f(x,y) \bar{g}(x,q^2y)\,h'(x,y) \nonumber
\\ & & \spat\ \ \ - \sum_{x=a}^{b} \sum_{y=c}^d f(q^2 x,y)\,\bar{g}(x,y)\,h'(x,y)
+ \sum_{x=q^{-2}a}^{q^{-2}b} \sum_{y=c}^d f(q^2x,y)\,\bar{g}(x,y)\,h'(x,y) \ .\nonumber
\end{eqnarray*}
In the above sum, most of the terms of the sums in the left column are cancelled by the terms of the sums in the right
column. What remains is
\begin{eqnarray*}
& & q^{-1}\,|\theta|^{\frac{1}{2}}\,\sgn(ac)\,C(a,b;c,d) \nonumber
\\ & & \spat = \sum_{x=a}^b f(q^2 x,q^{-2} d)\,\bar{g}(x,d)\,h'(x,q^{-2} d)
+ \sum_{y=q^{-2}c}^{d} f(q^2 a,y)\,\bar{g}(a,q^2 y)\,h'(a,y) \nonumber
\\ & & \spat\ \ \ - \sum_{x=q^{-2}a}^{b} f(q^2 x,c)\,\bar{g}(x,q^2 c)\,h'(x,c)
- \sum_{y=c}^d f(b,y)\,\bar{g}(q^{-2}b,q^2y)\,h'(q^{-2}b,y) \nonumber
\\ & & \spat\ \ \ - \sum_{x=a}^b f(x,q^{-2} d)\,\bar{g}(x,d)\,h'(x,q^{-2} d)
+ \sum_{x=a}^b f(x,c) \bar{g}(x,q^2 c)\,h'(x,c) \nonumber
\\ & & \spat\ \ \ - \sum_{y=c}^d f(q^2 a,y)\,\bar{g}(a,y)\,h'(a,y)
 + \sum_{y=c}^d f(b,y)\,\bar{g}(q^{-2}b,y)\,h'(q^{-2}b,y) \nonumber
\end{eqnarray*}

\begin{eqnarray}
& & \spat =  \sum_{x=a}^b (f(q^2 x,q^{-2} d)-f(x,q^{-2} d)\,) \,\bar{g}(x,d)\,h'(x,q^{-2} d) \nonumber
\\ & & \spat\ \ \ - \sum_{y=c}^d  f(b,y)\,(\bar{g}(q^{-2}b,q^2y)-\bar{g}(q^{-2}b,y)\,)\,h'(q^{-2}b,y) \nonumber
\\ & & \spat\ \ \ + f(a,c)\,\bar{g}(a,q^2 c) h'(a,c) - \sum_{x=q^{-2} a}^b (f(q^2 x,c) - f(x,c)\,)\,\bar{g}(x,q^2
c)\,h'(x,c)\nonumber
\\ & & \spat\ \ \ - f(q^2a,c)\,\bar{g}(a, c) h'(a,c) + \sum_{y=q^{-2} c}^d f(q^2 a,y)\,(\bar{g}(a,q^2 y)-
\bar{g}(a,y)\,)\,h'(a,y) \ .\label{coass.eq6}
\end{eqnarray}

\medskip

Set $I = \{\,v \in I_{q^2}^+ \mid v < q^2, v < |\theta|^{-1} \text{ and } v < |\theta|^{-\frac{1}{2}} q\,\}$. If $v \in
I$ we define $v' = q^2/v |\theta|$. Let us calculate the contribution of the 4 rectangles depicted in Fig.~1 and Fig.~2
to $\langle \gamma^\theta_l f , g \rangle - \langle f ,  (\gamma^\theta_r)^* g \rangle$. Therefore set $C(v) =
C(v,v';v,v') + C(-v,-q^2;v,v') + C(-v,-q^2;-v,-q^2) + C(v,v';-v,-q^2)$.

It is clear from Eq.~(\ref{coass.eq17}) and the dominated convergence theorem that $C(v)$ converges to $\langle
\gamma_l^\theta f , g \rangle - \langle f , (\gamma^\theta_r)^* g \rangle$ as $v \rightarrow 0$. \inlabel{coass.eq14}

\medskip

\begin{picture}(210,220)(-65,20)
\includegraphics[bb=140 600 595 820]{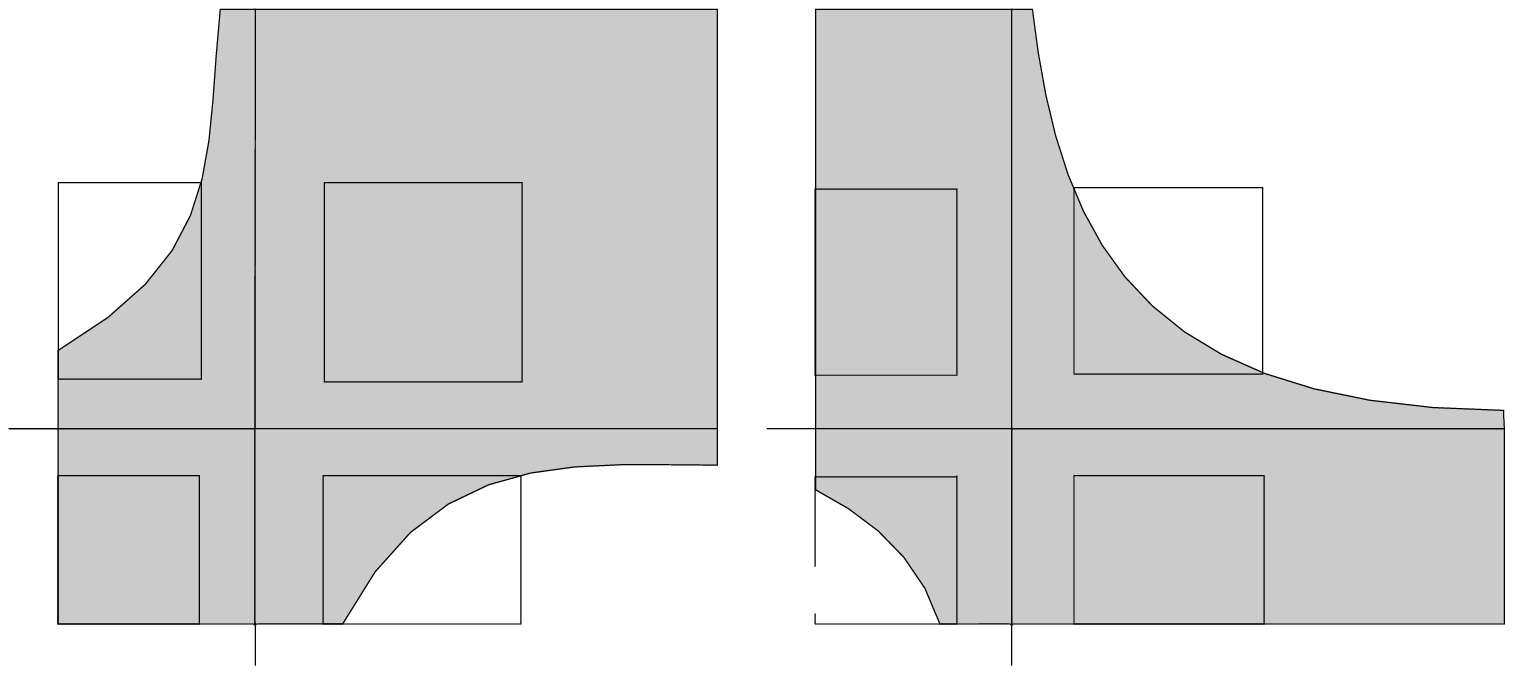}
\put(-345,190){$K(\theta)$} \put(-274,190){$K(\theta)$}

\put(-512,180){$\theta x y = - q^2$} \put(-205,180){$\theta x y = - q^2$}

\put(-360,63){$\theta x y = - q^2$} \put(-305,39){$\theta x y = - q^2$}

\put(-502,76){$-q^2$} \put(-283,76){$-q^2$}

\put(-445,34){$-q^2$} \put(-227,34){$-q^2$}

\put(-440,-7){Fig.~1 : $\theta > 0$} \put(-227,-7){Fig.~2 : $\theta < 0$}

\put(-424,106){$(v,v)$} \put(-399,148){$(v',v')$}

\put(-207,108){$(v,v)$} \put(-184,146){$(v',v')$}
\end{picture}

\vspace{2cm}

If $\theta < 0$, we set
$$
C_{\infty,1}(v)  =  \sum_{x=-v}^{-q^2} (\partial_1 f)(x,q^{-2} v') \,\bar{g}(x,v')\,h'(x,q^{-2} v')\,|x| -
\sum_{y=-v}^{-q^2} f(v',y)\,(\partial_2 \bar{g})(q^{-2} v',y)\,h'(q^{-2}v',y)\,|y|
$$
and if $\theta > 0$, we set
$$
C_{\infty,1}(v) = \sum_{x=v}^{v'} (\partial_1 f)(x,q^{-2} v')\,\bar{g}(x,v')\,h'(x,q^{-2} v')\,|x|
 - \sum_{y=v}^{v'} f(v',y)\,(\partial_2 \bar{g})(q^{-2}v',y)\,h'(q^{-2} v',y)\,|y| \ .
$$
Regardless of the sign of $\theta$, we set

\begin{eqnarray*}
C_{\infty,0}(v) & =  & \sgn(\theta)\,f(v',-\sgn(\theta)\,v)\,\bar{g}(q^{-2} v',-\sgn(\theta)\,q^2v)\,h'(q^{-2}
v',-\sgn(\theta)\,v)
\\ & &  -   \sgn(\theta)\,f(-\sgn(\theta)\,q^2 v,q^{-2} v')
\bar{g}(-\sgn(\theta)\,v,v')\,h'(-\sgn(\theta)\,v,q^{-2} v')
\end{eqnarray*}
and
\begin{eqnarray*}
C_{0,0}(v) & = & f(v,v) \bar{g}(v,q^2 v)\,h'(v,v) - f(q^2 v,v)\,\bar{g}(v,v)\,h'(v,v)
\\ & & - f(-v,v)\,\bar{g}(-v,q^2v)\,h'(-v,v) + f(-q^2 v,v)\,\bar{g}(-v,v)\,h'(-v,v)
\\ & & + f(-v,-v)\,\bar{g}(-v,-q^2 v)\,h'(-v,-v) - f(-q^2 v, -v)\,\bar{g}(-v,-v)\,h'(-v,-v)
\\ & & - f(v,-v)\,\bar{g}(v,-q^2v)\,h'(v,-v) + f(q^2 v,-v)\,\bar{g}(v,-v)\,h'(v,-v) \ ,
\end{eqnarray*}
\begin{eqnarray}
C_{0,1}(v) & = & \sum_{y=q^{-2}v}^{v'} f(q^2 v,y)\,(\partial_2 \bar{g})(v,y)\,h'(v,y)\,|y| - \sum_{y=q^{-2} v}^{v'}
f(-q^2 v,y)\,(\partial_2 \bar{g})(-v,y)\,h'(-v,y)\,|y| \ , \label{coass.mot2}
\\ C_{0,2}(v) & = &  \sum_{y=-q^{-2} v}^{-q^2} f(q^2 v,y)\,(\partial_2 \bar{g})(v,y)\,h'(v,y)\,|y|
- \sum_{y=-q^{-2} v}^{-q^2}  f(-q^2 v,y)\,(\partial_2 \bar{g})(-v,y)\,h'(-v,y)\,|y| \ , \nonumber
\\ C_{0,3}(v) & = & -
\sum_{x=q^{-2}v}^{v'} (\partial_1 f)(x,v)\,\bar{g}(x,q^2 v)\,h'(x,v)\,|x| + \sum_{x=q^{-2} v}^{v'} (\partial_1 f)(x,-
v) \,\bar{g}(x,-q^2 v)\,h'(x,-v)\,|x| \ , \nonumber
\\ C_{0,4}(v) & = & - \sum_{x=-q^{-2} v}^{-q^2} (\partial_1 f)(x,v)\,\bar{g}(x,q^2 v)\,h'(x,v)\,|x| + \sum_{x=-q^{-2}
v}^{-q^2} (\partial_1 f)(x,-v) \,\bar{g}(x,-q^2v)\,h'(x,-v)\,|x| \ . \nonumber
\end{eqnarray}
A bookkeeping exercise based on Eq.~(\ref{coass.eq6}) reveals that
\begin{equation}
q^{-1}\,|\theta|^{\frac{1}{2}}\,C(v) = C_{\infty,0}(v) + C_{\infty,1}(v) + \sum_{i=0}^4 C_{0,i}(v) \ .
\label{coass.eq7}
\end{equation}

In order to show all the nitty gritty work involved, let us calculate $C(-v,-q^2;v,v')$ in the case of $\theta < 0$
(keep Fig.~1 in mind). So in Eq.~(\ref{coass.eq6}) we have to set $a=-v$, $b=-q^2$, $c=v$ and $d=v'$.
\begin{trivlist}
\item[\ \,(1)] If $x \in I_{q^2} \cap [-q^2,-v]$, then $f(x,q^{-2} v')= 0$ since $(x,q^{-2}v') \notin K(\theta)$. If $x
\in I_{q^2} \cap [-q^2,-v)$ then also $f(q^2 x,q^{-2} v')= 0$ since $(q^2 x,q^{-2}v') \notin K(\theta)$. But note that
$(-q^2 v,q^{-2} v') \in K(\theta)$. So the first sum of  expression (\ref{coass.eq6}) equals $f(-q^2 v,q^{-2}
v')\,\bar{g}(-v,v')\,h'(-v,q^{-2}v')$.
\smallskip
\item[\ \,(2)] For all  $y \in I_{q^2} \cap [v,v']$, clearly $\bar{g}(-1,q^2 y) = \bar{g}(-1,y) =0$. Thus the second
sum in expression (\ref{coass.eq6}) equals 0.
\smallskip
\item[\ \,(3)] For all $x \in I_{q^2} \cap [-q^2,-q^{-2} v]$, we have that $f(q^2 x,v) - f(x,v) = - |x|\,(\partial_1
f)(x,v)$ and for all $y \in I_{q^2} \cap [v,v']$, we have that $\bar{g}(-v,q^2 y) - \bar{g}(-v,y) =  |y|\,(\partial_2
\bar{g})(-v,y)$.
\end{trivlist}
Putting all these results into Eq.~(\ref{coass.eq6}), we see that
\begin{eqnarray*}
& & \hspace{-5ex}- q^{-1}\,|\theta|^{\frac{1}{2}} \, C(-v,-q^2;v,v')
\\ & & \hspace{-7ex}\spat = f(-q^2 v,q^{-2} v')\,\bar{g}(-v,v')\,h'(-v,q^{-2}v') + f(-v,v)\,\bar{g}(-v,q^2 v)\,h'(-v,v) -  f(-q^2 v,v)\,\bar{g}(-v,v)\,h'(-v,v)
\\ & & \hspace{-5ex}\spat + \sum_{x=-q^{-2} v}^{-q^2} \, (\partial_1 f)(x,v) \, \bar{g}(x,q^2 v)\,h'(x,v)\,|x|
+ \sum_{y=q^{-2} v}^{v'}\, f(-q^2 v,y)\, (\partial_2 \bar{g})(-v,y) \,\,h'(-v,y)\,|y| \ .
\end{eqnarray*}
Similarly, one calculates $C(v,v';v,v')$, $C(-v,-q^2;-v,-q^2)$ and $C(v,v';-v,-q^2)$. Adding these four results
together, one finds Eq.~(\ref{coass.eq7}). The case $\theta < 0$ is treated similarly.

\medskip\smallskip

Now we are going to make specific choices for our functions $f$ and $g$. Therefore take $m,n \in \Z$, $p,t \in I_{q^2}$
such that $\sgn(p) = \sgn(t)=  \sgn(\theta)$ and define $f$ and $g$ such that for $(x,y) \in K(\theta)$,
\begin{equation} \label{coass.eq20}
 f(x,y) = \Psis{-q^2/p}{q^{2(1-m)}\,y/|p|}{\sgn(\theta)\,q^{2(1-m)}}\,\,\Psis{-\sgn(\theta)\,q^{2(1+m)}/y}{q^{2(1+m+k)}x}{\sgn(\theta)\,q^{2(1+k)}
y}
\end{equation}
if $\sgn(\theta)\,q^{-2m} y \in I_{q^2}$ and $f(x,y) = 0$ if $\sgn(\theta)\,q^{-2m} y \not\in I_{q^2}$, and on the
other hand,
\begin{equation} \label{coass.eq21}
g(x,y)  =
\Psis{-q^2/t}{q^{2(1-n)}\,x/|t|}{\sgn(\theta)\,q^{2(1-n)}}\,\,\Psis{-\sgn(\theta)\,q^{2(1+n)}/x}{q^{2(1+k+n)}y}{\sgn(\theta)\,q^{2(1+k)}
x}
\end{equation}
if $\sgn(\theta)\,q^{-2n}x \in I_{q^2}$ and $f(x,y) = 0$ if $\sgn(\theta)\,q^{-2n} x \not\in I_{q^2}$. By
Eq.~(\ref{coass.eq16}) and Result \ref{app.res2}, $F^k_{m,p}$,$G^k_{-n,t}$ are proportional to $f$,$g$ respectively. In
the next part we will show that for this choice of $f$ and $g$, $C(v) \rightarrow 0$ as $v \rightarrow 0$.

\smallskip

Notice that $g$ is obtained from $f$ by interchanging $x$ and $y$ and replacing the parameters $m$,$p$ by $n$,$t$.
There is also some symmetry with respect to interchanging $x$ and $y$ in the formulas defining $C_{\infty,0}(v)$ and
$C_{\infty,1}(v)$. Similarly, note the symmetry with respect to this variable interchange between $C_{0,1}(v)$,
$C_{0,2}(v)$ and $C_{0,3}(v)$, $C_{0,4}(v)$ respectively. This is of course in correspondence with the symmetry
observed in Eqs.~(\ref{coass.eq18}) and (\ref{coass.eq19}). We will make use of this symmetry to reduce our work.

\medskip

Let $a \in \C$. Then $(\,(aq^2z;q^2)_\infty - (a z;q^2)_\infty\,)/z =  a \, (a q^2 z;q^2)_\infty$ for all $z \in \C$.
Using this simple fact, a careful inspection learns that for $(x,y) \in K(\theta)$,
\begin{equation} \label{coass.eq22}
 (\partial_1 f)(x,y) =  q^{2(1+m+k)}\,\,
\Psis{-q^2/p}{q^{2(1-m)}\,y/|p|}{\sgn(\theta)\,q^{2(1-m)}}\,\,\Psis{-\sgn(\theta)\,q^{2(1+m)}/y}{q^{2(2+m+k)}x}{\sgn(\theta)\,q^{2(2+k)}
y} \nonumber
\end{equation}
if $\sgn(\theta)\,q^{-2m} y \in I_{q^2}$ and  $(\partial_1 f)(x,y) = 0$ if $\sgn(\theta)\,q^{-2m} y \not\in I_{q^2}$.
On the other hand,
\begin{equation} \label{coass.eq23}
(\partial_2 g)(x,y)  =  q^{2(1+n+k)}\,\,
\Psis{-q^2/t}{q^{2(1-n)}\,x/|t|}{\sgn(\theta)\,q^{2(1-n)}}\,\,\Psis{-\sgn(\theta)\,q^{2(1+n)}/x}{q^{2(2+n+k)}y}{\sgn(\theta)\,q^{2(2+k)}
x}
\end{equation}
if $\sgn(\theta)\,q^{-2n} x \in I_{q^2}$ and $(\partial_2 g)(x,y) = 0$ if $\sgn(\theta)\,q^{-2n} x \not\in I_{q^2}$.

\medskip

Now we will show that the different summands of $C(v)$ as described in $(\ref{coass.eq7})$ converge to 0 as
$v\rightarrow 0$.

\smallskip

\begin{trivlist}
\item[\ \,(1)] First we quickly check that $C_{0,0}(v) \rightarrow 0$ as $v \rightarrow 0$.

\smallskip

Let $(\,(x_r,y_r)\,)_{r=1}^\infty$ be a sequence in $K(\theta)$ that converges to $(0,0)$. Then $\sgn(\theta)\,q^{-2m}
y_r \in I_{q^2}$ and $\sgn(\theta)\,q^{-2n} x_r \in I_{q^2}$ for $r$ big enough. Therefore Eqs.~(\ref{coass.eq20}) and
(\ref{coass.eq21}) imply, as in the proof of Lemma \ref{app.lem2}, that
\begin{eqnarray*}
& & f(x_r,y_r) \rightarrow \Psi(-q^2/p;0;\sgn(\theta) q^{2(1-m)})\,\phi(0;-q^{2(2+m+k)}) \spat \text{ as } r
\rightarrow \infty\ ,
\\ & & g(x_r,y_r) \rightarrow \Psi(-q^2/t;0;\sgn(\theta) q^{2(1-n)})\,\phi(0;-q^{2(2+n+k)}) \spat \text{ as } r
\rightarrow \infty\ ,
\end{eqnarray*}
where we used notation (\ref{coass.eq32}). We have also that $h'(x_r,y_r) \rightarrow 1$ as $r \rightarrow \infty$.
From this all, we easily conclude that $C_{0,0}(v) \rightarrow 0$ as $v \rightarrow 0$.

\medskip

\item[\ \,(2)] Let us now deal with $C_{\infty,1}(v)$. First consider the case where $\theta < 0$. Assume for the moment that
$v \leq q^{-2m}/|\theta|$ and $v \leq q^{-2n}/|\theta|$. Then $q^{-2m} (q^{-2} v') \geq 1$ and $q^{-2n} (q^{-2} v')
\geq 1$. Hence, $\sgn(\theta)\,q^{-2m}\, (q^{-2} v')$ and  $\sgn(\theta)\,q^{-2n}\,(q^{-2}v')$ do not belong to
$I_{q^2}$. As a consequence, $(\partial_1 f)(x,q^{-2}v') = (\partial_2 g)(q^{-2}v',x) = 0$ for all $x \in I_{q^2} \cap
[-q^2,-v]$. We conclude that $C_{\infty,1}(v) = 0$ if $v \leq q^{-2m}/|\theta|$ and $v \leq q^{-2n}/|\theta|$. Thus
$C_{\infty,1}(v) \rightarrow 0$ as $v \rightarrow 0$.

\medskip

Now we look at the more challenging case where $\theta > 0$. Referring to Definition \ref{comult.def1},
Eq.~(\ref{coass.eq22}) and Result \ref{app.res2}, we get the existence of a constant $D > 0$, only depending on $m$,
$p$ and $\theta$ such that
\begin{eqnarray*}
& & |\,(\partial_1 f)(x,q^{-2}y)|^2 \,\,h'(x,q^{-2} y) \,x\,y
\\ & & \spat = D\,\,|a_{\kappa^{-1}(p)}(\kappa^{-1}(q^{-2m})
\kappa^{-1}(q^{-2} y), \kappa^{-1}(q^{-2} y))|^2\,\,|a_{\kappa^{-1}(q^{-2(m+1)}y)}(\kappa^{-1}(\theta
y)\,\kappa^{-1}(x), \kappa^{-1}(x))|^2
\end{eqnarray*}
for all $x,y \in I_{q^2}^+$. Therefore Proposition \ref{comult.ortho1} implies that
\begin{eqnarray*}
& & \sum_{y \in I_{q^2}^+} \sum_{x \in I_{q^2}^+} \, |\,(\partial_1 f)(x,q^{-2}y)|^2 \,\,h'(x,q^{-2} y) \,x\,y
  = D\sum_{y \in I_{q^2}^+} \,|a_{\kappa^{-1}(p)}(\kappa^{-1}(q^{-2m}) \kappa^{-1}(q^{-2} y),
\kappa^{-1}(q^{-2} y))|^2\,
\\ & & \spat\spat \times \ \sum_{x \in I_{q^2}^+}\,|a_{\kappa^{-1}(q^{-2(m+1)}y)}(\kappa^{-1}(\theta y)\,\kappa^{-1}(x), \kappa^{-1}(x))|^2
\\ & & \spat \leq D\sum_{y \in I_{q^2}^+} \,|a_{\kappa^{-1}(p)}(\kappa^{-1}(q^{-2m}) \kappa^{-1}(q^{-2} y), \kappa^{-1}(q^{-2}
y))|^2\ \leq \ D \ .
\end{eqnarray*}
Thus, $\sum_{y \in I_{q^2}^+} y\, \sum_{x \in I_{q^2}^+} \, |(\partial_1 f)(x,q^{-2}y)|^2 \,\,h'(x,q^{-2} y) \,x <
\infty$, which clearly implies that \newline $\sum_{x \in I_{q^2}^+} \, |\,(\partial_1 f)(x,q^{-2}y)|^2 \,\,h'(x,q^{-2}
y) \,x \rightarrow 0$ as $y \rightarrow \infty$.

\smallskip

On the other hand,
\begin{eqnarray*}
& & \sum_{\scriptsize \begin{array}{c} y \in I_{q^2}^+ \\ y \geq 1 \end{array}}\ \sum_{x \in I_{q^2}^+} |g(x,y)|^2\,
h'(x,q^{-2}y)\,x
 = \sum_{\scriptsize \begin{array}{c} y \in I_{q^2}^+ \\ y \geq 1 \end{array}}\ \sum_{x \in I_{q^2}^+} |g(x,y)|^2\,
\frac{(-x,-y;q^2)_\infty}{(-\theta x y;q^2)_\infty}\,(1+q^{-2} y)\,x
\\ & & \spat\spat\spat \leq  (1+q^{-2})\,\sum_{\scriptsize \begin{array}{c} y \in I_{q^2}^+ \\ y \geq 1 \end{array}}\ \sum_{x \in I_{q^2}^+} |g(x,y)|^2\,
\frac{(-x,-y;q^2)_\infty}{(-\theta x y;q^2)_\infty}\,x\,y \ < \infty \ ,
\end{eqnarray*}
since $g$ belongs to $L^2(K(\theta), \Om_{\theta})$. So we also get that $\sum_{x \in I_{q^2}^+} \, |g(x,y)|^2
\,\,h'(x,q^{-2}y) \,x \rightarrow 0$ as $y \rightarrow \infty$.

By the Cauchy-Schwarz inequality, we have
\begin{eqnarray*}
& & \bigl|\,\sum_{x=v}^{v'} (\partial_1 f)(x,q^{-2} v')\,\bar{g}(x,v')\,h'(x,q^{-2}v')\,x \,\bigr|^2
\\ & & \spat \leq \bigl(\,\sum_{x \in I_{q^2}^+} |\,(\partial_1 f)(x,q^{-2}v')|^2\,h'(x,q^{-2}v')\,x \,\bigr)\,
\bigl(\,\sum_{x \in I_{q^2}^+} |g(x,v')|^2\,h'(x,q^{-2} v')\,x \,\bigr) \ .
\end{eqnarray*}
Therefore the considerations above imply that
$$\sum_{x=v}^{v'} (\partial_1 f)(x,q^{-2} v')\,\bar{g}(x,v')\,h'(x,q^{-2} v')\,x \, \rightarrow 0 \hspace{1cm}
\text{as }\ v \rightarrow 0\ .$$  Hence, the obvious symmetry between $f$ and $g$ guarantees that also
$$\sum_{y=v}^{v'} f(v',y)\,(\partial_2 \bar{g})(q^{-2} v',y)\,h'(q^{-2} v',y)\,y\, \rightarrow 0  \hspace{1cm}
\text{as }\ v \rightarrow 0\ .$$ Thus, we conclude that $C_{\infty,1}(v) \rightarrow 0$ as $v \rightarrow 0$.

\medskip

\item[\ \,(3)] In the next step we prove that $C_{0,1}(v) + C_{0,2}(v)  \rightarrow 0$ as $v \rightarrow 0$.
\inlabel{coass.motivation}

\smallskip

Let $v \in I$ such that $v \leq  q^{2(n+1)}$. Then $|q^{-2n} v| \leq q^2$, thus $\pm \sgn(\theta)\,q^{-2n} v \in
I_{q^2}$. Therefore Eqs.~(\ref{coass.eq20}), (\ref{coass.eq23}) and Result \ref{app.res2} imply for $y \in I_{q^2} \cap
([-q^2,-q^{-2} v] \cup [q^{-2} v,v'])$,
\begin{eqnarray}
& & f(\pm q^2 v,y)\,(\partial_2 \bar{g})(\pm v,y)\,h'(\pm v,y)\,|y|\,
 =  q^{2(1+n+k)}\, (\mp v,-y;q^2)_\infty\,(\mp q^2\theta vy;q^2)_\infty^{-1}\,|y| \nonumber
\\ & & \spat\ \ \,\,\,\Psi(-q^2/p;q^{2(1-m)} y/|p|;\sgn(\theta)q^{2(1-m)})\,
\ \,\,\Psi(\mp 1/v;\sgn(\theta)q^{2(1+k)}y;\pm q^{2(2+m+k)}v)\, \nonumber
\\ & & \spat\ \ \,\,\,\Psi(-q^2/t;\pm q^{2(1-n)} v/|t|;\sgn(\theta)\,
q^{2(1-n)})\,\ \,\,\Psi(\mp\sgn(\theta)q^{2(1+n)}/v;q^{2(2+n+k)} y;\pm\sgn(\theta)q^{2(2+k)} v)\, \nonumber
\\ & & \label{alt.eq4}
\end{eqnarray}
if $\sgn(\theta) q^{-2m} y \in I_{q^2}$ and $f(\pm q^2 v,y)\,(\partial_2 \bar{g})(\pm v,y)\,h'(\pm v,y)\,|y|\, = 0$ if
$\sgn(\theta) q^{-2m} y \not\in I_{q^2}$ (since $f(\pm q^2 v,y)$ $= 0$ in this case).

\medskip

If $v \in I$  and  $y \in I_{q^2} \cap ([-q^2,-q^{-2} v] \cup [q^{-2} v,v'])$, then $(\pm v,y) \in K(\theta)$, thus
$\pm q^2\theta vy \geq -q^2$, so  $(\mp q^2\theta vy;q^2)_\infty \geq (q^2;q^2)_\infty$. So we get easily the existence
of a constant $C > 0$ such that
$$|q^{2(1+n+k)}\, (\mp v;q^2)_\infty \,(\mp q^2\theta vy;q^2)_\infty^{-1}\,\Psi(-q^2/t;\pm q^{2(1-n)}
v/|t|;\sgn(\theta)\, q^{2(1-n)})| \leq C$$ for  all  $v \in I$ such that $v \leq
 q^{2(n+1)}$ and  $y \in I_{q^2} \cap ([-q^2,-q^{-2} v] \cup [q^{-2} v,v'])$.

\smallskip

Moreover, Lemma \ref{alt.lem1} implies the existence of
\begin{trivlist}
\item[\ \,(a)] a number $D > 0$ such that $|\Psi(\mp 1/v;\sgn(\theta)q^{2(1+k)}y;\pm q^{2(2+m+k)}v)| \leq D$ for all $y \in
I_{q^2}$ such that $\sgn(\theta)\,q^{-2m}\,y \in I_{q^2}$ and $v \in I$ such that $v \leq
 q^{2(n+1)}$.
\smallskip
\item[\ \,(b)] a function $H_1 \in l^2(I_{q^2})^+$ such that
$|y|^{\frac{1}{2}}\,|(-y;q^2)_\infty|^{\frac{1}{2}}\,|\Psi(-q^2/p;q^{2(1-m)} y/|p|;\sgn(\theta)q^{2(1-m)})| \leq
H_1(y)$ for all $y \in I_{q^2}$
\smallskip
\item[\ \,(c)] a function $H_2 \in l^2(I_{q^2})^+$ such that
$|y|^{\frac{1}{2}}\,|(-y;q^2)_\infty|^{\frac{1}{2}}\,|\Psi(\mp\sgn(\theta)q^{2(1+n)}/v;q^{2(2+n+k)} y;$ \newline $
\pm\sgn(\theta)q^{2(2+k)} v)|$ $\leq H_2(y)$ for all $y \in I_{q^2}$ and $v \in I$ such that $v \leq
 q^{2(n+1)}$.
\end{trivlist}

\smallskip

Now define $H \in l^1(I_{q^2})^+$ as $H = C\,D\,H_1\, H_2$. Then Eq.~(\ref{alt.eq4}) implies that $$|\,f(\pm q^2
v,y)\,(\partial_2 \bar{g})(\pm v,y)\,h'(\pm v,y)\,|y|\,| \leq H(y)$$ for $v \in I$ such that $v \leq q^{2(n+1)}$ and $y
\in I_{q^2} \cap ([-q^2,-q^{-2} v] \cup [q^{-2} v,v'])$.

\smallskip

Define the function $G : I_{q^2} \rightarrow \C$ such that for $y \in I_{q^2}$,
\begin{eqnarray*}
G(y) & = & q^{2(1+n+k)}\,(-y;q^2)_\infty\,|y|\,\,\Psi(-q^2/t;0;\sgn(\theta)\,q^{2(1-n)})\,\,
\Psi(-q^2/p;q^{2(1-m)}y/|p|;\sgn(\theta)\,q^{2(1-m)})
\\ & & \,\phi(\sgn(\theta)\,q^{2(1+k)}y;-q^{2(2+m+k)})\,\,\phi(q^{2(2+n+k)}y;-q^{2(3+n+k)})
\end{eqnarray*}
if $\sgn(\theta)\,q^{-2m} y \in I_{q^2}$ and $G(y) = 0$ if $\sgn(\theta)\,q^{-2m} y \not\in I_{q^2}$ (here we use the
notation introduced in (\ref{coass.eq32})\,).

\smallskip

If we fix $y \in I_{q^2}$, we see that $y \in [-q^2,-q^{-2} v] \cup [q^{-2} v, v']$ if $v$ is small enough. By
Eq.~(\ref{alt.eq4}), we have moreover that  $f(\pm q^2 v,y)\,(\partial_2 \bar{g})(\pm v,y)\,h'(\pm v,y)\,|y|
\rightarrow G(y)$ as $v \rightarrow 0$. Therefore the dominated convergence theorem implies that $G \in l^1(I_{q^2})$
and
$$
\sum_{y=q^{-2}v}^{v'} \, f(\pm q^2 v,y)\,(\partial_2 \bar{g})(\pm v,y)\,h'(\pm v,y)\,|y| \ + \ \sum_{y=-q^2}^{-q^{-2}v}
\, f(\pm q^2 v,y)\,(\partial_2 \bar{g})(\pm v,y)\,h'(\pm v,y)\,|y| \rightarrow \sum_{y \in I_{q^2}}\, G(y)
$$
as $v \rightarrow 0$. Thus, it follows that $C_{0,1}(v) + C_{0,2}(v)  \rightarrow 0$ as $v \rightarrow 0$. The
aforementioned symmetry between $f$ and $g$ then guarantees that also $C_{0,3}(v) + C_{0,4}(v) \rightarrow 0$ as $v
\rightarrow 0$.

\smallskip

\item[\ \,(4)] In the last step we look at $C_{\infty,0}(v)$. First assume that $\theta < 0$. If $v \in I$ and $v \leq
q^{-2m}/|\theta|$ then $\sgn(\theta)\,q^{-2m}\, (q^{-2} v')  = - q^{-2m}/v |\theta| \leq -1$. Thus
Eq.~(\ref{coass.eq20}) implies that $f(-\sgn(\theta)\,q^2 v , q^{-2} v')$ \, $\bar{g}(-\sgn(\theta)\, v,v') \,
h'(-\sgn(\theta)\,v, q^{-2} v') = 0$. So in this case, $f(-\sgn(\theta)\,q^2 v , q^{-2} v') \, \bar{g}(-\sgn(\theta)\,
v,v')$ \newline $h'(-\sgn(\theta)\,v, q^{-2} v') \rightarrow 0$ as $v \rightarrow 0$.

\smallskip

Now we look at the more challenging case where $\theta > 0$. If $v \in I$ and $v \leq q^{2n+2}$, then $\sgn(\theta)
q^{-2m}\,(q^{-2} v')$ and $-q^{-2n} v$ clearly belong to $I_{q^2}$; thus Eqs.~(\ref{coass.eq20}) and (\ref{coass.eq21})
and Result \ref{app.res2} imply that
\begin{eqnarray*}
& & |f(-\sgn(\theta)\,q^2 v,q^{-2} v')\, \bar{g}(-\sgn(\theta)\,v,v')\,h'(-\sgn(\theta)\,v,q^{-2} v')| = (-v,-1/v
\theta;q^2)_\infty
\\ & & \spat\ \ \,\,\,(q^2;q^2)_\infty^{-1}\,\,|\,\Psi(-q^2/p;q^{2(1-m)}/v \theta p;q^{2(1-m)})\,|
\ \,|\,\Psi(1/v;q^{2(1+k)}/ v \theta;-q^{2(2+m+k)}v)\,|
\\ & & \spat\ \ \,\,|\,\Psi(-q^2/t;- q^{2(1-n)} v/t;
q^{2(1-n)})\,|\ \,|\,\Psi(q^{2(1+n)}/v;q^{2(2+n+k)}/v \theta;-q^{2(1+k)} v)\,| \ .
\end{eqnarray*}
By Lemma \ref{coass.lem4}, there exists a number $C > 0$ such that $|\,\Psi(-q^2/p;q^{2(1-m)} / v \theta
p;q^{2(1-m)})\,| \leq C$ for all $v \in I$. Therefore the above equality and Lemma \ref{alt.lem2} allow us to conclude
that
$$f(-\sgn(\theta)\,q^2 v,q^{-2} v')\, \bar{g}(-\sgn(\theta)\,v,v')\,h'(-\sgn(\theta)\,v,q^{-2} v')
\rightarrow 0 \spat \text{as }\ v \rightarrow 0 $$ also in this case. The symmetry between $f$ and $g$ thus guarantees
that also the first term of $C_{\infty,0}(v) \rightarrow 0$ as $v \rightarrow 0$. Therefore $C_{\infty,0}(v)
\rightarrow 0$ as $v \rightarrow 0$.
\end{trivlist}

\medskip

Together with Eq.~(\ref{coass.eq7}) the convergence results proven in parts (1),(2),(3) and (4) imply that $C(v)
\rightarrow 0$ as $v \rightarrow 0$. By (\ref{coass.eq14}), this implies that $\langle \gamma_l^\theta f , g \rangle -
\langle f , (\gamma^\theta_r)^* g \rangle = 0$. So, by our choice of $f$ and $g$, we find that $\langle \gamma_l^\theta
F^k_{m,p} , G^k_{-n,t} \rangle = \langle F^k_{m,p} , (\gamma^\theta_r)^* G^k_{-n,t} \rangle$.

Therefore the conditions stated in (\ref{coass.eq15}) imply that $\langle \gamma_l^\theta x' , y' \rangle = \langle x'
, (\gamma^\theta_r)^* y' \rangle $ for all $x' \in D(\gamma_l^\theta)$ and $y' \in D((\gamma_l^\theta)^*)$.
Consequently, $(\gamma_r^\theta)^* \subseteq (\gamma_l^\theta)^*$\ (*). Since $\gamma_l^\theta$ and $\gamma_r^\theta$
are normal, this implies that $D(\gamma_r^\theta) = D((\gamma_r^\theta)^*) \subseteq D((\gamma_l^\theta)^*) =
D(\gamma_l^\theta)$. But taking the adjoint of the inclusion (*), we arrive at the inclusion $\gamma_l^\theta \subseteq
\gamma_r^\theta$. Hence, $\gamma_l^\theta = \gamma_r^\theta$.
\end{proof}

\medskip

Since $U\,((\de \ot \io)\de(\gamma))\,U^* = M_\ze \ot \gamma_l$ and $U\,((\io \ot \de)\de(\gamma))\,U^* = M_\ze \ot
\gamma_r$ (see the remarks before lemma \ref{coass.lem3}) the previous proposition entails that $(\de \ot
\io)\de(\gamma) = (\io \ot \de)\de(\gamma)$ and we have proven Proposition \ref{coass.prop5}.

\medskip\medskip

We can finally explain why the comultiplication $\de$ is defined in such a way that $U_\theta^* \de(\gamma^*
\gamma)\!\!\restriction_{\ell_\theta'}\! U_\theta = 1 \ot L^{\sgn(\theta)}_\theta$ (see the discussions surrounding
Eq.~(\ref{comult.eq1}) and  after Proposition \ref{comult.prop1}). Recall that this choice is reflected in the presence
of the factor $s(x,y)$ in Definition \ref{comult.def1}.

In the proof above we showed that $C_{0,1}(v) \rightarrow 0$ as $v \rightarrow 0$, cf.~(\ref{coass.motivation}).
Looking at the defining formula (\ref{coass.mot2}) for $C_{0,1}(v)$, it is not too hard to imagine that  for this to be
true, we need at least that the functions $f$ and $g$ defined in Eqs.~(\ref{coass.eq20}) and (\ref{coass.eq21}) have
the same limit behavior when crossing the $Y$-axis, going from the first to the the second quadrant (for these
particular functions, these are just continuous transitions). Note that the relevance of $f$ and $g$ to the
coassociativity of $\de$ is obvious from Lemma \ref{coass.lem1} (and that all potential discontinuous transitions are
filtered out by $U$).

If one changes the comultiplication by leaving out the factor $s(x,y)$ in Definition \ref{comult.def1}, and thus opting
for the equality $U_\theta^* \de(\gamma^* \gamma)\!\!\restriction_{\ell_\theta'}\! U_\theta = 1 \ot L^1_\theta$, it
turns out that $f$ and $g$ have different limit transitions implying that $C_{0,1}(v)$ does not converge to $0$ as $v
\rightarrow 0$. From this it would follow that $(\de \ot \io)\de(\gamma)$ and $(\io \ot \de)\de(\gamma)$ have different
domains in this case.

\medskip\smallskip

From Lemma \ref{coass.lem1} we have two orthonormal bases for the same space $L^2(I_{q^2}^3, \Omega)$; $\bigl(
\tilde{F}_{m,p}^k \bigr)_{p\in I_{q^2},m,k\in \Z}$ and $\bigl( \tilde{G}_{m,p}^k \bigr)_{p\in I_{q^2},m,k\in \Z}$.
Hence, there exists a unitary operator $\cR$ on $L^2(I_{q^2}^3, \Omega)$ mapping the first basis on the second. Its
matrix elements are given by
$$
\cR(m_1,k_1,p_1; m_2,k_2,p_2) = \langle \tilde F_{m_1,p_1}^{k_1}, \tilde G_{m_2,p_2}^{k_2} \rangle_{L^2(I_{q^2}^3,
\Omega)}
$$
and these matrix elements of $\cR$ can be thought of as Racah coefficients. It would be desirable to find an explicit
expression for the Racah coefficients and to show that $\cR(m_1,k_1,p_1; m_2,k_2,p_2)=0$ for $p_1 \not= p_2$, since
this would immediately imply $\gamma_l=\gamma_r$ and hence yield an alternative proof for Proposition \ref{coass.prop3}
and hence of Proposition \ref{coass.prop5}. However, we have not been able to carry out this programme.

\bigskip\medskip

\sectie{Appendix}

\bigskip

In this appendix we collect the basic properties of the $\Psi$-functions defined in Eq.~(1) of the Introduction. Most
of them stem from special function theory, so the proofs in this section are mainly intended for people that are not
too familiar with this theory. As a general reference for $q$-hypergeometric functions we use \cite{Gas}. Let us fix a
number $0 < q < 1$.

\smallskip

Note first of all that the presence of the factor $q^{k^2}$ in the series (1) implies that if we keep two of the
parameters $a$,$b$,$z$ fixed, the series converges uniformly on compact subsets in the remaining parameter. As a
consequence, $\Psi(a;b;q,z)$ is analytic in the remaining parameter. Also, the function $\C^3 \rightarrow \C : (a,b,z)
\mapsto \Psi(a;b;q,z)$ is analytic.

The $\Psi$-functions are closely related to the $_1\vfi_1$-functions (see \cite[Def.~(1.2.22)]{Gas}) in the sense that
$\Psi(a;b;q,z)$ $= (b;q)_\infty\,\,\,_1\vfi_1(a;b;q,z)$ if $b \not\in q^{-\NN}$ (if $b \in q^{-\NN}$,
$_1\vfi_1(a;b;q,z)$ is not defined in general).  As a consequence, a lot of the identities that are known for
$_1\vfi_1$- functions extend by analytic continuation in the parameters to similar identities involving
$\Psi$-functions.

\medskip

A very basic but useful formula is known as the

\begin{result}[$\theta$-product identity] \label{app.res1}
Consider $a \in \C \setminus \{0\}$, $k \in \Z$. Then
\begin{equation} \label{theta}
 (a q^k;q)_\infty\,(q^{1-k}/a;q)_\infty = (-a)^{-k}
\, q^{-\frac{1}{2}k(k-1)}\,(a;q)_\infty\,(q/a;q)_\infty\ .
\end{equation}
\end{result}
\begin{proof}
If $a \in q^\Z$, both sides of the above equation are easily seen to be equal to 0. Now we look at the case where $a
\not\in q^\Z$. Suppose that $k \geq 0$. Then we have for $n \in \NE$,
$$
 \frac{(a q^k ;q)_n}{(a;q)_{n+k}} \, \frac{(\,q^{1-k}/a;q)_{n+k}}{(\,q/a;q)_n} = \frac{\prod_{i=0}^{k-1} (1 -
q^{-i}/a) }{\prod_{i=0}^{k-1} (1-a q^i)} = \frac{(-a)^{-k} \, \prod_{i=0}^{k-1} q^{-i}\, \prod_{i=0}^{k-1} (1-a
q^i)}{\prod_{i=0}^{k-1} (1-a q^i)} = (-a)^{-k}\, q^{-\frac{1}{2} k(k-1)}
$$
so if we let $n$ tend to infinity, Eq.~(\ref{theta}) follows. If $k \leq 0$, we apply Eq.~(\ref{theta}), where we
replace $k$ by $-k$ and $a$ by $q/a$.
\end{proof}

\smallskip

The most important low order $q$-hypergeometric functions are the $_2\vfi_1$ functions. If $a,b,c,z \in \C$ satisfy $c
\not\in q^{-\NN}$ and $|z| < 1$, these are defined as (see \cite[Def.~(1.2.22)]{Gas} and the discussion thereafter)
$$_2\vfi_1(a,b;c;q,z) = \,\vfisss{a}{b}{c}{q,z} = \sum_{n=0}^\infty \frac{(a;q)_n\,(b;q)_n}{(c;q)_n (q;q)_n}\,z^n\ . $$
If we keep $a$,$b$ and $c$ fixed, the function $z, |z| < 1  \mapsto \, _2\vfi_1(a,b;c;q,z)$ has an analytic extension
to the set $\C \setminus [1,\infty)$ (see the beginning of \cite[Sec.~4.3]{Gas}). Of course, the value of this
extension at $z \in \C \setminus [1,\infty)$ is also denoted by $_2\vfi_1(a,b;c;q,z)$. We will not use this fact
directly in this paper but point to this extension because it is used in \cite{Cic}.

\medskip

A connection between $_2\vfi_1$-functions and $_1\vfi_1$-functions is obtained by the following basic limit transition.

\begin{lemma} \label{app.lem2}
Consider $a,b,c,z \in \C$ such that $c \not\in q^{-\NN}$ and $(x_i)_{i=1}^\infty$ a sequence in $\C \setminus \{0\}$
 that converges to $0$ and satisfies $|x_i z| < 1$ for all $i \in \NE$. Then $\bigl(\,_2\vfi_1(a,b/x_i;c;q,x_i z)\,
\bigr)_{i=1}^\infty$ converges to $_1\vfi_1(a;c;q,bz)$.
\end{lemma}
\begin{proof}
By definition, we have that $_2\vfi_1(a,b/x_i;c;q,x_i z) = \sum_{n=0}^\infty \frac{(a;q)_n\,}{(c;q)_n (q;q)_n}\,z^n\,
(b/x_i;q)_n\,x_i^n$ and \newline $(b/x_i;q)_n\,x_i^n = \prod_{k=0}^{n-1} (x_i - bq^k)$. It follows that the $n$-th term
of this series converges to the number $\frac{(a;q)_n\,}{(c;q)_n (q;q)_n}\,(-1)^n q^{\frac{1}{2}n(n-1)}\,(bz)^n$ as $i
\rightarrow \infty$. Take $\vep > 0$ such that $|\vep z| < 1$. So there exists $i_0 \in \NE$ such that $|x_i| < \vep/2$
for all $i \in \Z_{\geq i_0}$. It is not so difficult to see that  there exists $M > 0$ such that $|(b/x_i;q)_n\,x_i^n
\vep^{-n}| \leq M$ for all $n \in \NN$ and $i \in \Z_{\geq i_0}$. Thus, since $\sum_{n=0}^\infty
\left|\,\frac{(a;q)_n\,}{(c;q)_n (q;q)_n}\,(\vep z)^n\,\right| < \infty$, the dominated convergence theorem implies
that the sequence $\bigl(\,_2\vfi_1(a,b/x_i;c;q,x_i z)\,\bigr)_{i=1}^\infty$ converges to
$\sum_{n=0}^\infty\,\frac{(a;q)_n\,}{(c;q)_n (q;q)_n}\,(-1)^n q^{\frac{1}{2}n(n-1)}\,(bz)^n =$ $_1\vfi_1(a;c;q,bz)$.
\end{proof}

\smallskip

This limit transition is used to get a more direct relationship between $_2\vfi_1$-functions and $\Psi$-functions.

\begin{lemma} \label{app.lem1}
Let $a,b,z \in \C$ such that $|z| < 1$. Then $(z;q)_\infty\,_2\vfi_1(a,b;0;q,z) = \Psi(a;az;q,bz)$.
\end{lemma}
\begin{proof}
First suppose that $a \not= 0$ and $az \not\in q^{-\NN}$. Take a sequence $(x_i)_{i=1}^\infty$ in $\C \setminus (\{0\}
\cup q^{-\NN})$ such that $(x_i)_{i=1}^\infty \rightarrow 0$ and $|x_i/a| < 1$ for all $i \in \NE$.

Then Heine's transformation formula  \cite[Eq.~(1.4.5)]{Gas} implies that $_2\vfi_1(a,b;x_i;q,z) = (x_i/a,az;q)_\infty$
\newline $(x_i,z;q)_\infty^{-1}\, _2\vfi_1(abz/x_i,a;az;q,x_i/a)$. If we let $i$ tend to $\infty$ in this equality, we obtain
by the previous lemma the equality $_2\vfi_1(a,b;0;q,z) = (az;q)_\infty\,(z;q)_\infty^{-1}\, _1\vfi_1(a;az;q,bz) =
(z;q)_\infty^{-1}\, \Psi(a;az;q,bz)$. The general result follows by analytic continuation.
\end{proof}

\smallskip

The following formula will be used throughout the paper.

\begin{result} \label{app.res2}
Consider $a,b,z \in \C$ such that $b \not= 0$. Then $\Psi(a;b;q,z) = \Psi(az/b;z;q,b)$.
\end{result}
\begin{proof}
First assume that $b$ and $z$ do not belong to $q^{-\NN}$. Take a sequence $(x_i)_{i=1}^\infty$ in $\C \setminus (\{0\}
\cup q^{-\NN})$ such that $(x_i)_{i=1}^\infty \rightarrow 0$, $|z x_i| < 1$ and $|b x_i| < 1$ for all $i \in \NE$. Let
$i \in \NE$ and apply Heine's transformation formula  \cite[Eq.~(1.4.5)]{Gas} for $a \rightsquigarrow a$, $b
\rightsquigarrow 1/x_i$, $c \rightsquigarrow b$ and $z \rightsquigarrow x_i z$. Thus, $_2\vfi_1(a,1/x_i;b;q,x_i z) = (b
x_i,z;q)_\infty\,\,(b,x_i z;q)_\infty^{-1}\,\, _2\vfi_1(az/b,1/x_i;z;q,b x_i)$. If we let $i$ tend to $\infty$ in this
equality, Lemma \ref{app.lem2} implies that $_1\vfi_1(a;b;q,z) = (z;q)_\infty\,(b;q)_\infty^{-1}$
$_1\vfi_1(az/b;z;q,b)$. Thus, $\Psi(a;b;q,z) = \Psi(az/b;z;q,b)$, still under the assumption that $b,z \not\in
q^{-\NN}$. The general formula now follows from analytic continuation.
\end{proof}

\medskip

Roughly speaking, $q$-contiguous relations are relations between expressions of the form $\Psi(ra;sb;q,tz)$ where
$r,s,t \in \{1,q,q^{-1}\}$. We need the following two.

\begin{lemma} \label{app.contiguous} Consider $a,b,z \in \C$, then
\begin{eqnarray}
\Psi(a;b;q,z) & = &  (1-a)\,\, \Psi(qa;b;q,z) + a \,\,\Psi(a;b;q,qz)\ . \label{app.eq1}
\\ \Psi(a;b;q,z)  & = & (a-b)\,\,\Psi(a;qb;q,qz) + (1-a)\,\, \Psi(qa;qb;q,z)\,\ .  \label{app.eq2}
\end{eqnarray}
\end{lemma}
\begin{proof}
By definition, we have
\begin{eqnarray*}
& & (1-a)\,\, \Psi(qa;b;q,z) =  \sum_{n=0}^\infty (1-a)\,\frac{(qa;q)_n\,(q^n b;q)_\infty}{(q;q)_n}
\,(-1)^n\,q^{\frac{1}{2}n(n-1)}\,z^n
\\ & & \spat  =  \sum_{n=0}^\infty \frac{(a;q)_n\,(1-a q^n)\,(q^n b;q)_\infty}{(q;q)_n} \,(-1)^n\,q^{\frac{1}{2}n(n-1)}\,z^n
\\ & & \spat  =  \sum_{n=0}^\infty \frac{(a;q)_n\,(q^n b;q)_\infty}{(q;q)_n} \,(-1)^n\,q^{\frac{1}{2}n(n-1)}\,z^n
- a\ \sum_{n=0}^\infty \frac{(a;q)_n\,(q^n b;q)_\infty}{(q;q)_n} \,(-1)^n\,q^{\frac{1}{2}n(n-1)}\,(q z)^n
\\ & & \spat = \Psi(a;b;q,z) - a \,\,\Psi(a;b;q,qz) \ .
\end{eqnarray*}
On the other hand,
\begin{eqnarray*}
& &  (a-b)\,\,\Psi(a;qb;q,qz) + (1-a)\,\, \Psi(qa;qb;q,z)
\\ & & \spat\spat =  \sum_{n=0}^\infty \,(a-b)\,\frac{(a;q)_n\,(q^{n+1} b;q)_\infty}{(q;q)_n} \,(-1)^n\,q^{\frac{1}{2}n(n-1)}\,(qz)^n
\\ & & \spat\spat\ \ \  + \sum_{n=0}^\infty\,(1-a)\, \frac{(qa;q)_n\,(q^{n+1} b;q)_\infty}{(q;q)_n} \,(-1)^n\,q^{\frac{1}{2}n(n-1)}\,z^n
\\ & & \spat\spat =  \sum_{n=0}^\infty \,(a-b)q^n\,\frac{(a;q)_n\,(q^{n+1} b;q)_\infty}{(q;q)_n} \,(-1)^n\,q^{\frac{1}{2}n(n-1)}\,z^n
\\ & & \spat\spat\ \ \  + \sum_{n=0}^\infty\,(1-q^n a)\, \frac{(a;q)_n\,(q^{n+1} b;q)_\infty}{(q;q)_n}
\,(-1)^n\,q^{\frac{1}{2}n(n-1)}\,z^n
\\ & & \spat\spat =  \sum_{n=0}^\infty \,(1-q^n b)\,\frac{(a;q)_n\,(q^{n+1} b;q)_\infty}{(q;q)_n} \,(-1)^n\,q^{\frac{1}{2}n(n-1)}\,z^n
\\ & & \spat\spat =  \sum_{n=0}^\infty \,\frac{(a;q)_n\,(q^n b;q)_\infty}{(q;q)_n} \,(-1)^n\,q^{\frac{1}{2}n(n-1)}\,z^n
 = \Psi(a;b;q,z) \ .
\end{eqnarray*}
\end{proof}

We need the following transformation formulas.

\smallskip

\begin{proposition} \label{app.prop1}
Consider $a,b,z \in \C$ and $k \in \Z$. Then
\begin{trivlist}
\item[\ \,(1)] $(q^{k+1}/a;q)_\infty\,\Psi(aq^{-k};q^{-k+1};q,z) = (q/a;q)_\infty\, (a z/q)^k\,\,
\Psi(a;q^{k+1};q,zq^k)$ \ \ if $a,z \not= 0$,
\smallskip
\item[\ \,(2)] $(q^{k+1} a/z;q)_\infty \, \Psi(q^{-k};a;q,z) = (z/q)^k\,(q^k a;q)_\infty\,\Psi(q^{-k};qa/z;q,q^2/z)$ \ \ if
$z \not=0$ and $k \geq 0$.
\end{trivlist}
\end{proposition}
\vspace{0.5ex}
\begin{trivlist}
\item[\ \,{\it Proof.} (1)] Suppose first that $k \geq 0$. Note that $(q^{n+1-k};q)_\infty = 0$ for all $n \in \Z_{\leq k-1}$. Hence,
\begin{eqnarray*}
& & \hspace{-4ex} (q^{k+1}/a;q)_\infty\,\Psi(aq^{-k};q^{-k+1};q,z) = (q^{1+k}/a;q)_\infty\,\sum_{n=k}^\infty
\frac{(aq^{-k};q)_n (q^{n+1-k};q)_\infty}{(q;q)_n} \, (-1)^n\, q^{\frac{1}{2}n(n-1)}\,z^n
\\ & & \spat \hspace{-5ex} = (q^{k+1}/a;q)_\infty\,\sum_{n=0}^\infty \frac{(aq^{-k};q)_{n+k}
(q^{n+1};q)_\infty}{(q;q)_{n+k}} \, (-1)^{n+k}\, q^{\frac{1}{2}(n+k)(n+k-1)}\,z^{n+k}
\\ & & \spat \hspace{-5ex} = (-1)^k\,q^{\frac{1}{2}k(k-1)}\,\, z^k\,(q^{k+1}/a;q)_\infty\,(a/q;q^{-1})_k\,\sum_{n=0}^\infty \frac{(a;q)_n
(q^{n+k+1};q)_\infty}{(q;q)_n} \, (-1)^n  \,q^{nk}\, q^{\frac{1}{2}n(n-1)}\,z^n
\\ & & \spat \hspace{-5ex} = (-1)^k\,q^{\frac{1}{2}k(k-1)}\,\, z^k\,(q^{k+1}/a;q)_\infty\,(q/a;q)_k\,(-a)^k\, q^{-k}\,q^{-\frac{1}{2}k(k-1)}  \,\, \Psi(a;q^{k+1};q,zq^k)
\\ & & \spat \hspace{-5ex} = (q/a;q)_\infty\, (a z/q)^k\,\, \Psi(a;q^{k+1};q,zq^k) \ .
\end{eqnarray*}
If $k < 0$, we apply (1) where we replace $k$ by $-k$, $a$ by $a q^{-k}$ and $z$ by $z q^k$ to obtain (1) in this case.
\smallskip
\item[\ \,(2)] First we assume that $a \not= 0$ and $|q^{k+1} a/z| < 1$. Since $(q^{-k};q)_n = 0$ if $n \geq k+1$, the
series terminates and
\begin{eqnarray}
& & \hspace{-7ex} (q^{k+1} a/z;q)_\infty \, \Psi(q^{-k};a;q,z) = (q^{k+1} a/z;q)_\infty \,
\sum_{n=0}^k\,\,\frac{(q^{-k};q)_n\, (q^n a;q)_\infty}{(q;q)_n}\,(-1)^n\,q^{\frac{1}{2}n(n-1)}\, z^n \nonumber
\\ & & \hspace{-7ex} \spat = (q^{k+1} a/z;q)_\infty \, \sum_{n=0}^k\,\,\frac{(q^{-k};q)_{k-n}\, (q^{k-n}
a;q)_\infty}{(q;q)_{k-n}}\,(-1)^{k-n}\,q^{\frac{1}{2}(k-n)(k-n-1)}\, z^{k-n} \nonumber
\\ & & \hspace{-7ex} \spat = (q^{k+1} a/z;q)_\infty \,(-1)^k\,q^{\frac{1}{2}k(k-1)}\,z^k\, \sum_{n=0}^k\,\,\frac{(q^{-k};q)_{k-n}\, (q^{k-n}
a;q)_\infty}{(q;q)_{k-n}}\,(-1)^n\,q^{-kn}\,q^{\frac{1}{2}n(n+1)}\, z^{-n} \ . \label{app.eq3}
\end{eqnarray}
Now,
\begin{eqnarray*}
(q^{k-n} a;q)_\infty & = & (q^k a;q)_\infty\,(q^{k-1} a;q^{-1})_n = (q^k a;q)_\infty\,\,q^{nk}
\,q^{-\frac{1}{2}n(n+1)}\,(-a)^n \,(q^{1-k}/a;q)_n
\\ & &
\\ \frac{(q^{-k};q)_{k-n}}{(q;q)_{k-n}} & = & \frac{\prod_{i=-k}^{-1}(1-q^i)}{\prod_{i=-n}^{-1} (1-q^i)}\,\frac{\prod_{i=k-n+1}^k (1-q^i)}{\prod_{i=1}^k
(1-q^i)} = (-1)^k\,q^{-\frac{1}{2} k(k+1)} \,
\frac{(-1)^n\,q^{kn}\,q^{-\frac{1}{2}n(n-1)}\,(q^{-k};q)_n}{(-1)^n\,q^{-\frac{1}{2}n(n+1)}\,(q;q)_n}
\\ & = & (-1)^k\,q^{-\frac{1}{2} k(k+1)} \, q^{(k+1)n}\,
\frac{(q^{-k};q)_n}{(q;q)_n} \,
\end{eqnarray*}
which upon substitution in Eq.~(\ref{app.eq3}) implies that
\begin{eqnarray*}
& & \hspace{-1ex} (q^{k+1} a/z;q)_\infty \, \Psi(q^{-k};a;q,z)  = (q^{k+1} a/z;q)_\infty \,(q^k a;q)_\infty \,(z/q)^k\,
\sum_{n=0}^k\,\,\frac{(q^{-k};q)_n\, (q^{1-k}/a;q)_n}{(q;q)_n}\, (q^{k+1} a/z)^n
\\ & & \spat = (q^{k+1} a/z;q)_\infty \,(q^k a;q)_\infty \,(z/q)^k\,_2\vfi_1(q^{-k},q^{1-k}/a;0;q,q^{k+1}a/z)
\\ & & \spat = (q^k a;q)_\infty \,(z/q)^k\,\,\,\Psi(q^{-k};qa/z;q,q^2/z) \ ,
\end{eqnarray*}
where we used Lemma \ref{app.lem1} in the last step. The general result now follows by analytic continuation. \qed
\end{trivlist}

\bigskip

In most cases it is not necessary to know the precise expression for the $\Psi$-functions, sometimes only the
asymptotic behavior matters. We need the following estimate (see  \cite[Lem.~2.8]{ErikJasper2}).

\smallskip

\begin{lemma} \label{coass.lem4}
Consider $a,z \in \C$, $k \in \Z_{\geq 0}$. Then $|\Psi(a;q^{1-k};q,z)| \leq (-|a|;q)_\infty\,(-|z|;q)_\infty \,
|z|^k\,q^{\frac{1}{2}k(k-1)}$.
\end{lemma}
\begin{proof}
For $n \in \NN$, we have $|(a;q)_n| \leq (-|a|;q)_n \leq (-|a|;q)_\infty$. Also notice that $(q^{n-k+1};q)_\infty = 0$
for all $n \in \NN$ satisfying $n < k$. Hence, by the definition of $\Psi$, we see that
\begin{eqnarray*}
& & \left|\,\Psi(a;q^{1-k};q,z)\,\right|   \leq  (-|a|;q)_\infty \, \sum_{n=0}^\infty
\frac{(q^{n-k+1};q)_\infty}{(q;q)_n}\,q^{\frac{1}{2}n(n-1)}\,|z|^n
\\ & & \spat  =   (-|a|;q)_\infty \, \sum_{n=k}^\infty
\frac{(q^{n-k+1};q)_\infty}{(q;q)_n}\,q^{\frac{1}{2}n(n-1)}\,|z|^n  =   (-|a|;q)_\infty \, \sum_{n=0}^\infty
\frac{(q^{n+1};q)_\infty}{(q;q)_{n+k}}\,q^{\frac{1}{2}(n+k)(n+k-1)}\,|z|^{n+k}
\\ & & \spat =  (-|a|;q)_\infty \,|z|^{k}\, \sum_{n=0}^\infty
\frac{(q^{n+k+1};q)_\infty}{(q;q)_n}\,q^{\frac{1}{2}n(n-1)}\,q^{nk}\,q^{\frac{1}{2}k(k-1)}\,|z|^n
\\ & & \spat \leq
(-|a|;q)_\infty \,|z|^{k}\,q^{\frac{1}{2}k(k-1)}\, \sum_{n=0}^\infty
\frac{1}{(q;q)_n}\,q^{\frac{1}{2}n(n-1)}\,(q^{k}\,|z|)^n
 =  (-|a|;q)_\infty \, (-|z| q^{k};q)_\infty \,|z|^{k}\,q^{\frac{1}{2}k(k-1)}
\end{eqnarray*}
where we used  \cite[Eq.~1.3.16]{Gas} in the last equality. Now the lemma follows.
\end{proof}

\medskip

The last two estimates play a vital role in the proof of the coassociativity of the comultiplication.

\begin{lemma} \label{alt.lem1}
Consider $a,b,c \in \C$, $p \in \qz$ and $\be > 0$. Let $m \in [1,\infty]$,  $M >0$,  $\vep
> 0$ and define the $q$-interval $J$ as $J = \{\,x \in \qz \mid |x| \leq M \text{ or } \sgn(x) =
\sgn(p)\,\}$. Then there exist $N
>0$ and $f \in l^m(J)^+$ such that
$$
| \Psi(a/\lambda;p x;q,b\lambda) | \leq N \hspace{6ex} \text{ and } \hspace{6ex}
|x|^{\be}\,\,|(cx;q)_\infty|^{\frac{1}{2}} \,\, | \Psi(a/\lambda;p x;q,b\lambda) | \leq f(x)
$$
for all $x \in J$ and $\lambda \in \C \setminus \{0\}$ such that $|\lambda| \leq \vep$.
\end{lemma}
\begin{proof} Choose $k \in \Z$ such that $|p| = q^k$.  Take $r,s > 0$ such that $r > |b|(\vep+|a|)$ and $s > |c|$.
Let $\lambda \in \C \setminus \{0\}$ such that $|\lambda| \leq \vep$ and $x \in J$. Since $(p x q^r ;q)_\infty \geq 0$
for all $r \in \Z$, we get that
\begin{eqnarray}
& &  |\,\Psi(a/\lambda;px;q,b\lambda)\,|   \leq \sum_{n=0}^\infty \frac{|(a/\lambda;q)_n|\, (p x
q^n;q)_\infty}{(q;q)_n}\, q^{\frac{1}{2}n(n-1)}\,|b\lambda|^n \nonumber
\\ & & \spat  =  \sum_{n=0}^\infty \frac{(p x
q^n;q)_\infty}{(q;q)_n}\,(\,\prod_{i=0}^{n-1} |\lambda-q^i a|\,) q^{\frac{1}{2}n(n-1)}\,|b|^n  \leq \sum_{n=0}^\infty
\frac{(p x q^n;q)_\infty}{(q;q)_n}\, q^{\frac{1}{2}n(n-1)}\,(|b|\,(|\lambda|+|a|))^n \nonumber
\\ & & \spat  \leq   \Psi(0;p
x;q,-r)\ . \label{app.eq24}
\end{eqnarray}
This implies the existence of $C > 0$ such that $|\,\Psi(a/\lambda;px;q,b\lambda)\,| \leq C$  and
$|(cx;q)_\infty|^{\frac{1}{2}} \leq C$ for all $x \in J$ such that $|x| \leq \max\{M,q^2/|p|\}$ and $\lambda \in \C
\setminus \{0\}$ such that $|\lambda| \leq \vep$. \!\!\!\inlabel{alt.eq1}

\smallskip

Define the function $f : J \rightarrow \R^+$ such that for $x \in J$, we have that $f(x) = C^2\,|x|^{\be}$ if $|x| \leq
\max\{M,q^2/|p|\}$ and
$$f(x) = (-r;q)_\infty\,(-s;q)_\infty^{\frac{1}{2}}\,(-q/s;q)_\infty^{\frac{1}{2}} \, q^{\frac{1}{2}k(k-1)}\,q^{(k+\be)t}\,
r^{1-k-t}\,s^{-\frac{t}{2}}\,q^{\frac{1}{4}t(t-1)}$$ if $|x| > \max\{M,q^2/|p|\}$ and $t \in \Z$ is such that $|x| =
q^t$. Then it is clear that $f \in l^m(J)^+$.

\smallskip

Take $x \in J$ such that $|x| > \max\{M,q^2/|p|\}$. Choose $t \in \Z$ such that $|x| = q^t$. Since $x \in J$ and $|x|
> M$, we have that $\sgn(x) = \sgn(p)$, hence $px = q^{k+t}$. It is also clear that $t \leq 1-k$, thus $1-k-t \geq 0$. By estimate (\ref{app.eq24}) and Lemma \ref{coass.lem4},  we get for
all $\lambda \in \C \setminus \{0\}$ such that $|\lambda| \leq \vep$,
\begin{eqnarray}
 |\,\Psi(a/\lambda;px;q,b\lambda)\,| & \leq & \Psi(0;p x;q,-r) = \Psi(0;q^{1-(1-k-t)};-r) \leq
(-r;q)_\infty\,r^{1-k-t}\,q^{\frac{1}{2}(k+t-1)(k+t)} \nonumber
\\ & = & (-r;q)_\infty\,r^{1-k-t}\,q^{\frac{1}{2}k(k-1)}\,q^{kt}\, q^{\frac{1}{2}t(t-1)} \ .\label{alt.eq2}
\end{eqnarray}
Also, Lemma \ref{app.res1} guarantees that
\begin{eqnarray*}
 |\,(cx;q)_\infty\,|^{\frac{1}{2}} & \leq & (-|cx| ;q)^{\frac{1}{2}} \leq (-s q^t;q)_\infty^{\frac{1}{2}} =
\frac{(-s;q)_\infty^{\frac{1}{2}}\,(-q/s;q)_\infty^{\frac{1}{2}}}{(-q^{1-t}/s;q)_\infty^{\frac{1}{2}}} \,
s^{-\frac{t}{2}} \,q^{-\frac{1}{4}t(t-1)}
\\ & \leq & (-s;q)_\infty^{\frac{1}{2}}\,(-q/s;q)_\infty^{\frac{1}{2}} \,
s^{-\frac{t}{2}} \,q^{-\frac{1}{4}t(t-1)} \ .
\end{eqnarray*}
This estimate, together with the estimates in (\ref{alt.eq1}) and (\ref{alt.eq2}) imply that
$|x|^{\be}\,\,|(cx;q)_\infty|^{\frac{1}{2}} \,\, | \Psi(a/\lambda;p x;q,b\lambda) |$ $ \leq f(x)$ for all $x \in J$ and
$\lambda \in \C \setminus \{0\}$ such that $|\lambda| \leq \vep$.

Estimate (\ref{alt.eq2}) also implies the existence of $D > 0$ such that $| \Psi(a/\lambda;p x;q,b\lambda) | \leq D$
for all $x \in J$ such that $|x| > \max\{M,q^2/|p|\}$ and $\lambda \in \C \setminus \{0\}$ such that $|\lambda| \leq
\vep$. So if we set $N = \max\{C,D\}$, we find by the estimate in (\ref{alt.eq1}) that $| \Psi(a/\lambda;p
x;q,b\lambda) | \leq N$ for all $x \in J$ and $\lambda \in \C \setminus \{0\}$ such that $|\lambda| \leq \vep$.
\end{proof}

\medskip

The next result is an easy consequence of the previous lemma.

\smallskip

\begin{lemma} \label{alt.lem2}
Consider $a,b,c \in \C$, $p \in \qz$. Let $(x_i)_{i \in I}$ be a net in $\qz$ such that $\sgn(x_i) = \sgn(p)$ for all
$i \in I$ and $(x_i)_{i \in I} \rightarrow \infty$. Then the net $\bigl(\, |(cx_i;q)_\infty|^{\frac{1}{2}} \,\, |\Psi(a
x_i;p x_i;q,b/ x_i)|\bigr)_{i \in I}$ converges to $0$.
\end{lemma}

\end{document}